%
%
 
\font \aufont=  cmr10 at 12pt  
\font\titfont=  cmr10 at 18pt 
\font\footfont=cmr10 at 8  pt
   

 \magnification=1200

%
%

\def\fmtversion{2.0}
\catcode`\@=11
\ifx\amstexloaded@\relax\catcode`\@=\active
 \endinput\else\let\amstexloaded@\relax\fi
\def\W@{\immediate\write\sixt@@n}
\def\CR@{\W@{}\W@{AmS-TeX - Version \fmtversion}\W@{}
\W@{COPYRIGHT 1985, 1990 - AMERICAN MATHEMATICAL SOCIETY}
\W@{Use of this macro package is not restricted provided}
\W@{each use is acknowledged upon publication.}\W@{}}
\CR@
\everyjob{\CR@}
\toksdef\toks@@=2
\long\def\rightappend@#1\to#2{\toks@{\\{#1}}\toks@@
 =\expandafter{#2}\xdef#2{\the\toks@@\the\toks@}\toks@{}\toks@@{}}
\def\alloclist@{}
\newif\ifalloc@
\def\showallocations{{\def\\{\immediate\write\m@ne}\alloclist@}\alloc@true}
\def\alloc@#1#2#3#4#5{\global\advance\count1#1by\@ne
 \ch@ck#1#4#2\allocationnumber=\count1#1
 \global#3#5=\allocationnumber
 \edef\next@{\string#5=\string#2\the\allocationnumber}%
 \expandafter\rightappend@\next@\to\alloclist@}
\newcount\count@@
\newcount\count@@@
\def\FN@{\futurelet\next}
\def\DN@{\def\next@}
\def\DNii@{\def\nextii@}
\def\RIfM@{\relax\ifmmode}
\def\RIfMIfI@{\relax\ifmmode\ifinner}
\def\setboxz@h{\setbox\z@\hbox}
\def\wdz@{\wd\z@}
\def\boxz@{\box\z@}
\def\setbox@ne{\setbox\@ne}
\def\wd@ne{\wd\@ne}
\def\iterate{\body\expandafter\iterate\else\fi}
\newlinechar=`\^^J
\def\err@#1{\errmessage{AmS-TeX error: #1}}
\newhelp\defaulthelp@{Sorry, I already gave what help I could...^^J
Maybe you should try asking a human?^^J
An error might have occurred before I noticed any problems.^^J
``If all else fails, read the instructions.''}
\def\Err@#1{\errhelp\defaulthelp@\errmessage{AmS-TeX error: #1}}
\def\eat@#1{}
\def\in@#1#2{\def\in@@##1#1##2##3\in@@{\ifx\in@##2\in@false\else\in@true\fi}%
 \in@@#2#1\in@\in@@}
\newif\ifin@
\def\space@.{\futurelet\space@\relax}
\space@. %
\newhelp\athelp@
{Only certain combinations beginning with @ make sense to me.^^J
Perhaps you wanted \string\@\space for a printed @?^^J
I've ignored the character or group after @.}
\def\futureletnextat@{\futurelet\next\at@}
{\catcode`\@=\active
\lccode`\Z=`\@ \lccode`\I=`\I \lowercase{%
\gdef@{\csname futureletnextatZ\endcsname}
\expandafter\gdef\csname atZ\endcsname                                      
 {\ifcat\noexpand\next a\def\next{\csname atZZ\endcsname}\else
 \ifcat\noexpand\next0\def\next{\csname atZZ\endcsname}\else
 \ifcat\noexpand\next\relax\def\next{\csname atZZZ\endcsname}\else
 \def\next{\csname atZZZZ\endcsname}\fi\fi\fi\next}
\expandafter\gdef\csname atZZ\endcsname#1{\expandafter                      
 \ifx\csname #1Zat\endcsname\relax\def\next
 {\errhelp\expandafter=\csname athelpZ\endcsname
 \csname errZ\endcsname{Invalid use of \string@}}\else
 \def\next{\csname #1Zat\endcsname}\fi\next}
\expandafter\gdef\csname atZZZ\endcsname#1{\expandafter                     
 \ifx\csname \string#1ZZat\endcsname\relax\def\next
 {\errhelp\expandafter=\csname athelpZ\endcsname
 \csname errZ\endcsname{Invalid use of \string@}}\else
 \def\next{\csname \string#1ZZat\endcsname}\fi\next}
\expandafter\gdef\csname atZZZZ\endcsname#1{\errhelp                        
 \expandafter=\csname athelpZ\endcsname
 \csname errZ\endcsname{Invalid use of \string@}}}}                         
\def\atdef@#1{\expandafter\def\csname #1@at\endcsname}
\def\atdef@@#1{\expandafter\def\csname \string#1@@at\endcsname}
\newhelp\defahelp@{If you typed \string\define\space cs instead of
\string\define\string\cs\space^^J
I've substituted an inaccessible control sequence so that your^^J
definition will be completed without mixing me up too badly.^^J
If you typed \string\define{\string\cs} the inaccessible control sequence^^J
was defined to be \string\cs, and the rest of your^^J
definition appears as input.}
\newhelp\defbhelp@{I've ignored your definition, because it might^^J
conflict with other uses that are important to me.}
\def\define{\FN@\define@}
\def\define@{\ifcat\noexpand\next\relax
 \expandafter\define@@\else\errhelp\defahelp@                               
 \err@{\string\define\space must be followed by a control
 sequence}\expandafter\def\expandafter\nextii@\fi}                          
\def\undefined@@@@@@@@@@{}
\def\preloaded@@@@@@@@@@{}
\def\next@@@@@@@@@@{}
\def\define@@#1{\ifx#1\relax\errhelp\defbhelp@                              
 \err@{\string#1\space is already defined}\DN@{\DNii@}\else
 \expandafter\ifx\csname\expandafter\eat@\string                            
 #1@@@@@@@@@@\endcsname\undefined@@@@@@@@@@\errhelp\defbhelp@
 \err@{\string#1\space can't be defined}\DN@{\DNii@}\else
 \expandafter\ifx\csname\expandafter\eat@\string#1\endcsname\relax          
 \global\let#1\undefined\DN@{\def#1}\else\errhelp\defbhelp@
 \err@{\string#1\space is already defined}\DN@{\DNii@}\fi
 \fi\fi\next@}
\let\redefine\def
\def\predefine#1#2{\let#1#2}
\def\undefine#1{\let#1\undefined}
\newdimen\captionwidth@
\captionwidth@\hsize
\advance\captionwidth@-1.5in
\def\pagewidth#1{\hsize#1\relax
 \captionwidth@\hsize\advance\captionwidth@-1.5in}
\def\pageheight#1{\vsize#1\relax}
\def\hcorrection#1{\advance\hoffset#1\relax}
\def\vcorrection#1{\advance\voffset#1\relax}

\let\graveaccent\`
\let\acuteaccent\'
\let\tildeaccent\~
\let\hataccent\^
\let\underscore\_
\let\B\=
\let\D\.
\let\ic@\/
\def\/{\unskip\ic@}
\def\textfonti{\the\textfont\@ne}
\def\t#1#2{{\edef\next@{\the\font}\textfonti\accent"7F \next@#1#2}}
\def~{\unskip\nobreak\ \ignorespaces}
\def\.{.\spacefactor\@m}
\atdef@;{\leavevmode\null;}
\atdef@:{\leavevmode\null:}
\atdef@?{\leavevmode\null?}
\def\@{\char64 }
\def\relaxnext@{\let\next\relax}
\atdef@-{\relaxnext@\leavevmode
 \DN@{\ifx\next-\DN@-{\FN@\nextii@}\else
  \DN@{\leavevmode\hbox{-}}\fi\next@}%
 \DNii@{\ifx\next-\DN@-{\leavevmode\hbox{---}}\else
  \DN@{\leavevmode\hbox{--}}\fi\next@}%
 \FN@\next@}
\def\srdr@{\kern.16667em}
\def\drsr@{\kern.02778em}
\def\sldl@{\kern.02778em}
\def\dlsl@{\kern.16667em}
\atdef@"{\unskip\relaxnext@
 \DN@{\ifx\next\space@\DN@. {\FN@\nextii@}\else
  \DN@.{\FN@\nextii@}\fi\next@.}%
 \DNii@{\ifx\next`\DN@`{\FN@\nextiii@}\else
  \ifx\next\lq\DN@\lq{\FN@\nextiii@}\else
  \DN@####1{\FN@\nextiv@}\fi\fi\next@}%
 \def\nextiii@{\ifx\next`\DN@`{\sldl@``}\else\ifx\next\lq
  \DN@\lq{\sldl@``}\else\DN@{\dlsl@`}\fi\fi\next@}%
 \def\nextiv@{\ifx\next'\DN@'{\srdr@''}\else
  \ifx\next\rq\DN@\rq{\srdr@''}\else\DN@{\drsr@'}\fi\fi\next@}%
 \FN@\next@}

\def\textfontii{\the\textfont\tw@}
\def\lbrace@{\delimiter"4266308 }
\def\rbrace@{\delimiter"5267309 }
\def\{{\RIfM@\lbrace@\else{\textfontii f}\spacefactor\@m\fi}
\def\}{\RIfM@\rbrace@\else
 \let\@sf\empty\ifhmode\edef\@sf{\spacefactor\the\spacefactor}\fi
 {\textfontii g}\@sf\relax\fi}
\let\lbrace\{
\let\rbrace\}
\def\AmSTeX{{\textfontii A}\kern-.1667em\lower.5ex\hbox
 {\textfontii M}\kern-.125em{\textfontii S}-\TeX}
\def\vmodeerr@#1{\Err@{\string#1\space not allowed between paragraphs}}
\def\mathmodeerr@#1{\Err@{\string#1\space not allowed in math mode}}
\def\linebreak{\RIfM@\mathmodeerr@\linebreak\else
 \ifhmode\unskip\unkern\break\else\vmodeerr@\linebreak\fi\fi}

\newskip\saveskip@
\def\allowlinebreak{\RIfM@\mathmodeerr@\allowlinebreak\else
 \ifhmode\saveskip@\lastskip\unskip
 \allowbreak\ifdim\saveskip@>\z@\hskip\saveskip@\fi
 \else\vmodeerr@\allowlinebreak\fi\fi}
\def\nolinebreak{\RIfM@\mathmodeerr@\nolinebreak\else
 \ifhmode\saveskip@\lastskip\unskip
 \nobreak\ifdim\saveskip@>\z@\hskip\saveskip@\fi
 \else\vmodeerr@\nolinebreak\fi\fi}
\def\newline{\relaxnext@
 \DN@{\RIfM@\expandafter\mathmodeerr@\expandafter\newline\else
  \ifhmode\ifx\next\par\else
  \expandafter\unskip\expandafter\null\expandafter\hfill\expandafter\break\fi
  \else
  \expandafter\vmodeerr@\expandafter\newline\fi\fi}%
 \FN@\next@}
\def\dmatherr@#1{\Err@{\string#1\space not allowed in display math mode}}
\def\nondmatherr@#1{\Err@{\string#1\space not allowed in non-display math
 mode}}
\def\onlydmatherr@#1{\Err@{\string#1\space allowed only in display math mode}}
\def\nonmatherr@#1{\Err@{\string#1\space allowed only in math mode}}
\def\mathbreak{\RIfMIfI@\break\else
 \dmatherr@\mathbreak\fi\else\nonmatherr@\mathbreak\fi}
\def\nomathbreak{\RIfMIfI@\nobreak\else
 \dmatherr@\nomathbreak\fi\else\nonmatherr@\nomathbreak\fi}
\def\allowmathbreak{\RIfMIfI@\allowbreak\else
 \dmatherr@\allowmathbreak\fi\else\nonmatherr@\allowmathbreak\fi}
\def\pagebreak{\RIfM@
 \ifinner\nondmatherr@\pagebreak\else\postdisplaypenalty-\@M\fi
 \else\ifvmode\removelastskip\break\else\vadjust{\break}\fi\fi}
\def\nopagebreak{\RIfM@
 \ifinner\nondmatherr@\nopagebreak\else\postdisplaypenalty\@M\fi
 \else\ifvmode\nobreak\else\vadjust{\nobreak}\fi\fi}
\def\nonvmodeerr@#1{\Err@{\string#1\space not allowed within a paragraph
 or in math}}
\def\vnonvmode@#1#2{\relaxnext@\DNii@{\ifx\next\par\DN@{#1}\else
 \DN@{#2}\fi\next@}%
 \ifvmode\DN@{#1}\else
 \DN@{\FN@\nextii@}\fi\next@}
\def\newpage{\vnonvmode@{\vfill\break}{\nonvmodeerr@\newpage}}
\def\smallpagebreak{\vnonvmode@\smallbreak{\nonvmodeerr@\smallpagebreak}}
\def\medpagebreak{\vnonvmode@\medbreak{\nonvmodeerr@\medpagebreak}}
\def\bigpagebreak{\vnonvmode@\bigbreak{\nonvmodeerr@\bigpagebreak}}
\def\NoBlackBoxes{\global\overfullrule\z@}
\def\BlackBoxes{\global\overfullrule5\p@}
\def\Invalid@#1{\def#1{\Err@{\Invalid@@\string#1}}}
\def\Invalid@@{Invalid use of }
\Invalid@\caption
\Invalid@\captionwidth
\newdimen\smallcaptionwidth@
\def\topspace{\mid@false\ins@}
\def\midspace{\mid@true\ins@}
\newif\ifmid@
\def\captionfont@{}
\def\ins@#1{\relaxnext@\allowbreak
 \smallcaptionwidth@\captionwidth@\gdef\thespace@{#1}%
 \DN@{\ifx\next\space@\DN@. {\FN@\nextii@}\else
  \DN@.{\FN@\nextii@}\fi\next@.}%
 \DNii@{\ifx\next\caption\DN@\caption{\FN@\nextiii@}%
  \else\let\next@\nextiv@\fi\next@}%
 \def\nextiv@{\vnonvmode@
  {\ifmid@\expandafter\midinsert\else\expandafter\topinsert\fi
   \vbox to\thespace@{}\endinsert}
  {\ifmid@\nonvmodeerr@\midspace\else\nonvmodeerr@\topspace\fi}}%
 \def\nextiii@{\ifx\next\captionwidth\expandafter\nextv@
  \else\expandafter\nextvi@\fi}%
 \def\nextv@\captionwidth##1##2{\smallcaptionwidth@##1\relax\nextvi@{##2}}%
 \def\nextvi@##1{\def\thecaption@{\captionfont@##1}%
  \DN@{\ifx\next\space@\DN@. {\FN@\nextvii@}\else
   \DN@.{\FN@\nextvii@}\fi\next@.}%
  \FN@\next@}%
 \def\nextvii@{\vnonvmode@
  {\ifmid@\expandafter\midinsert\else
  \expandafter\topinsert\fi\vbox to\thespace@{}\nobreak\smallskip
  \setboxz@h{\noindent\ignorespaces\thecaption@\unskip}%
  \ifdim\wdz@>\smallcaptionwidth@\centerline{\vbox{\hsize\smallcaptionwidth@
   \noindent\ignorespaces\thecaption@\unskip}}%
  \else\centerline{\boxz@}\fi\endinsert}
  {\ifmid@\nonvmodeerr@\midspace
  \else\nonvmodeerr@\topspace\fi}}%
 \FN@\next@}
\def\newcodes@{\catcode`\\=12 \catcode`\{=12 \catcode`\}=12 \catcode`\#=12
 \catcode`\%=12\relax}
\def\oldcodes@{\catcode`\\=0 \catcode`\{=1 \catcode`\}=2 \catcode`\#=6
 \catcode`\%=14\relax}
\def\comment{\newcodes@\endlinechar=10 \comment@}
{\lccode`\0=`\\
\lowercase{\gdef\comment@#1^^J{\comment@@#10endcomment\comment@@@}%
\gdef\comment@@#10endcomment{\FN@\comment@@@}%
\gdef\comment@@@#1\comment@@@{\ifx\next\comment@@@\let\next\comment@
 \else\def\next{\oldcodes@\endlinechar=`\^^M\relax}%
 \fi\next}}}
\def\pr@m@s{\ifx'\next\DN@##1{\prim@s}\else\let\next@\egroup\fi\next@}
\def\prime{{\null\prime@\null}}
\mathchardef\prime@="0230
\let\dsize\displaystyle
\let\tsize\textstyle
\let\ssize\scriptstyle

\def\,{\RIfM@\mskip\thinmuskip\relax\else\kern.16667em\fi}
\def\!{\RIfM@\mskip-\thinmuskip\relax\else\kern-.16667em\fi}
\let\thinspace\,
\let\negthinspace\!
\def\medspace{\RIfM@\mskip\medmuskip\relax\else\kern.222222em\fi}
\def\negmedspace{\RIfM@\mskip-\medmuskip\relax\else\kern-.222222em\fi}
\def\thickspace{\RIfM@\mskip\thickmuskip\relax\else\kern.27777em\fi}
\let\;\thickspace
\def\negthickspace{\RIfM@\mskip-\thickmuskip\relax\else
 \kern-.27777em\fi}
\atdef@,{\RIfM@\mskip.1\thinmuskip\else\leavevmode\null,\fi}
\atdef@!{\RIfM@\mskip-.1\thinmuskip\else\leavevmode\null!\fi}
\atdef@.{\RIfM@&&\else\leavevmode.\spacefactor3000 \fi}
\def\and{\DOTSB\;\mathchar"3026 \;}

\def\frac#1#2{{#1\over#2}}

\newdimen\ex@
\ex@.2326ex
\Invalid@\thickness
\def\thickfrac{\relaxnext@
 \DN@{\ifx\next\thickness\let\next@\nextii@\else
 \DN@{\nextii@\thickness1}\fi\next@}%
 \DNii@\thickness##1##2##3{{##2\above##1\ex@##3}}%
 \FN@\next@}

\def\thickfracwithdelims#1#2{\relaxnext@\def\ldelim@{#1}\def\rdelim@{#2}%
 \DN@{\ifx\next\thickness\let\next@\nextii@\else
 \DN@{\nextii@\thickness1}\fi\next@}%
 \DNii@\thickness##1##2##3{{##2\abovewithdelims
 \ldelim@\rdelim@##1\ex@##3}}%
 \FN@\next@}
\def\binom#1#2{{#1\choose#2}}

\def\:{\nobreak\hskip.1111em\mathpunct{}\nonscript\mkern-\thinmuskip{:}\hskip
 .3333emplus.0555em\relax}
\def\snug{\unskip\kern-\mathsurround}
\def\topsmash{\top@true\bot@false\smash@}
\def\botsmash{\top@false\bot@true\smash@}
\newif\iftop@
\newif\ifbot@
\def\smash{\top@true\bot@true\smash@}
\def\smash@{\RIfM@\expandafter\mathpalette\expandafter\mathsm@sh\else
 \expandafter\makesm@sh\fi}
\def\finsm@sh{\iftop@\ht\z@\z@\fi\ifbot@\dp\z@\z@\fi\leavevmode\boxz@}
\def\LimitsOnSums{\global\let\slimits@\displaylimits}
\def\NoLimitsOnSums{\global\let\slimits@\nolimits}
\LimitsOnSums
\mathchardef\coprod@="1360       \def\coprod{\DOTSB\coprod@\slimits@}
\mathchardef\bigvee@="1357       \def\bigvee{\DOTSB\bigvee@\slimits@}
\mathchardef\bigwedge@="1356     \def\bigwedge{\DOTSB\bigwedge@\slimits@}
\mathchardef\biguplus@="1355     \def\biguplus{\DOTSB\biguplus@\slimits@}
\mathchardef\bigcap@="1354       \def\bigcap{\DOTSB\bigcap@\slimits@}
\mathchardef\bigcup@="1353       \def\bigcup{\DOTSB\bigcup@\slimits@}
\mathchardef\prod@="1351         \def\prod{\DOTSB\prod@\slimits@}
\mathchardef\sum@="1350          \def\sum{\DOTSB\sum@\slimits@}
\mathchardef\bigotimes@="134E    \def\bigotimes{\DOTSB\bigotimes@\slimits@}
\mathchardef\bigoplus@="134C     \def\bigoplus{\DOTSB\bigoplus@\slimits@}
\mathchardef\bigodot@="134A      \def\bigodot{\DOTSB\bigodot@\slimits@}
\mathchardef\bigsqcup@="1346     \def\bigsqcup{\DOTSB\bigsqcup@\slimits@}
\def\LimitsOnInts{\global\let\ilimits@\displaylimits}
\def\NoLimitsOnInts{\global\let\ilimits@\nolimits}
\NoLimitsOnInts
\def\int{\DOTSI\intop\ilimits@}
\def\oint{\DOTSI\ointop\ilimits@}
\def\intic@{\mathchoice{\hskip.5em}{\hskip.4em}{\hskip.4em}{\hskip.4em}}
\def\negintic@{\mathchoice
 {\hskip-.5em}{\hskip-.4em}{\hskip-.4em}{\hskip-.4em}}
\def\intkern@{\mathchoice{\!\!\!}{\!\!}{\!\!}{\!\!}}
\def\intdots@{\mathchoice{\plaincdots@}
 {{\cdotp}\mkern1.5mu{\cdotp}\mkern1.5mu{\cdotp}}
 {{\cdotp}\mkern1mu{\cdotp}\mkern1mu{\cdotp}}
 {{\cdotp}\mkern1mu{\cdotp}\mkern1mu{\cdotp}}}
\newcount\intno@
\def\iint{\DOTSI\intno@\tw@\FN@\ints@}
\def\iiint{\DOTSI\intno@\thr@@\FN@\ints@}
\def\iiiint{\DOTSI\intno@4 \FN@\ints@}
\def\idotsint{\DOTSI\intno@\z@\FN@\ints@}
\def\ints@{\findlimits@\ints@@}
\newif\iflimtoken@
\newif\iflimits@
\def\findlimits@{\limtoken@true\ifx\next\limits\limits@true
 \else\ifx\next\nolimits\limits@false\else
 \limtoken@false\ifx\ilimits@\nolimits\limits@false\else
 \ifinner\limits@false\else\limits@true\fi\fi\fi\fi}
\def\multint@{\int\ifnum\intno@=\z@\intdots@                                
 \else\intkern@\fi                                                          
 \ifnum\intno@>\tw@\int\intkern@\fi                                         
 \ifnum\intno@>\thr@@\int\intkern@\fi                                       
 \int}                                                                      
\def\multintlimits@{\intop\ifnum\intno@=\z@\intdots@\else\intkern@\fi
 \ifnum\intno@>\tw@\intop\intkern@\fi
 \ifnum\intno@>\thr@@\intop\intkern@\fi\intop}
\def\ints@@{\iflimtoken@                                                    
 \def\ints@@@{\iflimits@\negintic@\mathop{\intic@\multintlimits@}\limits    
  \else\multint@\nolimits\fi                                                
  \eat@}                                                                    
 \else                                                                      
 \def\ints@@@{\iflimits@\negintic@
  \mathop{\intic@\multintlimits@}\limits\else
  \multint@\nolimits\fi}\fi\ints@@@}
\def\LimitsOnNames{\global\let\nlimits@\displaylimits}
\def\NoLimitsOnNames{\global\let\nlimits@\nolimits@}
\LimitsOnNames
\def\nolimits@{\relaxnext@
 \DN@{\ifx\next\limits\DN@\limits{\nolimits}\else
  \let\next@\nolimits\fi\next@}%
 \FN@\next@}
\def\newmcodes@{\mathcode`\'="0027 \mathcode`\*="002A \mathcode`\.="613A
 \mathcode`\-="002D \mathcode`\/="002F \mathcode`\:="603A }
\def\operatorname#1{\mathop{\newmcodes@\kern\z@\fam\z@#1}\nolimits@}
\def\operatornamewithlimits#1{\mathop{\newmcodes@\kern\z@\fam\z@#1}\nlimits@}
\def\qopname@#1{\mathop{\fam\z@#1}\nolimits@}
\def\qopnamewl@#1{\mathop{\fam\z@#1}\nlimits@}
\def\arccos{\qopname@{arccos}}
\def\arcsin{\qopname@{arcsin}}
\def\arctan{\qopname@{arctan}}
\def\arg{\qopname@{arg}}
\def\cos{\qopname@{cos}}
\def\cosh{\qopname@{cosh}}
\def\cot{\qopname@{cot}}
\def\coth{\qopname@{coth}}
\def\csc{\qopname@{csc}}
\def\deg{\qopname@{deg}}
\def\det{\qopnamewl@{det}}
\def\dim{\qopname@{dim}}
\def\exp{\qopname@{exp}}
\def\gcd{\qopnamewl@{gcd}}
\def\hom{\qopname@{hom}}
\def\inf{\qopnamewl@{inf}}
\def\injlim{\qopnamewl@{inj\,lim}}
\def\ker{\qopname@{ker}}
\def\lg{\qopname@{lg}}
\def\lim{\qopnamewl@{lim}}
\def\liminf{\qopnamewl@{lim\,inf}}
\def\limsup{\qopnamewl@{lim\,sup}}
\def\ln{\qopname@{ln}}
\def\log{\qopname@{log}}
\def\max{\qopnamewl@{max}}
\def\min{\qopnamewl@{min}}
\def\Pr{\qopnamewl@{Pr}}
\def\projlim{\qopnamewl@{proj\,lim}}
\def\sec{\qopname@{sec}}
\def\sin{\qopname@{sin}}
\def\sinh{\qopname@{sinh}}
\def\sup{\qopnamewl@{sup}}
\def\tan{\qopname@{tan}}
\def\tanh{\qopname@{tanh}}
\def\varinjlim{\mathop{\vtop{\ialign{##\crcr
 \hfil\rm lim\hfil\crcr\noalign{\nointerlineskip}\rightarrowfill\crcr
 \noalign{\nointerlineskip\kern-\ex@}\crcr}}}}
\def\varprojlim{\mathop{\vtop{\ialign{##\crcr
 \hfil\rm lim\hfil\crcr\noalign{\nointerlineskip}\leftarrowfill\crcr
 \noalign{\nointerlineskip\kern-\ex@}\crcr}}}}
\def\varliminf{\mathop{\underline{\vrule height\z@ depth.2exwidth\z@
 \hbox{\rm lim}}}}

\newdimen\buffer@
\buffer@\fontdimen13 \tenex
\newdimen\buffer
\buffer\buffer@

\def\ResetBuffer{\fontdimen13 \tenex\buffer@\global\buffer\buffer@}
\def\shave#1{\mathop{\hbox{$\m@th\fontdimen13 \tenex\z@                     
 \displaystyle{#1}$}}\fontdimen13 \tenex\buffer}

\Invalid@\\
\def\Let@{\relax\iffalse{\fi\let\\=\cr\iffalse}\fi}
\Invalid@\vspace
\def\vspace@{\def\vspace##1{\crcr\noalign{\vskip##1\relax}}}
\def\multilimits@{\bgroup\vspace@\Let@
 \baselineskip\fontdimen10 \scriptfont\tw@
 \advance\baselineskip\fontdimen12 \scriptfont\tw@
 \lineskip\thr@@\fontdimen8 \scriptfont\thr@@
 \lineskiplimit\lineskip
 \vbox\bgroup\ialign\bgroup\hfil$\m@th\scriptstyle{##}$\hfil\crcr}
\def\Sb{_\multilimits@}
\def\endSb{\crcr\egroup\egroup\egroup}
\def\Sp{^\multilimits@}

\def\spreadlines#1{\RIfMIfI@\onlydmatherr@\spreadlines\else
 \openup#1\relax\fi\else\onlydmatherr@\spreadlines\fi}
\def\Mathstrut@{\copy\Mathstrutbox@}
\newbox\Mathstrutbox@
\setbox\Mathstrutbox@\null
\setboxz@h{$\m@th($}
\ht\Mathstrutbox@\ht\z@
\dp\Mathstrutbox@\dp\z@
\newdimen\spreadmlines@
\def\spreadmatrixlines#1{\RIfMIfI@
 \onlydmatherr@\spreadmatrixlines\else
 \spreadmlines@#1\relax\fi\else\onlydmatherr@\spreadmatrixlines\fi}
\def\matrix{\null\,\vcenter\bgroup\Let@\vspace@
 \normalbaselines\openup\spreadmlines@\ialign
 \bgroup\hfil$\m@th##$\hfil&&\quad\hfil$\m@th##$\hfil\crcr
 \Mathstrut@\crcr\noalign{\kern-\baselineskip}}
\def\endmatrix{\crcr\Mathstrut@\crcr\noalign{\kern-\baselineskip}\egroup
 \egroup\,}
\def\format{\crcr\egroup\iffalse{\fi\ifnum`}=0 \fi\format@}
\newtoks\hashtoks@
\hashtoks@{#}
\def\format@#1\\{\def\preamble@{#1}%
 \def\l{$\m@th\the\hashtoks@$\hfil}%
 \def\c{\hfil$\m@th\the\hashtoks@$\hfil}%
 \def\r{\hfil$\m@th\the\hashtoks@$}%
 \edef\Preamble@{\preamble@}\ifnum`{=0 \fi\iffalse}\fi
 \ialign\bgroup\span\Preamble@\crcr}
\def\smallmatrix{\null\,\vcenter\bgroup\vspace@\Let@
 \baselineskip9\ex@\lineskip\ex@
 \ialign\bgroup\hfil$\m@th\scriptstyle{##}$\hfil&&\thickspace\hfil
 $\m@th\scriptstyle{##}$\hfil\crcr}
\def\endsmallmatrix{\crcr\egroup\egroup\,}

\newmuskip\dotsspace@
\dotsspace@1.5mu
\def\strip@#1 {#1}
\def\spacehdots#1\for#2{\multispan{#2}\xleaders
 \hbox{$\m@th\mkern\strip@#1 \dotsspace@.\mkern\strip@#1 \dotsspace@$}\hfill}
\def\hdotsfor#1{\spacehdots\@ne\for{#1}}
\def\multispan@#1{\omit\mscount#1\unskip\loop\ifnum\mscount>\@ne\sp@n\repeat}
\def\spaceinnerhdots#1\for#2\after#3{\multispan@{\strip@#2 }#3\xleaders
 \hbox{$\m@th\mkern\strip@#1 \dotsspace@.\mkern\strip@#1 \dotsspace@$}\hfill}
\def\innerhdotsfor#1\after#2{\spaceinnerhdots\@ne\for#1\after{#2}}
\def\cases{\bgroup\spreadmlines@\jot\left\{\,\matrix\format\l&\quad\l\\}
\def\endcases{\endmatrix\right.\egroup}
\newif\ifinany@
\newif\ifinalign@
\newif\ifingather@
\def\strut@{\copy\strutbox@}
\newbox\strutbox@
\setbox\strutbox@\hbox{\vrule height8\p@ depth3\p@ width\z@}
\def\topaligned{\null\,\vtop\aligned@}
\def\botaligned{\null\,\vbox\aligned@}
\def\aligned{\null\,\vcenter\aligned@}
\def\aligned@{\bgroup\vspace@\Let@
 \ifinany@\else\openup\jot\fi\ialign
 \bgroup\hfil\strut@$\m@th\displaystyle{##}$&
 $\m@th\displaystyle{{}##}$\hfil\crcr}
\def\endaligned{\crcr\egroup\egroup}

\def\alignedat#1{\null\,\vcenter\bgroup\doat@{#1}\vspace@\Let@
 \ifinany@\else\openup\jot\fi\ialign\bgroup\span\preamble@@\crcr}
\newcount\atcount@
\def\doat@#1{\toks@{\hfil\strut@$\m@th
 \displaystyle{\the\hashtoks@}$&$\m@th\displaystyle
 {{}\the\hashtoks@}$\hfil}
 \atcount@#1\relax\advance\atcount@\m@ne                                    
 \loop\ifnum\atcount@>\z@\toks@=\expandafter{\the\toks@&\hfil$\m@th
 \displaystyle{\the\hashtoks@}$&$\m@th
 \displaystyle{{}\the\hashtoks@}$\hfil}\advance
  \atcount@\m@ne\repeat                                                     
 \xdef\preamble@{\the\toks@}\xdef\preamble@@{\preamble@}}

\def\gathered{\null\,\vcenter\bgroup\vspace@\Let@
 \ifinany@\else\openup\jot\fi\ialign
 \bgroup\hfil\strut@$\m@th\displaystyle{##}$\hfil\crcr}
\def\endgathered{\crcr\egroup\egroup}
\newif\iftagsleft@
\def\TagsOnLeft{\global\tagsleft@true}
\def\TagsOnRight{\global\tagsleft@false}
\TagsOnLeft
\newif\ifmathtags@
\def\TagsAsMath{\global\mathtags@true}
\def\TagsAsText{\global\mathtags@false}
\TagsAsText
\def\tagform@#1{\hbox{\rm(\ignorespaces#1\unskip)}}
\def\thetag{\leavevmode\tagform@}
\def\tag#1$${\iftagsleft@\leqno\else\eqno\fi                                
 \maketag@#1\maketag@                                                       
 $$}                                                                        
\def\maketag@{\FN@\maketag@@}
\def\maketag@@{\ifx\next"\expandafter\maketag@@@\else\expandafter\maketag@@@@
 \fi}
\def\maketag@@@"#1"#2\maketag@{\hbox{\rm#1}}                                
\def\maketag@@@@#1\maketag@{\ifmathtags@\tagform@{$\m@th#1$}\else
 \tagform@{#1}\fi}
\interdisplaylinepenalty\@M
\def\allowdisplaybreaks{\RIfMIfI@
 \onlydmatherr@\allowdisplaybreaks\else
 \interdisplaylinepenalty\z@\fi\else\onlydmatherr@\allowdisplaybreaks\fi}
\Invalid@\allowdisplaybreak
\Invalid@\displaybreak
\Invalid@\intertext
\def\allowdisplaybreak@{\def\allowdisplaybreak{\crcr\noalign{\allowbreak}}}
\def\displaybreak@{\def\displaybreak{\crcr\noalign{\break}}}
\def\intertext@{\def\intertext##1{\crcr\noalign{\vskip\belowdisplayskip
 \vbox{\normalbaselines\noindent##1}\vskip\abovedisplayskip}}}
\newskip\centering@
\centering@\z@ plus\@m\p@
\def\align{\relax\ifingather@\DN@{\csname align (in
  \string\gather)\endcsname}\else
 \ifmmode\ifinner\DN@{\onlydmatherr@\align}\else
  \let\next@\align@\fi
 \else\DN@{\onlydmatherr@\align}\fi\fi\next@}
\newhelp\andhelp@
{An extra & here is so disastrous that you should probably exit^^J
and fix things up.}
\newif\iftag@
\newcount\and@
\def\align@{\inalign@true\inany@true
 \vspace@\allowdisplaybreak@\displaybreak@\intertext@
 \def\tag{\global\tag@true\ifnum\and@=\z@\DN@{&&}\else
          \DN@{&}\fi\next@}%
 \iftagsleft@\DN@{\csname align \endcsname}\else
  \DN@{\csname align \space\endcsname}\fi\next@}
\def\Tag@{\iftag@\else\errhelp\andhelp@\err@{Extra & on this line}\fi}
\newdimen\lwidth@
\newdimen\rwidth@
\newdimen\maxlwidth@
\newdimen\maxrwidth@
\newdimen\totwidth@
\def\measure@#1\endalign{\lwidth@\z@\rwidth@\z@\maxlwidth@\z@\maxrwidth@\z@
 \global\and@\z@                                                            
 \setbox@ne\vbox                                                            
  {\everycr{\noalign{\global\tag@false\global\and@\z@}}\Let@                
  \halign{\setboxz@h{$\m@th\displaystyle{\@lign##}$}
   \global\lwidth@\wdz@                                                     
   \ifdim\lwidth@>\maxlwidth@\global\maxlwidth@\lwidth@\fi                  
   \global\advance\and@\@ne                                                 
   &\setboxz@h{$\m@th\displaystyle{{}\@lign##}$}\global\rwidth@\wdz@        
   \ifdim\rwidth@>\maxrwidth@\global\maxrwidth@\rwidth@\fi                  
   \global\advance\and@\@ne                                                
   &\Tag@
   \eat@{##}\crcr#1\crcr}}
 \totwidth@\maxlwidth@\advance\totwidth@\maxrwidth@}                       
\def\displ@y@{\global\dt@ptrue\openup\jot
 \everycr{\noalign{\global\tag@false\global\and@\z@\ifdt@p\global\dt@pfalse
 \vskip-\lineskiplimit\vskip\normallineskiplimit\else
 \penalty\interdisplaylinepenalty\fi}}}
\def\black@#1{\noalign{\ifdim#1>\displaywidth
 \dimen@\prevdepth\nointerlineskip                                          
 \vskip-\ht\strutbox@\vskip-\dp\strutbox@                                   
 \vbox{\noindent\hbox to#1{\strut@\hfill}}
 \prevdepth\dimen@                                                          
 \fi}}
\expandafter\def\csname align \space\endcsname#1\endalign
 {\measure@#1\endalign\global\and@\z@                                       
 \ifingather@\everycr{\noalign{\global\and@\z@}}\else\displ@y@\fi           
 \Let@\tabskip\centering@                                                   
 \halign to\displaywidth
  {\hfil\strut@\setboxz@h{$\m@th\displaystyle{\@lign##}$}
  \global\lwidth@\wdz@\boxz@\global\advance\and@\@ne                        
  \tabskip\z@skip                                                           
  &\setboxz@h{$\m@th\displaystyle{{}\@lign##}$}
  \global\rwidth@\wdz@\boxz@\hfill\global\advance\and@\@ne                  
  \tabskip\centering@                                                       
  &\setboxz@h{\@lign\strut@\maketag@##\maketag@}
  \dimen@\displaywidth\advance\dimen@-\totwidth@
  \divide\dimen@\tw@\advance\dimen@\maxrwidth@\advance\dimen@-\rwidth@     
  \ifdim\dimen@<\tw@\wdz@\llap{\vtop{\normalbaselines\null\boxz@}}
  \else\llap{\boxz@}\fi                                                    
  \tabskip\z@skip                                                          
  \crcr#1\crcr                                                             
  \black@\totwidth@}}                                                      
\newdimen\lineht@
\expandafter\def\csname align \endcsname#1\endalign{\measure@#1\endalign
 \global\and@\z@
 \ifdim\totwidth@>\displaywidth\let\displaywidth@\totwidth@\else
  \let\displaywidth@\displaywidth\fi                                        
 \ifingather@\everycr{\noalign{\global\and@\z@}}\else\displ@y@\fi
 \Let@\tabskip\centering@\halign to\displaywidth
  {\hfil\strut@\setboxz@h{$\m@th\displaystyle{\@lign##}$}%
  \global\lwidth@\wdz@\global\lineht@\ht\z@                                 
  \boxz@\global\advance\and@\@ne
  \tabskip\z@skip&\setboxz@h{$\m@th\displaystyle{{}\@lign##}$}%
  \global\rwidth@\wdz@\ifdim\ht\z@>\lineht@\global\lineht@\ht\z@\fi         
  \boxz@\hfil\global\advance\and@\@ne
  \tabskip\centering@&\kern-\displaywidth@                                  
  \setboxz@h{\@lign\strut@\maketag@##\maketag@}%
  \dimen@\displaywidth\advance\dimen@-\totwidth@
  \divide\dimen@\tw@\advance\dimen@\maxlwidth@\advance\dimen@-\lwidth@
  \ifdim\dimen@<\tw@\wdz@
   \rlap{\vbox{\normalbaselines\boxz@\vbox to\lineht@{}}}\else
   \rlap{\boxz@}\fi
  \tabskip\displaywidth@\crcr#1\crcr\black@\totwidth@}}
\expandafter\def\csname align (in \string\gather)\endcsname
  #1\endalign{\vcenter{\align@#1\endalign}}
\Invalid@\endalign
\newif\ifxat@
\def\alignat{\RIfMIfI@\DN@{\onlydmatherr@\alignat}\else
 \DN@{\csname alignat \endcsname}\fi\else
 \DN@{\onlydmatherr@\alignat}\fi\next@}
\newif\ifmeasuring@
\newbox\savealignat@
\expandafter\def\csname alignat \endcsname#1#2\endalignat                   
 {\inany@true\xat@false
 \def\tag{\global\tag@true\count@#1\relax\multiply\count@\tw@
  \xdef\tag@{}\loop\ifnum\count@>\and@\xdef\tag@{&\tag@}\advance\count@\m@ne
  \repeat\tag@}%
 \vspace@\allowdisplaybreak@\displaybreak@\intertext@
 \displ@y@\measuring@true                                                   
 \setbox\savealignat@\hbox{$\m@th\displaystyle\Let@
  \attag@{#1}
  \vbox{\halign{\span\preamble@@\crcr#2\crcr}}$}%
 \measuring@false                                                           
 \Let@\attag@{#1}
 \tabskip\centering@\halign to\displaywidth
  {\span\preamble@@\crcr#2\crcr                                             
  \black@{\wd\savealignat@}}}                                               
\Invalid@\endalignat
\def\xalignat{\RIfMIfI@
 \DN@{\onlydmatherr@\xalignat}\else
 \DN@{\csname xalignat \endcsname}\fi\else
 \DN@{\onlydmatherr@\xalignat}\fi\next@}
\expandafter\def\csname xalignat \endcsname#1#2\endxalignat
 {\inany@true\xat@true
 \def\tag{\global\tag@true\def\tag@{}\count@#1\relax\multiply\count@\tw@
  \loop\ifnum\count@>\and@\xdef\tag@{&\tag@}\advance\count@\m@ne\repeat\tag@}%
 \vspace@\allowdisplaybreak@\displaybreak@\intertext@
 \displ@y@\measuring@true\setbox\savealignat@\hbox{$\m@th\displaystyle\Let@
 \attag@{#1}\vbox{\halign{\span\preamble@@\crcr#2\crcr}}$}%
 \measuring@false\Let@
 \attag@{#1}\tabskip\centering@\halign to\displaywidth
 {\span\preamble@@\crcr#2\crcr\black@{\wd\savealignat@}}}
\def\attag@#1{\let\Maketag@\maketag@\let\TAG@\Tag@                          
 \let\Tag@=0\let\maketag@=0
 \ifmeasuring@\def\llap@##1{\setboxz@h{##1}\hbox to\tw@\wdz@{}}%
  \def\rlap@##1{\setboxz@h{##1}\hbox to\tw@\wdz@{}}\else
  \let\llap@\llap\let\rlap@\rlap\fi                                         
 \toks@{\hfil\strut@$\m@th\displaystyle{\@lign\the\hashtoks@}$\tabskip\z@skip
  \global\advance\and@\@ne&$\m@th\displaystyle{{}\@lign\the\hashtoks@}$\hfil
  \ifxat@\tabskip\centering@\fi\global\advance\and@\@ne}
 \iftagsleft@
  \toks@@{\tabskip\centering@&\Tag@\kern-\displaywidth
   \rlap@{\@lign\maketag@\the\hashtoks@\maketag@}%
   \global\advance\and@\@ne\tabskip\displaywidth}\else
  \toks@@{\tabskip\centering@&\Tag@\llap@{\@lign\maketag@
   \the\hashtoks@\maketag@}\global\advance\and@\@ne\tabskip\z@skip}\fi      
 \atcount@#1\relax\advance\atcount@\m@ne
 \loop\ifnum\atcount@>\z@
 \toks@=\expandafter{\the\toks@&\hfil$\m@th\displaystyle{\@lign
  \the\hashtoks@}$\global\advance\and@\@ne
  \tabskip\z@skip&$\m@th\displaystyle{{}\@lign\the\hashtoks@}$\hfil\ifxat@
  \tabskip\centering@\fi\global\advance\and@\@ne}\advance\atcount@\m@ne
 \repeat                                                                    
 \xdef\preamble@{\the\toks@\the\toks@@}
 \xdef\preamble@@{\preamble@}
 \let\maketag@\Maketag@\let\Tag@\TAG@}                                      
\Invalid@\endxalignat
\def\xxalignat{\RIfMIfI@
 \DN@{\onlydmatherr@\xxalignat}\else\DN@{\csname xxalignat
  \endcsname}\fi\else
 \DN@{\onlydmatherr@\xxalignat}\fi\next@}
\expandafter\def\csname xxalignat \endcsname#1#2\endxxalignat{\inany@true
 \vspace@\allowdisplaybreak@\displaybreak@\intertext@
 \displ@y\setbox\savealignat@\hbox{$\m@th\displaystyle\Let@
 \xxattag@{#1}\vbox{\halign{\span\preamble@@\crcr#2\crcr}}$}%
 \Let@\xxattag@{#1}\tabskip\z@skip\halign to\displaywidth
 {\span\preamble@@\crcr#2\crcr\black@{\wd\savealignat@}}}
\def\xxattag@#1{\toks@{\tabskip\z@skip\hfil\strut@
 $\m@th\displaystyle{\the\hashtoks@}$&%
 $\m@th\displaystyle{{}\the\hashtoks@}$\hfil\tabskip\centering@&}%
 \atcount@#1\relax\advance\atcount@\m@ne\loop\ifnum\atcount@>\z@
 \toks@=\expandafter{\the\toks@&\hfil$\m@th\displaystyle{\the\hashtoks@}$%
  \tabskip\z@skip&$\m@th\displaystyle{{}\the\hashtoks@}$\hfil
  \tabskip\centering@}\advance\atcount@\m@ne\repeat
 \xdef\preamble@{\the\toks@\tabskip\z@skip}\xdef\preamble@@{\preamble@}}
\Invalid@\endxxalignat
\newdimen\gwidth@
\newdimen\gmaxwidth@
\def\gmeasure@#1\endgather{\gwidth@\z@\gmaxwidth@\z@\setbox@ne\vbox{\Let@
 \halign{\setboxz@h{$\m@th\displaystyle{##}$}\global\gwidth@\wdz@
 \ifdim\gwidth@>\gmaxwidth@\global\gmaxwidth@\gwidth@\fi
 &\eat@{##}\crcr#1\crcr}}}
\def\gather{\RIfMIfI@\DN@{\onlydmatherr@\gather}\else
 \ingather@true\inany@true\def\tag{&}%
 \vspace@\allowdisplaybreak@\displaybreak@\intertext@
 \displ@y\Let@
 \iftagsleft@\DN@{\csname gather \endcsname}\else
  \DN@{\csname gather \space\endcsname}\fi\fi
 \else\DN@{\onlydmatherr@\gather}\fi\next@}
\expandafter\def\csname gather \space\endcsname#1\endgather
 {\gmeasure@#1\endgather\tabskip\centering@
 \halign to\displaywidth{\hfil\strut@\setboxz@h{$\m@th\displaystyle{##}$}%
 \global\gwidth@\wdz@\boxz@\hfil&
 \setboxz@h{\strut@{\maketag@##\maketag@}}%
 \dimen@\displaywidth\advance\dimen@-\gwidth@
 \ifdim\dimen@>\tw@\wdz@\llap{\boxz@}\else
 \llap{\vtop{\normalbaselines\null\boxz@}}\fi
 \tabskip\z@skip\crcr#1\crcr\black@\gmaxwidth@}}
\newdimen\glineht@
\expandafter\def\csname gather \endcsname#1\endgather{\gmeasure@#1\endgather
 \ifdim\gmaxwidth@>\displaywidth\let\gdisplaywidth@\gmaxwidth@\else
 \let\gdisplaywidth@\displaywidth\fi\tabskip\centering@\halign to\displaywidth
 {\hfil\strut@\setboxz@h{$\m@th\displaystyle{##}$}%
 \global\gwidth@\wdz@\global\glineht@\ht\z@\boxz@\hfil&\kern-\gdisplaywidth@
 \setboxz@h{\strut@{\maketag@##\maketag@}}%
 \dimen@\displaywidth\advance\dimen@-\gwidth@
 \ifdim\dimen@>\tw@\wdz@\rlap{\boxz@}\else
 \rlap{\vbox{\normalbaselines\boxz@\vbox to\glineht@{}}}\fi
 \tabskip\gdisplaywidth@\crcr#1\crcr\black@\gmaxwidth@}}
\newif\ifctagsplit@
\def\CenteredTagsOnSplits{\global\ctagsplit@true}
\def\TopOrBottomTagsOnSplits{\global\ctagsplit@false}
\TopOrBottomTagsOnSplits
\def\split{\relax\ifinany@\let\next@\insplit@\else
 \ifmmode\ifinner\def\next@{\onlydmatherr@\split}\else
 \let\next@\outsplit@\fi\else
 \def\next@{\onlydmatherr@\split}\fi\fi\next@}
\def\insplit@{\global\setbox\z@\vbox\bgroup\vspace@\Let@\ialign\bgroup
 \hfil\strut@$\m@th\displaystyle{##}$&$\m@th\displaystyle{{}##}$\hfill\crcr}
\def\endsplit{\crcr\egroup\egroup\iftagsleft@\expandafter\lendsplit@\else
 \expandafter\rendsplit@\fi}
\def\rendsplit@{\global\setbox9 \vbox
 {\unvcopy\z@\global\setbox8 \lastbox\unskip}
 \setbox@ne\hbox{\unhcopy8 \unskip\global\setbox\tw@\lastbox
 \unskip\global\setbox\thr@@\lastbox}
 \global\setbox7 \hbox{\unhbox\tw@\unskip}
 \ifinalign@\ifctagsplit@                                                   
  \gdef\split@{\hbox to\wd\thr@@{}&
   \vcenter{\vbox{\moveleft\wd\thr@@\boxz@}}}
 \else\gdef\split@{&\vbox{\moveleft\wd\thr@@\box9}\crcr
  \box\thr@@&\box7}\fi                                                      
 \else                                                                      
  \ifctagsplit@\gdef\split@{\vcenter{\boxz@}}\else
  \gdef\split@{\box9\crcr\hbox{\box\thr@@\box7}}\fi
 \fi
 \split@}                                                                   
\def\lendsplit@{\global\setbox9\vtop{\unvcopy\z@}
 \setbox@ne\vbox{\unvcopy\z@\global\setbox8\lastbox}
 \setbox@ne\hbox{\unhcopy8\unskip\setbox\tw@\lastbox
  \unskip\global\setbox\thr@@\lastbox}
 \ifinalign@\ifctagsplit@                                                   
  \gdef\split@{\hbox to\wd\thr@@{}&
  \vcenter{\vbox{\moveleft\wd\thr@@\box9}}}
  \else                                                                     
  \gdef\split@{\hbox to\wd\thr@@{}&\vbox{\moveleft\wd\thr@@\box9}}\fi
 \else
  \ifctagsplit@\gdef\split@{\vcenter{\box9}}\else
  \gdef\split@{\box9}\fi
 \fi\split@}
\def\outsplit@#1$${\align\insplit@#1\endalign$$}
\newdimen\multlinegap@
\multlinegap@1em
\newdimen\multlinetaggap@
\multlinetaggap@1em
\def\MultlineGap#1{\global\multlinegap@#1\relax}
\def\multlinegap#1{\RIfMIfI@\onlydmatherr@\multlinegap\else
 \multlinegap@#1\relax\fi\else\onlydmatherr@\multlinegap\fi}
\def\nomultlinegap{\multlinegap{\z@}}
\def\multline{\RIfMIfI@
 \DN@{\onlydmatherr@\multline}\else
 \DN@{\multline@}\fi\else
 \DN@{\onlydmatherr@\multline}\fi\next@}
\newif\iftagin@
\def\tagin@#1{\tagin@false\in@\tag{#1}\ifin@\tagin@true\fi}
\def\multline@#1$${\inany@true\vspace@\allowdisplaybreak@\displaybreak@
 \tagin@{#1}\iftagsleft@\DN@{\multline@l#1$$}\else
 \DN@{\multline@r#1$$}\fi\next@}
\newdimen\mwidth@
\def\rmmeasure@#1\endmultline{%
 \def\shoveleft##1{##1}\def\shoveright##1{##1}
 \setbox@ne\vbox{\Let@\halign{\setboxz@h
  {$\m@th\@lign\displaystyle{}##$}\global\mwidth@\wdz@
  \crcr#1\crcr}}}
\newdimen\mlineht@
\newif\ifzerocr@
\newif\ifonecr@
\def\lmmeasure@#1\endmultline{\global\zerocr@true\global\onecr@false
 \everycr{\noalign{\ifonecr@\global\onecr@false\fi
  \ifzerocr@\global\zerocr@false\global\onecr@true\fi}}
  \def\shoveleft##1{##1}\def\shoveright##1{##1}%
 \setbox@ne\vbox{\Let@\halign{\setboxz@h
  {$\m@th\@lign\displaystyle{}##$}\ifonecr@\global\mwidth@\wdz@
  \global\mlineht@\ht\z@\fi\crcr#1\crcr}}}
\newbox\mtagbox@
\newdimen\ltwidth@
\newdimen\rtwidth@
\def\multline@l#1$${\iftagin@\DN@{\lmultline@@#1$$}\else
 \DN@{\setbox\mtagbox@\null\ltwidth@\z@\rtwidth@\z@
  \lmultline@@@#1$$}\fi\next@}
\def\lmultline@@#1\endmultline\tag#2$${%
 \setbox\mtagbox@\hbox{\maketag@#2\maketag@}
 \lmmeasure@#1\endmultline\dimen@\mwidth@\advance\dimen@\wd\mtagbox@
 \advance\dimen@\multlinetaggap@                                            
 \ifdim\dimen@>\displaywidth\ltwidth@\z@\else\ltwidth@\wd\mtagbox@\fi       
 \lmultline@@@#1\endmultline$$}
\def\lmultline@@@{\displ@y
 \def\shoveright##1{##1\hfilneg\hskip\multlinegap@}%
 \def\shoveleft##1{\setboxz@h{$\m@th\displaystyle{}##1$}%
  \setbox@ne\hbox{$\m@th\displaystyle##1$}%
  \hfilneg
  \iftagin@
   \ifdim\ltwidth@>\z@\hskip\ltwidth@\hskip\multlinetaggap@\fi
  \else\hskip\multlinegap@\fi\hskip.5\wd@ne\hskip-.5\wdz@##1}
  \halign\bgroup\Let@\hbox to\displaywidth
   {\strut@$\m@th\displaystyle\hfil{}##\hfil$}\crcr
   \hfilneg                                                                 
   \iftagin@                                                                
    \ifdim\ltwidth@>\z@                                                     
     \box\mtagbox@\hskip\multlinetaggap@                                    
    \else
     \rlap{\vbox{\normalbaselines\hbox{\strut@\box\mtagbox@}%
     \vbox to\mlineht@{}}}\fi                                               
   \else\hskip\multlinegap@\fi}                                             
\def\multline@r#1$${\iftagin@\DN@{\rmultline@@#1$$}\else
 \DN@{\setbox\mtagbox@\null\ltwidth@\z@\rtwidth@\z@
  \rmultline@@@#1$$}\fi\next@}
\def\rmultline@@#1\endmultline\tag#2$${\ltwidth@\z@
 \setbox\mtagbox@\hbox{\maketag@#2\maketag@}%
 \rmmeasure@#1\endmultline\dimen@\mwidth@\advance\dimen@\wd\mtagbox@
 \advance\dimen@\multlinetaggap@
 \ifdim\dimen@>\displaywidth\rtwidth@\z@\else\rtwidth@\wd\mtagbox@\fi
 \rmultline@@@#1\endmultline$$}
\def\rmultline@@@{\displ@y
 \def\shoveright##1{##1\hfilneg\iftagin@\ifdim\rtwidth@>\z@
  \hskip\rtwidth@\hskip\multlinetaggap@\fi\else\hskip\multlinegap@\fi}%
 \def\shoveleft##1{\setboxz@h{$\m@th\displaystyle{}##1$}%
  \setbox@ne\hbox{$\m@th\displaystyle##1$}%
  \hfilneg\hskip\multlinegap@\hskip.5\wd@ne\hskip-.5\wdz@##1}%
 \halign\bgroup\Let@\hbox to\displaywidth
  {\strut@$\m@th\displaystyle\hfil{}##\hfil$}\crcr
 \hfilneg\hskip\multlinegap@}
\def\endmultline{\iftagsleft@\expandafter\lendmultline@\else
 \expandafter\rendmultline@\fi}
\def\lendmultline@{\hfilneg\hskip\multlinegap@\crcr\egroup}
\def\rendmultline@{\iftagin@                                                
 \ifdim\rtwidth@>\z@                                                        
  \hskip\multlinetaggap@\box\mtagbox@                                       
 \else\llap{\vtop{\normalbaselines\null\hbox{\strut@\box\mtagbox@}}}\fi     
 \else\hskip\multlinegap@\fi                                                
 \hfilneg\crcr\egroup}
\def\bmod{\mskip-\medmuskip\mkern5mu\mathbin{\fam\z@ mod}\penalty900
 \mkern5mu\mskip-\medmuskip}
\def\pmod#1{\allowbreak\ifinner\mkern8mu\else\mkern18mu\fi
 ({\fam\z@ mod}\,\,#1)}
\def\pod#1{\allowbreak\ifinner\mkern8mu\else\mkern18mu\fi(#1)}
\def\mod#1{\allowbreak\ifinner\mkern12mu\else\mkern18mu\fi{\fam\z@ mod}\,\,#1}
\newcount\cfraccount@
\def\cfrac{\bgroup\bgroup\advance\cfraccount@\@ne\strut
 \iffalse{\fi\def\\{\over\displaystyle}\iffalse}\fi}
\def\lcfrac{\bgroup\bgroup\advance\cfraccount@\@ne\strut
 \iffalse{\fi\def\\{\hfill\over\displaystyle}\iffalse}\fi}
\def\rcfrac{\bgroup\bgroup\advance\cfraccount@\@ne\strut\hfill
 \iffalse{\fi\def\\{\over\displaystyle}\iffalse}\fi}
\def\gloop@#1\repeat{\gdef\body{#1}\iterate}
\def\endcfrac{\gloop@\ifnum\cfraccount@>\z@\global\advance\cfraccount@\m@ne
 \egroup\hskip-\nulldelimiterspace\egroup\repeat}
\def\binrel@#1{\setboxz@h{\thinmuskip0mu
  \medmuskip\m@ne mu\thickmuskip\@ne mu$#1\m@th$}%
 \setbox@ne\hbox{\thinmuskip0mu\medmuskip\m@ne mu\thickmuskip
  \@ne mu${}#1{}\m@th$}%
 \setbox\tw@\hbox{\hskip\wd@ne\hskip-\wdz@}}
\def\overset#1\to#2{\binrel@{#2}\ifdim\wd\tw@<\z@
 \mathbin{\mathop{\kern\z@#2}\limits^{#1}}\else\ifdim\wd\tw@>\z@
 \mathrel{\mathop{\kern\z@#2}\limits^{#1}}\else
 {\mathop{\kern\z@#2}\limits^{#1}}{}\fi\fi}
\def\underset#1\to#2{\binrel@{#2}\ifdim\wd\tw@<\z@
 \mathbin{\mathop{\kern\z@#2}\limits_{#1}}\else\ifdim\wd\tw@>\z@
 \mathrel{\mathop{\kern\z@#2}\limits_{#1}}\else
 {\mathop{\kern\z@#2}\limits_{#1}}{}\fi\fi}
\def\oversetbrace#1\to#2{\overbrace{#2}^{#1}}
\def\undersetbrace#1\to#2{\underbrace{#2}_{#1}}
\def\sideset#1\and#2\to#3{%
 \setbox@ne\hbox{$\dsize{\vphantom{#3}}#1{#3}\m@th$}%
 \setbox\tw@\hbox{$\dsize{#3}#2\m@th$}%
 \hskip\wd@ne\hskip-\wd\tw@\mathop{\hskip\wd\tw@\hskip-\wd@ne
  {\vphantom{#3}}#1{#3}#2}}
\def\rightarrowfill@#1{$#1\m@th\mathord-\mkern-6mu\cleaders
 \hbox{$#1\mkern-2mu\mathord-\mkern-2mu$}\hfill
 \mkern-6mu\mathord\rightarrow$}
\def\leftarrowfill@#1{$#1\m@th\mathord\leftarrow\mkern-6mu\cleaders
 \hbox{$#1\mkern-2mu\mathord-\mkern-2mu$}\hfill\mkern-6mu\mathord-$}
\def\leftrightarrowfill@#1{$#1\m@th\mathord\leftarrow\mkern-6mu\cleaders
 \hbox{$#1\mkern-2mu\mathord-\mkern-2mu$}\hfill
 \mkern-6mu\mathord\rightarrow$}
\def\overrightarrow{\mathpalette\overrightarrow@}
\def\overrightarrow@#1#2{\vbox{\ialign{##\crcr\rightarrowfill@#1\crcr
 \noalign{\kern-\ex@\nointerlineskip}$\m@th\hfil#1#2\hfil$\crcr}}}

\def\overleftarrow{\mathpalette\overleftarrow@}
\def\overleftarrow@#1#2{\vbox{\ialign{##\crcr\leftarrowfill@#1\crcr
 \noalign{\kern-\ex@\nointerlineskip}$\m@th\hfil#1#2\hfil$\crcr}}}
\def\overleftrightarrow{\mathpalette\overleftrightarrow@}
\def\overleftrightarrow@#1#2{\vbox{\ialign{##\crcr\leftrightarrowfill@#1\crcr
 \noalign{\kern-\ex@\nointerlineskip}$\m@th\hfil#1#2\hfil$\crcr}}}
\def\underrightarrow{\mathpalette\underrightarrow@}
\def\underrightarrow@#1#2{\vtop{\ialign{##\crcr$\m@th\hfil#1#2\hfil$\crcr
 \noalign{\nointerlineskip}\rightarrowfill@#1\crcr}}}

\def\underleftarrow{\mathpalette\underleftarrow@}
\def\underleftarrow@#1#2{\vtop{\ialign{##\crcr$\m@th\hfil#1#2\hfil$\crcr
 \noalign{\nointerlineskip}\leftarrowfill@#1\crcr}}}
\def\underleftrightarrow{\mathpalette\underleftrightarrow@}
\def\underleftrightarrow@#1#2{\vtop{\ialign{##\crcr$\m@th\hfil#1#2\hfil$\crcr
 \noalign{\nointerlineskip}\leftrightarrowfill@#1\crcr}}}
\let\DOTSI\relax
\let\DOTSB\relax

\newif\ifmath@
{\uccode`7=`\\ \uccode`8=`m \uccode`9=`a \uccode`0=`t \uccode`!=`h
 \uppercase{\gdef\math@#1#2#3#4#5#6\math@{\global\math@false\ifx 7#1\ifx 8#2%
 \ifx 9#3\ifx 0#4\ifx !#5\xdef\meaning@{#6}\global\math@true\fi\fi\fi\fi\fi}}}
\newif\ifmathch@
{\uccode`7=`c \uccode`8=`h \uccode`9=`\"
 \uppercase{\gdef\mathch@#1#2#3#4#5#6\mathch@{\global\mathch@false
  \ifx 7#1\ifx 8#2\ifx 9#5\global\mathch@true\xdef\meaning@{9#6}\fi\fi\fi}}}
\newcount\classnum@
\def\getmathch@#1.#2\getmathch@{\classnum@#1 \divide\classnum@4096
 \ifcase\number\classnum@\or\or\gdef\thedots@{\dotsb@}\or
 \gdef\thedots@{\dotsb@}\fi}
\newif\ifmathbin@
{\uccode`4=`b \uccode`5=`i \uccode`6=`n
 \uppercase{\gdef\mathbin@#1#2#3{\relaxnext@
  \DNii@##1\mathbin@{\ifx\space@\next\global\mathbin@true\fi}%
 \global\mathbin@false\DN@##1\mathbin@{}%
 \ifx 4#1\ifx 5#2\ifx 6#3\DN@{\FN@\nextii@}\fi\fi\fi\next@}}}
\newif\ifmathrel@
{\uccode`4=`r \uccode`5=`e \uccode`6=`l
 \uppercase{\gdef\mathrel@#1#2#3{\relaxnext@
  \DNii@##1\mathrel@{\ifx\space@\next\global\mathrel@true\fi}%
 \global\mathrel@false\DN@##1\mathrel@{}%
 \ifx 4#1\ifx 5#2\ifx 6#3\DN@{\FN@\nextii@}\fi\fi\fi\next@}}}
\newif\ifmacro@
{\uccode`5=`m \uccode`6=`a \uccode`7=`c
 \uppercase{\gdef\macro@#1#2#3#4\macro@{\global\macro@false
  \ifx 5#1\ifx 6#2\ifx 7#3\global\macro@true
  \xdef\meaning@{\macro@@#4\macro@@}\fi\fi\fi}}}
\def\macro@@#1->#2\macro@@{#2}
\newif\ifDOTS@
\newcount\DOTSCASE@
{\uccode`6=`\\ \uccode`7=`D \uccode`8=`O \uccode`9=`T \uccode`0=`S
 \uppercase{\gdef\DOTS@#1#2#3#4#5{\global\DOTS@false\DN@##1\DOTS@{}%
  \ifx 6#1\ifx 7#2\ifx 8#3\ifx 9#4\ifx 0#5\let\next@\DOTS@@\fi\fi\fi\fi\fi
  \next@}}}
{\uccode`3=`B \uccode`4=`I \uccode`5=`X
 \uppercase{\gdef\DOTS@@#1{\relaxnext@
  \DNii@##1\DOTS@{\ifx\space@\next\global\DOTS@true\fi}%
  \DN@{\FN@\nextii@}%
  \ifx 3#1\global\DOTSCASE@\z@\else
  \ifx 4#1\global\DOTSCASE@\@ne\else
  \ifx 5#1\global\DOTSCASE@\tw@\else\DN@##1\DOTS@{}%
  \fi\fi\fi\next@}}}
\newif\ifnot@
{\uccode`5=`\\ \uccode`6=`n \uccode`7=`o \uccode`8=`t
 \uppercase{\gdef\not@#1#2#3#4{\relaxnext@
  \DNii@##1\not@{\ifx\space@\next\global\not@true\fi}%
 \global\not@false\DN@##1\not@{}%
 \ifx 5#1\ifx 6#2\ifx 7#3\ifx 8#4\DN@{\FN@\nextii@}\fi\fi\fi
 \fi\next@}}}
\newif\ifkeybin@
\def\keybin@{\keybin@true
 \ifx\next+\else\ifx\next=\else\ifx\next<\else\ifx\next>\else\ifx\next-\else
 \ifx\next*\else\ifx\next:\else\keybin@false\fi\fi\fi\fi\fi\fi\fi}
\def\dots{\RIfM@\expandafter\mdots@\else\expandafter\tdots@\fi}
\def\tdots@{\unskip\relaxnext@
 \DN@{$\m@th\mathinner{\ldotp\ldotp\ldotp}\,
   \ifx\next,\,$\else\ifx\next.\,$\else\ifx\next;\,$\else\ifx\next:\,$\else
   \ifx\next?\,$\else\ifx\next!\,$\else$ \fi\fi\fi\fi\fi\fi}%
 \ \FN@\next@}
\def\mdots@{\FN@\mdots@@}
\def\mdots@@{\gdef\thedots@{\dotso@}
 \ifx\next\boldkey\gdef\thedots@\boldkey{\boldkeydots@}\else                
 \ifx\next\boldsymbol\gdef\thedots@\boldsymbol{\boldsymboldots@}\else       
 \ifx,\next\gdef\thedots@{\dotsc}
 \else\ifx\not\next\gdef\thedots@{\dotsb@}
 \else\keybin@
 \ifkeybin@\gdef\thedots@{\dotsb@}
 \else\xdef\meaning@{\meaning\next..........}\xdef\meaning@@{\meaning@}
  \expandafter\math@\meaning@\math@
  \ifmath@
   \expandafter\mathch@\meaning@\mathch@
   \ifmathch@\expandafter\getmathch@\meaning@\getmathch@\fi                 
  \else\expandafter\macro@\meaning@@\macro@                                 
  \ifmacro@                                                                
   \expandafter\not@\meaning@\not@\ifnot@\gdef\thedots@{\dotsb@}
  \else\expandafter\DOTS@\meaning@\DOTS@
  \ifDOTS@
   \ifcase\number\DOTSCASE@\gdef\thedots@{\dotsb@}%
    \or\gdef\thedots@{\dotsi}\else\fi                                      
  \else\expandafter\math@\meaning@\math@                                   
  \ifmath@\expandafter\mathbin@\meaning@\mathbin@
  \ifmathbin@\gdef\thedots@{\dotsb@}
  \else\expandafter\mathrel@\meaning@\mathrel@
  \ifmathrel@\gdef\thedots@{\dotsb@}
  \fi\fi\fi\fi\fi\fi\fi\fi\fi\fi\fi\fi
 \thedots@}
\def\plainldots@{\mathinner{\ldotp\ldotp\ldotp}}
\def\plaincdots@{\mathinner{\cdotp\cdotp\cdotp}}
\def\dotsi{\!\plaincdots@}
\let\dotsb@\plaincdots@
\newif\ifextra@
\newif\ifrightdelim@
\def\rightdelim@{\global\rightdelim@true                                    
 \ifx\next)\else                                                            
 \ifx\next]\else
 \ifx\next\rbrack\else
 \ifx\next\}\else
 \ifx\next\rbrace\else
 \ifx\next\rangle\else
 \ifx\next\rceil\else
 \ifx\next\rfloor\else
 \ifx\next\rgroup\else
 \ifx\next\rmoustache\else
 \ifx\next\right\else
 \ifx\next\bigr\else
 \ifx\next\biggr\else
 \ifx\next\Bigr\else                                                        
 \ifx\next\Biggr\else\global\rightdelim@false
 \fi\fi\fi\fi\fi\fi\fi\fi\fi\fi\fi\fi\fi\fi\fi}
\def\extra@{%
 \global\extra@false\rightdelim@\ifrightdelim@\global\extra@true            
 \else\ifx\next$\global\extra@true                                          
 \else\xdef\meaning@{\meaning\next..........}
 \expandafter\macro@\meaning@\macro@\ifmacro@                               
 \expandafter\DOTS@\meaning@\DOTS@
 \ifDOTS@
 \ifnum\DOTSCASE@=\tw@\global\extra@true                                    
 \fi\fi\fi\fi\fi}
\newif\ifbold@
\def\dotso@{\relaxnext@
 \ifbold@
  \let\next\delayed@
  \DNii@{\extra@\plainldots@\ifextra@\,\fi}%
 \else
  \DNii@{\DN@{\extra@\plainldots@\ifextra@\,\fi}\FN@\next@}%
 \fi
 \nextii@}
\def\extrap@#1{%
 \ifx\next,\DN@{#1\,}\else
 \ifx\next;\DN@{#1\,}\else
 \ifx\next.\DN@{#1\,}\else\extra@
 \ifextra@\DN@{#1\,}\else
 \let\next@#1\fi\fi\fi\fi\next@}
\def\ldots{\DN@{\extrap@\plainldots@}%
 \FN@\next@}
\def\cdots{\DN@{\extrap@\plaincdots@}%
 \FN@\next@}

\def\dotsc{\relaxnext@
 \DN@{\ifx\next;\plainldots@\,\else
  \ifx\next.\plainldots@\,\else\extra@\plainldots@
  \ifextra@\,\fi\fi\fi}%
 \FN@\next@}
\def\cdot{\mathchar"2201 }
\def\longrightarrow{\DOTSB\relbar\joinrel\rightarrow}

\def\mapsto{\DOTSB\mapstochar\rightarrow}

\def\dddot#1{{\mathop{#1}\limits^{\vbox to-1.4\ex@{\kern-\tw@\ex@
 \hbox{\rm...}\vss}}}}
\def\ddddot#1{{\mathop{#1}\limits^{\vbox to-1.4\ex@{\kern-\tw@\ex@
 \hbox{\rm....}\vss}}}}
\def\sphat{^{\mathchoice{}{}%
 {\,\,\botsmash{\hbox{\lower4\ex@\hbox{$\m@th\widehat{\null}$}}}}%
 {\,\botsmash{\hbox{\lower3\ex@\hbox{$\m@th\hat{\null}$}}}}}}

\def\spacute{^{\!\botsmash{\hbox{\lower\@ne ex\hbox{\'{}}}}}}
\def\spgrave{^{\mathchoice{}{}{}{\!}%
 \botsmash{\hbox{\lower\@ne ex\hbox{\`{}}}}}}
\def\spdot{^{\hbox{\raise\ex@\hbox{\rm.}}}}
\def\spddot{^{\hbox{\raise\ex@\hbox{\rm..}}}}
\def\spdddot{^{\hbox{\raise\ex@\hbox{\rm...}}}}
\def\spddddot{^{\hbox{\raise\ex@\hbox{\rm....}}}}
\def\spbreve{^{\!\botsmash{\hbox{\lower4\ex@\hbox{\u{}}}}}}

\def\textonlyfont@#1#2{\def#1{\RIfM@
 \Err@{Use \string#1\space only in text}\else#2\fi}}
\textonlyfont@\rm\tenrm
\textonlyfont@\it\tenit
\textonlyfont@\sl\tensl
\textonlyfont@\bf\tenbf
\def\oldnos#1{\RIfM@{\mathcode`\,="013B \fam\@ne#1}\else
 \leavevmode\hbox{$\m@th\mathcode`\,="013B \fam\@ne#1$}\fi}
\def\text{\RIfM@\expandafter\text@\else\expandafter\text@@\fi}
\def\text@@#1{\leavevmode\hbox{#1}}
\def\mathhexbox@#1#2#3{\text{$\m@th\mathchar"#1#2#3$}}
\def\dag{{\mathhexbox@279}}
\def\ddag{{\mathhexbox@27A}}
\def\S{{\mathhexbox@278}}
\def\P{{\mathhexbox@27B}}
\newif\iffirstchoice@
\firstchoice@true
\def\text@#1{\mathchoice
 {\hbox{\everymath{\displaystyle}\def\textfonti{\the\textfont\@ne}%
  \def\textfontii{\the\textfont\tw@}\textdef@@ T#1}}
 {\hbox{\firstchoice@false
  \everymath{\textstyle}\def\textfonti{\the\textfont\@ne}%
  \def\textfontii{\the\textfont\tw@}\textdef@@ T#1}}
 {\hbox{\firstchoice@false
  \everymath{\scriptstyle}\def\textfonti{\the\scriptfont\@ne}%
  \def\textfontii{\the\scriptfont\tw@}\textdef@@ S\rm#1}}
 {\hbox{\firstchoice@false
  \everymath{\scriptscriptstyle}\def\textfonti
  {\the\scriptscriptfont\@ne}%
  \def\textfontii{\the\scriptscriptfont\tw@}\textdef@@ s\rm#1}}}
\def\textdef@@#1{\textdef@#1\rm\textdef@#1\bf\textdef@#1\sl\textdef@#1\it}
\def\rmfam{0}
\def\textdef@#1#2{%
 \DN@{\csname\expandafter\eat@\string#2fam\endcsname}%
 \if S#1\edef#2{\the\scriptfont\next@\relax}%
 \else\if s#1\edef#2{\the\scriptscriptfont\next@\relax}%
 \else\edef#2{\the\textfont\next@\relax}\fi\fi}
\scriptfont\itfam\tenit \scriptscriptfont\itfam\tenit
\scriptfont\slfam\tensl \scriptscriptfont\slfam\tensl
\newif\iftopfolded@
\newif\ifbotfolded@
\def\topfoldedtext{\topfolded@true\botfolded@false\foldedtext@}
\def\botfoldedtext{\botfolded@true\topfolded@false\foldedtext@}
\def\foldedtext{\topfolded@false\botfolded@false\foldedtext@}
\Invalid@\foldedwidth
\def\foldedtext@{\relaxnext@
 \DN@{\ifx\next\foldedwidth\let\next@\nextii@\else
  \DN@{\nextii@\foldedwidth{.3\hsize}}\fi\next@}%
 \DNii@\foldedwidth##1##2{\setbox\z@\vbox
  {\normalbaselines\hsize##1\relax
  \tolerance1600 \noindent\ignorespaces##2}\ifbotfolded@\boxz@\else
  \iftopfolded@\vtop{\unvbox\z@}\else\vcenter{\boxz@}\fi\fi}%
 \FN@\next@}
\def\bold{\RIfM@\expandafter\bold@\else
 \expandafter\nonmatherr@\expandafter\bold\fi}
\def\bold@#1{{\bold@@{#1}}}
\def\bold@@#1{\fam\bffam\relax#1}
\def\slanted{\RIfM@\expandafter\slanted@\else
 \expandafter\nonmatherr@\expandafter\slanted\fi}
\def\slanted@#1{{\slanted@@{#1}}}
\def\slanted@@#1{\fam\slfam\relax#1}
\def\roman{\RIfM@\expandafter\roman@\else
 \expandafter\nonmatherr@\expandafter\roman\fi}
\def\roman@#1{{\roman@@{#1}}}
\def\roman@@#1{\fam\rmfam\relax#1}
\def\italic{\RIfM@\expandafter\italic@\else
 \expandafter\nonmatherr@\expandafter\italic\fi}
\def\italic@#1{{\italic@@{#1}}}
\def\italic@@#1{\fam\itfam\relax#1}
\def\Cal{\RIfM@\expandafter\Cal@\else
 \expandafter\nonmatherr@\expandafter\Cal\fi}
\def\Cal@#1{{\Cal@@{#1}}}
\def\Cal@@#1{\noaccents@\fam\tw@#1}
\mathchardef\Gamma="0000
\mathchardef\Delta="0001
\mathchardef\Theta="0002
\mathchardef\Lambda="0003
\mathchardef\Xi="0004
\mathchardef\Pi="0005
\mathchardef\Sigma="0006
\mathchardef\Upsilon="0007
\mathchardef\Phi="0008
\mathchardef\Psi="0009
\mathchardef\Omega="000A
\mathchardef\varGamma="0100
\mathchardef\varDelta="0101
\mathchardef\varTheta="0102
\mathchardef\varLambda="0103
\mathchardef\varXi="0104
\mathchardef\varPi="0105
\mathchardef\varSigma="0106
\mathchardef\varUpsilon="0107
\mathchardef\varPhi="0108
\mathchardef\varPsi="0109
\mathchardef\varOmega="010A
\newif\ifmsamloaded@
\newif\ifmsbmloaded@
\newif\ifeufmloaded@
\let\alloc@@\alloc@
\def\hexnumber@#1{\ifcase#1 0\or 1\or 2\or 3\or 4\or 5\or 6\or 7\or 8\or
 9\or A\or B\or C\or D\or E\or F\fi}
\edef\bffam@{\hexnumber@\bffam}
\def\loadmsam{\msamloaded@true
 \font@\tenmsa=msam10
 \font@\sevenmsa=msam7
 \font@\fivemsa=msam5
 \alloc@@8\fam\chardef\sixt@@n\msafam
 \textfont\msafam=\tenmsa
 \scriptfont\msafam=\sevenmsa
 \scriptscriptfont\msafam=\fivemsa
 \edef\msafam@{\hexnumber@\msafam}%
 \mathchardef\dabar@"0\msafam@39
 \def\dashrightarrow{\mathrel{\dabar@\dabar@\mathchar"0\msafam@4B}}%
 \def\dashleftarrow{\mathrel{\mathchar"0\msafam@4C\dabar@\dabar@}}%
 \let\dasharrow\dashrightarrow
 \def\ulcorner{\delimiter"4\msafam@70\msafam@70 }
 \def\urcorner{\delimiter"5\msafam@71\msafam@71 }
 \def\llcorner{\delimiter"4\msafam@78\msafam@78 }
 \def\lrcorner{\delimiter"5\msafam@79\msafam@79 }
 \def\yen{{\mathhexbox@\msafam@55 }}
 \def\checkmark{{\mathhexbox@\msafam@58 }}
 \def\circledR{{\mathhexbox@\msafam@72 }}
 \def\maltese{{\mathhexbox@\msafam@7A }}}
\def\loadmsbm{\msbmloaded@true
 \font@\tenmsb=msbm10
 \font@\sevenmsb=msbm7
 \font@\fivemsb=msbm5
 \alloc@@8\fam\chardef\sixt@@n\msbfam
 \textfont\msbfam=\tenmsb
 \scriptfont\msbfam=\sevenmsb
 \scriptscriptfont\msbfam=\fivemsb
 \edef\msbfam@{\hexnumber@\msbfam}%
 }
\def\widehat#1{\ifmsbmloaded@
  \setboxz@h{$\m@th#1$}\ifdim\wdz@>\tw@ em\mathaccent"0\msbfam@5B{#1}\else
  \mathaccent"0362{#1}\fi
 \else\mathaccent"0362{#1}\fi}
\def\widetilde#1{\ifmsbmloaded@
  \setboxz@h{$\m@th#1$}\ifdim\wdz@>\tw@ em\mathaccent"0\msbfam@5D{#1}\else
  \mathaccent"0365{#1}\fi
 \else\mathaccent"0365{#1}\fi}
\def\newsymbol#1#2#3#4#5{\define#1{}\let\next@\relax
 \ifnum#2=\@ne\ifmsamloaded@\let\next@\msafam@\fi\else
 \ifnum#2=\tw@\ifmsbmloaded@\let\next@\msbfam@\fi\fi\fi
 \ifx\next@\relax
  \ifnum#2>\tw@\Err@{\Invalid@@\string\newsymbol}\else
  \ifnum#2=\@ne\Err@{You must first \string\loadmsam}\else
   \Err@{You must first \string\loadmsbm}\fi\fi
 \else
  \mathchardef#1="#3\next@#4#5
 \fi}

\def\Bbb{\RIfM@\expandafter\Bbb@\else
 \expandafter\nonmatherr@\expandafter\Bbb\fi}
\def\Bbb@#1{{\Bbb@@{#1}}}
\def\Bbb@@#1{\noaccents@\fam\msbfam\relax#1}
\def\loadeufm{\eufmloaded@true
 \font@\teneufm=eufm10
 \font@\seveneufm=eufm7
 \font@\fiveeufm=eufm5
 \alloc@@8\fam\chardef\sixt@@n\eufmfam
 \textfont\eufmfam=\teneufm
 \scriptfont\eufmfam=\seveneufm
 \scriptscriptfont\eufmfam=\fiveeufm}
\def\frak{\RIfM@\expandafter\frak@\else
 \expandafter\nonmatherr@\expandafter\frak\fi}
\def\frak@#1{{\frak@@{#1}}}
\def\frak@@#1{\fam\eufmfam\relax#1}

\newif\ifcmmibloaded@
\newif\ifcmbsyloaded@
\def\loadbold{\cmmibloaded@true\cmbsyloaded@true
 \font@\tencmmib=cmmib10 \font@\sevencmmib=cmmib7 \font@\fivecmmib=cmmib5
 \skewchar\tencmmib='177 \skewchar\sevencmmib='177 \skewchar\fivecmmib='177
 \alloc@@8\fam\chardef\sixt@@n\cmmibfam
 \textfont\cmmibfam=\tencmmib
 \scriptfont\cmmibfam=\sevencmmib
 \scriptscriptfont\cmmibfam=\fivecmmib
 \edef\cmmibfam@{\hexnumber@\cmmibfam}%
 \font@\tencmbsy=cmbsy10 \font@\sevencmbsy=cmbsy7 \font@\fivecmbsy=cmbsy5
 \skewchar\tencmbsy='60 \skewchar\sevencmbsy='60 \skewchar\fivecmbsy='60
 \alloc@@8\fam\chardef\sixt@@n\cmbsyfam
 \textfont\cmbsyfam=\tencmbsy
 \scriptfont\cmbsyfam=\sevencmbsy
 \scriptscriptfont\cmbsyfam=\fivecmbsy
 \edef\cmbsyfam@{\hexnumber@\cmbsyfam}}
\def\mathchari@#1#2#3{\ifcmmibloaded@\mathchar"#1\cmmibfam@#2#3 \else
 \Err@{First bold symbol font not loaded}\fi}
\def\mathcharii@#1#2#3{\ifcmbsyloaded@\mathchar"#1\cmbsyfam@#2#3 \else
 \Err@{Second bold symbol font not loaded}\fi}
\def\boldkey#1{\ifcat\noexpand#1A%
  \ifcmmibloaded@{\fam\cmmibfam#1}\else
   \Err@{First bold symbol font not loaded}\fi
 \else
 \ifx#1!\mathchar"5\bffam@21 \else
 \ifx#1(\mathchar"4\bffam@28 \else\ifx#1)\mathchar"5\bffam@29 \else
 \ifx#1+\mathchar"2\bffam@2B \else\ifx#1:\mathchar"3\bffam@3A \else
 \ifx#1;\mathchar"6\bffam@3B \else\ifx#1=\mathchar"3\bffam@3D \else
 \ifx#1?\mathchar"5\bffam@3F \else\ifx#1[\mathchar"4\bffam@5B \else
 \ifx#1]\mathchar"5\bffam@5D \else
 \ifx#1,\mathchari@63B \else
 \ifx#1-\mathcharii@200 \else
 \ifx#1.\mathchari@03A \else
 \ifx#1/\mathchari@03D \else
 \ifx#1<\mathchari@33C \else
 \ifx#1>\mathchari@33E \else
 \ifx#1*\mathcharii@203 \else
 \ifx#1|\mathcharii@06A \else
 \ifx#10\bold0\else\ifx#11\bold1\else\ifx#12\bold2\else\ifx#13\bold3\else
 \ifx#14\bold4\else\ifx#15\bold5\else\ifx#16\bold6\else\ifx#17\bold7\else
 \ifx#18\bold8\else\ifx#19\bold9\else
  \Err@{\string\boldkey\space can't be used with #1}%
 \fi\fi\fi\fi\fi\fi\fi\fi\fi\fi\fi\fi\fi\fi\fi
 \fi\fi\fi\fi\fi\fi\fi\fi\fi\fi\fi\fi\fi\fi}
\def\boldsymbol#1{%
 \DN@{\Err@{You can't use \string\boldsymbol\space with \string#1}#1}%
 \ifcat\noexpand#1A%
   \let\next@\relax
  \ifcmmibloaded@{\fam\cmmibfam#1}\else\Err@{First bold symbol
   font not loaded}\fi
 \else
  \xdef\meaning@{\meaning#1.........}%
  \expandafter\math@\meaning@\math@
  \ifmath@
   \expandafter\mathch@\meaning@\mathch@
   \ifmathch@
    \expandafter\boldsymbol@@\meaning@\boldsymbol@@
   \fi
  \else
   \expandafter\macro@\meaning@\macro@
   \expandafter\delim@\meaning@\delim@
   \ifdelim@
    \expandafter\delim@@\meaning@\delim@@
   \else
    \boldsymbol@{#1}%
   \fi
  \fi
 \fi
 \next@}
\def\mathhexboxii@#1#2{\ifcmbsyloaded@\mathhexbox@{\cmbsyfam@}{#1}{#2}\else
  \Err@{Second bold symbol font not loaded}\fi}
\def\boldsymbol@#1{\let\next@\relax\let\next#1%
 \ifx\next\cdot\mathcharii@201 \else
 \ifx\next\prime{{\null\mathcharii@030 \null}}\else
 \ifx\next\lbrack\mathchar"4\bffam@5B \else
 \ifx\next\rbrack\mathchar"5\bffam@5D \else
 \ifx\next\{\mathcharii@466 \else
 \ifx\next\lbrace\mathcharii@466 \else
 \ifx\next\}\mathcharii@567 \else
 \ifx\next\rbrace\mathcharii@567 \else
 \ifx\next\surd{{\mathcharii@170}}\else
 \ifx\next\S{{\mathhexboxii@78}}\else
 \ifx\next\P{{\mathhexboxii@7B}}\else
 \ifx\next\dag{{\mathhexboxii@79}}\else
 \ifx\next\ddag{{\mathhexboxii@7A}}\else
 \DN@{\Err@{You can't use \string\boldsymbol\space with \string#1}#1}%
 \fi\fi\fi\fi\fi\fi\fi\fi\fi\fi\fi\fi\fi}
\def\boldsymbol@@#1.#2\boldsymbol@@{\classnum@#1 \count@@@\classnum@        
 \divide\classnum@4096 \count@\classnum@                                    
 \multiply\count@4096 \advance\count@@@-\count@ \count@@\count@@@           
 \divide\count@@@\@cclvi \count@\count@@                                    
 \multiply\count@@@\@cclvi \advance\count@@-\count@@@                       
 \divide\count@@@\@cclvi                                                    
 \multiply\classnum@4096 \advance\classnum@\count@@                         
 \ifnum\count@@@=\z@                                                        
  \count@"\bffam@ \multiply\count@\@cclvi
  \advance\classnum@\count@
  \DN@{\mathchar\number\classnum@}%
 \else
  \ifnum\count@@@=\@ne                                                      
   \ifcmmibloaded@
   \count@"\cmmibfam@ \multiply\count@\@cclvi
   \advance\classnum@\count@
   \DN@{\mathchar\number\classnum@}%
   \else\DN@{\Err@{First bold symbol font not loaded}}\fi
  \else
   \ifnum\count@@@=\tw@                                                    
  \ifcmbsyloaded@
    \count@"\cmbsyfam@ \multiply\count@\@cclvi
    \advance\classnum@\count@
    \DN@{\mathchar\number\classnum@}%
  \else\DN@{\Err@{Second bold symbol font not loaded}}\fi
  \fi
 \fi
\fi}
\newif\ifdelim@
\newcount\delimcount@
{\uccode`6=`\\ \uccode`7=`d \uccode`8=`e \uccode`9=`l
 \uppercase{\gdef\delim@#1#2#3#4#5\delim@
  {\delim@false\ifx 6#1\ifx 7#2\ifx 8#3\ifx 9#4\delim@true
   \xdef\meaning@{#5}\fi\fi\fi\fi}}}
\def\delim@@#1"#2#3#4#5#6\delim@@{\if#32%
\let\next@\relax
 \ifcmbsyloaded@
 \mathcharii@#2#4#5 \else\Err@{Second bold family not loaded}\fi\fi}
\def\vert{\delimiter"026A30C }
\def\Vert{\delimiter"026B30D }
\let\|\Vert

\def\boldkeydots@#1{\bold@true\let\next=#1\let\delayed@=#1\mdots@@
 \boldkey#1\bold@false}  
\def\boldsymboldots@#1{\bold@true\let\next#1\let\delayed@#1\mdots@@
 \boldsymbol#1\bold@false}
\newif\ifeufbloaded@
\def\loadeufb{\eufbloaded@true
 \font@\teneufb=eufb10
 \font@\seveneufb=eufb7
 \font@\fiveeufb=eufb5
 \alloc@@8\fam\chardef\sixt@@n\eufbfam
 \textfont\eufbfam=\teneufb
 \scriptfont\eufbfam=\seveneufb
 \scriptscriptfont\eufbfam=\fiveeufb
 \edef\eufbfam@{\hexnumber@\eufbfam}}
\newif\ifeusmloaded@
\def\loadeusm{\eusmloaded@true
 \font@\teneusm=eusm10
 \font@\seveneusm=eusm7
 \font@\fiveeusm=eusm5
 \alloc@@8\fam\chardef\sixt@@n\eusmfam
 \textfont\eusmfam=\teneusm
 \scriptfont\eusmfam=\seveneusm
 \scriptscriptfont\eusmfam=\fiveeusm
 \edef\eusmfam@{\hexnumber@\eusmfam}}
\newif\ifeusbloaded@
\def\loadeusb{\eusbloaded@true
 \font@\teneusb=eusb10
 \font@\seveneusb=eusb7
 \font@\fiveeusb=eusb5
 \alloc@@8\fam\chardef\sixt@@n\eusbfam
 \textfont\eusbfam=\teneusb
 \scriptfont\eusbfam=\seveneusb
 \scriptscriptfont\eusbfam=\fiveeusb
 \edef\eusbfam@{\hexnumber@\eusbfam}}
\newif\ifeurmloaded@
\def\loadeurm{\eurmloaded@true
 \font@\teneurm=eurm10
 \font@\seveneurm=eurm7
 \font@\fiveeurm=eurm5
 \alloc@@8\fam\chardef\sixt@@n\eurmfam
 \textfont\eurmfam=\teneurm
 \scriptfont\eurmfam=\seveneurm
 \scriptscriptfont\eurmfam=\fiveeurm
 \edef\eurmfam@{\hexnumber@\eurmfam}}
\newif\ifeurbloaded@
\def\loadeurb{\eurbloaded@true
 \font@\teneurb=eurb10
 \font@\seveneurb=eurb7
 \font@\fiveeurb=eurb5
 \alloc@@8\fam\chardef\sixt@@n\eurbfam
 \textfont\eurbfam=\teneurb
 \scriptfont\eurbfam=\seveneurb
 \scriptscriptfont\eurbfam=\fiveeurb
 \edef\eurbfam@{\hexnumber@\eurbfam}}
\def\accentclass@{7}
\def\noaccents@{\def\accentclass@{0}}
\def\makeacc@#1#2{\def#1{\mathaccent"\accentclass@#2 }}
\makeacc@\hat{05E}
\makeacc@\check{014}
\makeacc@\tilde{07E}
\makeacc@\acute{013}
\makeacc@\grave{012}
\makeacc@\dot{05F}
\makeacc@\ddot{07F}
\makeacc@\breve{015}
\makeacc@\bar{016}

\newcount\skewcharcount@
\newcount\familycount@
\def\theskewchar@{\familycount@\@ne
 \global\skewcharcount@\the\skewchar\textfont\@ne                           
 \ifnum\fam>\m@ne\ifnum\fam<16
  \global\familycount@\the\fam\relax
  \global\skewcharcount@\the\skewchar\textfont\the\fam\relax\fi\fi          
 \ifnum\skewcharcount@>\m@ne
  \ifnum\skewcharcount@<128
  \multiply\familycount@256
  \global\advance\skewcharcount@\familycount@
  \global\advance\skewcharcount@28672
  \mathchar\skewcharcount@\else
  \global\skewcharcount@\m@ne\fi\else
 \global\skewcharcount@\m@ne\fi}                                            
\newcount\pointcount@
\def\getpoints@#1.#2\getpoints@{\pointcount@#1 }
\newdimen\accentdimen@
\newcount\accentmu@
\def\dimentomu@{\multiply\accentdimen@ 100
 \expandafter\getpoints@\the\accentdimen@\getpoints@
 \multiply\pointcount@18
 \divide\pointcount@\@m
 \global\accentmu@\pointcount@}
\def\Makeacc@#1#2{\def#1{\RIfM@\DN@{\mathaccent@
 {"\accentclass@#2 }}\else\DN@{\nonmatherr@{#1}}\fi\next@}}
\def\unbracefonts@{\let\Cal@\Cal@@\let\roman@\roman@@\let\bold@\bold@@
 \let\slanted@\slanted@@}
\def\mathaccent@#1#2{\ifnum\fam=\m@ne\xdef\thefam@{1}\else
 \xdef\thefam@{\the\fam}\fi                                                 
 \accentdimen@\z@                                                           
 \setboxz@h{\unbracefonts@$\m@th\fam\thefam@\relax#2$}
 \ifdim\accentdimen@=\z@\DN@{\mathaccent#1{#2}}
  \setbox@ne\hbox{\unbracefonts@$\m@th\fam\thefam@\relax#2\theskewchar@$}
  \setbox\tw@\hbox{$\m@th\ifnum\skewcharcount@=\m@ne\else
   \mathchar\skewcharcount@\fi$}
  \global\accentdimen@\wd@ne\global\advance\accentdimen@-\wdz@
  \global\advance\accentdimen@-\wd\tw@                                     
  \global\multiply\accentdimen@\tw@
  \dimentomu@\global\advance\accentmu@\@ne                                 
 \else\DN@{{\mathaccent#1{#2\mkern\accentmu@ mu}%
    \mkern-\accentmu@ mu}{}}\fi                                             
 \next@}\Makeacc@\Hat{05E}
\Makeacc@\Check{014}
\Makeacc@\Tilde{07E}
\Makeacc@\Acute{013}
\Makeacc@\Grave{012}
\Makeacc@\Dot{05F}
\Makeacc@\Ddot{07F}
\Makeacc@\Breve{015}
\Makeacc@\Bar{016}
\def\Vec{\RIfM@\DN@{\mathaccent@{"017E }}\else
 \DN@{\nonmatherr@\Vec}\fi\next@}
\def\newbox@{\alloc@4\box\chardef\insc@unt}
\def\accentedsymbol#1#2{\expandafter\newbox@\csname\expandafter
  \eat@\string#1@box\endcsname
 \expandafter\setbox\csname\expandafter\eat@
  \string#1@box\endcsname\hbox{$\m@th#2$}\define
  #1{\expandafter\copy\csname\expandafter\eat@\string#1@box\endcsname{}}}
\def\sqrt#1{\radical"270370 {#1}}
\let\underline@\underline
\let\overline@\overline
\def\underline#1{\underline@{#1}}
\def\overline#1{\overline@{#1}}
\Invalid@\leftroot
\Invalid@\uproot
\newcount\uproot@
\newcount\leftroot@
\def\root{\relaxnext@
  \DN@{\ifx\next\uproot\let\next@\nextii@\else
   \ifx\next\leftroot\let\next@\nextiii@\else
   \let\next@\plainroot@\fi\fi\next@}%
  \DNii@\uproot##1{\uproot@##1\relax\FN@\nextiv@}%
  \def\nextiv@{\ifx\next\space@\DN@. {\FN@\nextv@}\else
   \DN@.{\FN@\nextv@}\fi\next@.}%
  \def\nextv@{\ifx\next\leftroot\let\next@\nextvi@\else
   \let\next@\plainroot@\fi\next@}%
  \def\nextvi@\leftroot##1{\leftroot@##1\relax\plainroot@}%
   \def\nextiii@\leftroot##1{\leftroot@##1\relax\FN@\nextvii@}%
  \def\nextvii@{\ifx\next\space@
   \DN@. {\FN@\nextviii@}\else
   \DN@.{\FN@\nextviii@}\fi\next@.}%
  \def\nextviii@{\ifx\next\uproot\let\next@\nextix@\else
   \let\next@\plainroot@\fi\next@}%
  \def\nextix@\uproot##1{\uproot@##1\relax\plainroot@}%
  \bgroup\uproot@\z@\leftroot@\z@\FN@\next@}
\def\plainroot@#1\of#2{\setbox\rootbox\hbox{$\m@th\scriptscriptstyle{#1}$}%
 \mathchoice{\r@@t\displaystyle{#2}}{\r@@t\textstyle{#2}}
 {\r@@t\scriptstyle{#2}}{\r@@t\scriptscriptstyle{#2}}\egroup}
\def\r@@t#1#2{\setboxz@h{$\m@th#1\sqrt{#2}$}%
 \dimen@\ht\z@\advance\dimen@-\dp\z@
 \setbox@ne\hbox{$\m@th#1\mskip\uproot@ mu$}\advance\dimen@ by1.667\wd@ne
 \mkern-\leftroot@ mu\mkern5mu\raise.6\dimen@\copy\rootbox
 \mkern-10mu\mkern\leftroot@ mu\boxz@}
\def\boxed#1{\setboxz@h{$\m@th\displaystyle{#1}$}\dimen@.4\ex@
 \advance\dimen@3\ex@\advance\dimen@\dp\z@
 \hbox{\lower\dimen@\hbox{%
 \vbox{\hrule height.4\ex@
 \hbox{\vrule width.4\ex@\hskip3\ex@\vbox{\vskip3\ex@\boxz@\vskip3\ex@}%
 \hskip3\ex@\vrule width.4\ex@}\hrule height.4\ex@}%
 }}}
\let\ampersand@\relax
\newdimen\minaw@
\minaw@11.11128\ex@
\newdimen\minCDaw@
\minCDaw@2.5pc
\def\minCDarrowwidth#1{\RIfMIfI@\onlydmatherr@\minCDarrowwidth
 \else\minCDaw@#1\relax\fi\else\onlydmatherr@\minCDarrowwidth\fi}
\newif\ifCD@
\def\CD{\bgroup\vspace@\relax\let\ampersand@&\iffalse}\fi
 \CD@true\vcenter\bgroup\Let@\tabskip\z@skip\baselineskip20\ex@
 \lineskip3\ex@\lineskiplimit3\ex@\halign\bgroup
 &\hfill$\m@th##$\hfill\crcr}
\def\endCD{\crcr\egroup\egroup\egroup}
\newdimen\bigaw@
\atdef@>#1>#2>{\ampersand@                                                  
 \setboxz@h{$\m@th\ssize\;{#1}\;\;$}
 \setbox@ne\hbox{$\m@th\ssize\;{#2}\;\;$}
 \setbox\tw@\hbox{$\m@th#2$}
 \ifCD@\global\bigaw@\minCDaw@\else\global\bigaw@\minaw@\fi                 
 \ifdim\wdz@>\bigaw@\global\bigaw@\wdz@\fi
 \ifdim\wd@ne>\bigaw@\global\bigaw@\wd@ne\fi                                
 \ifCD@\hskip.5em\fi                                                        
 \ifdim\wd\tw@>\z@
  \mathrel{\mathop{\hbox to\bigaw@{\rightarrowfill}}\limits^{#1}_{#2}}
 \else\mathrel{\mathop{\hbox to\bigaw@{\rightarrowfill}}\limits^{#1}}\fi    
 \ifCD@\hskip.5em\fi                                                       
 \ampersand@}                                                              
\atdef@<#1<#2<{\ampersand@\setboxz@h{$\m@th\ssize\;\;{#1}\;$}%
 \setbox@ne\hbox{$\m@th\ssize\;\;{#2}\;$}\setbox\tw@\hbox{$\m@th#2$}%
 \ifCD@\global\bigaw@\minCDaw@\else\global\bigaw@\minaw@\fi
 \ifdim\wdz@>\bigaw@\global\bigaw@\wdz@\fi
 \ifdim\wd@ne>\bigaw@\global\bigaw@\wd@ne\fi
 \ifCD@\hskip.5em\fi
 \ifdim\wd\tw@>\z@
  \mathrel{\mathop{\hbox to\bigaw@{\leftarrowfill}}\limits^{#1}_{#2}}\else
  \mathrel{\mathop{\hbox to\bigaw@{\leftarrowfill}}\limits^{#1}}\fi
 \ifCD@\hskip.5em\fi\ampersand@}
\atdef@)#1)#2){\ampersand@
 \setboxz@h{$\m@th\ssize\;{#1}\;\;$}%
 \setbox@ne\hbox{$\m@th\ssize\;{#2}\;\;$}%
 \setbox\tw@\hbox{$\m@th#2$}%
 \ifCD@
 \global\bigaw@\minCDaw@\else\global\bigaw@\minaw@\fi
 \ifdim\wdz@>\bigaw@\global\bigaw@\wdz@\fi
 \ifdim\wd@ne>\bigaw@\global\bigaw@\wd@ne\fi
 \ifCD@\hskip.5em\fi
 \ifdim\wd\tw@>\z@
  \mathrel{\mathop{\hbox to\bigaw@{\rightarrowfill}}\limits^{#1}_{#2}}%
 \else\mathrel{\mathop{\hbox to\bigaw@{\rightarrowfill}}\limits^{#1}}\fi
 \ifCD@\hskip.5em\fi
 \ampersand@}
\atdef@(#1(#2({\ampersand@\setboxz@h{$\m@th\ssize\;\;{#1}\;$}%
 \setbox@ne\hbox{$\m@th\ssize\;\;{#2}\;$}\setbox\tw@\hbox{$\m@th#2$}%
 \ifCD@\global\bigaw@\minCDaw@\else\global\bigaw@\minaw@\fi
 \ifdim\wdz@>\bigaw@\global\bigaw@\wdz@\fi
 \ifdim\wd@ne>\bigaw@\global\bigaw@\wd@ne\fi
 \ifCD@\hskip.5em\fi
 \ifdim\wd\tw@>\z@
  \mathrel{\mathop{\hbox to\bigaw@{\leftarrowfill}}\limits^{#1}_{#2}}\else
  \mathrel{\mathop{\hbox to\bigaw@{\leftarrowfill}}\limits^{#1}}\fi
 \ifCD@\hskip.5em\fi\ampersand@}
\atdef@ A#1A#2A{\llap{$\m@th\vcenter{\hbox
 {$\ssize#1$}}$}\Big\uparrow\rlap{$\m@th\vcenter{\hbox{$\ssize#2$}}$}&&}
\atdef@ V#1V#2V{\llap{$\m@th\vcenter{\hbox
 {$\ssize#1$}}$}\Big\downarrow\rlap{$\m@th\vcenter{\hbox{$\ssize#2$}}$}&&}
\atdef@={&\hskip.5em\mathrel
 {\vbox{\hrule width\minCDaw@\vskip3\ex@\hrule width
 \minCDaw@}}\hskip.5em&}
\atdef@|{\Big\Vert&&}
\atdef@@\vert{\Big\Vert&&}
\def\pretend#1\haswidth#2{\setboxz@h{$\m@th\scriptstyle{#2}$}\hbox
 to\wdz@{\hfill$\m@th\scriptstyle{#1}$\hfill}}
\def\pmb{\RIfM@\expandafter\mathpalette\expandafter\pmb@\else
 \expandafter\pmb@@\fi}
\def\pmb@@#1{\leavevmode\setboxz@h{#1}\kern-.025em\copy\z@\kern-\wdz@
 \kern-.05em\copy\z@\kern-\wdz@\kern-.025em\raise.0433em\boxz@}
\def\binrel@@#1{\ifdim\wd2<\z@\mathbin{#1}\else\ifdim\wd\tw@>\z@
 \mathrel{#1}\else{#1}\fi\fi}
\newdimen\pmbraise@
\def\pmb@#1#2{\setbox\thr@@\hbox{$\m@th#1{#2}$}%
 \setbox4 \hbox{$\m@th#1\mkern.7794mu$}\pmbraise@\wd4
 \binrel@{#2}\binrel@@{\mkern-.45mu\copy\thr@@\kern-\wd\thr@@
 \mkern-.9mu\copy\thr@@\kern-\wd\thr@@\mkern-.45mu\raise\pmbraise@\box\thr@@}}
\def\documentstyle#1{\input #1.sty\relax}
\font\dummyft@=dummy
\fontdimen1 \dummyft@=\z@
\fontdimen2 \dummyft@=\z@
\fontdimen3 \dummyft@=\z@
\fontdimen4 \dummyft@=\z@
\fontdimen5 \dummyft@=\z@
\fontdimen6 \dummyft@=\z@
\fontdimen7 \dummyft@=\z@
\fontdimen8 \dummyft@=\z@
\fontdimen9 \dummyft@=\z@
\fontdimen10 \dummyft@=\z@
\fontdimen11 \dummyft@=\z@
\fontdimen12 \dummyft@=\z@
\fontdimen13 \dummyft@=\z@
\fontdimen14 \dummyft@=\z@
\fontdimen15 \dummyft@=\z@
\fontdimen16 \dummyft@=\z@
\fontdimen17 \dummyft@=\z@
\fontdimen18 \dummyft@=\z@
\fontdimen19 \dummyft@=\z@
\fontdimen20 \dummyft@=\z@
\fontdimen21 \dummyft@=\z@
\fontdimen22 \dummyft@=\z@
\def\fontlist@{\\{\tenrm}\\{\sevenrm}\\{\fiverm}\\{\teni}\\{\seveni}%
 \\{\fivei}\\{\tensy}\\{\sevensy}\\{\fivesy}\\{\tenex}\\{\tenbf}\\{\sevenbf}%
 \\{\fivebf}\\{\tensl}\\{\tenit}}
\def\font@#1=#2 {\rightappend@#1\to\fontlist@\font#1=#2 }
\def\dodummy@{{\def\\##1{\global\let##1\dummyft@}\fontlist@}}
\def\nopages@{\output={\setbox\z@\box255 \deadcycles\z@}%
 \alloc@5\toks\toksdef\@cclvi\output}
\let\galleys\nopages@
\newif\ifsyntax@
\newcount\countxviii@
\def\syntax{\syntax@true\dodummy@\countxviii@\count18
 \loop\ifnum\countxviii@>\m@ne\textfont\countxviii@=\dummyft@
 \scriptfont\countxviii@=\dummyft@\scriptscriptfont\countxviii@=\dummyft@
 \advance\countxviii@\m@ne\repeat                                           
 \dummyft@\tracinglostchars\z@\nopages@\frenchspacing\hbadness\@M}
\def\S@{S } \def\G@{G } \def\P@{P }
\newif\ifbadans@
\def\printoptions{\W@{Do you want S(yntax check),
  G(alleys) or P(ages)?^^JType S, G or P, follow by <return>: }\loop
 \read\m@ne to\ans@
 \xdef\next@{\def\noexpand\Ans@{\ans@}}\uppercase\expandafter{\next@}
 \ifx\Ans@\S@\badans@false\syntax\else
 \ifx\Ans@\G@\badans@false\galleys\else
 \ifx\Ans@\P@\badans@false\else
 \badans@true\fi\fi\fi
 \ifbadans@\W@{Type S, G or P, follow by <return>: }%
 \repeat}
\def\alloc@#1#2#3#4#5{\global\advance\count1#1by\@ne
 \ch@ck#1#4#2\allocationnumber=\count1#1
 \global#3#5=\allocationnumber
 \ifalloc@\wlog{\string#5=\string#2\the\allocationnumber}\fi}
\def\document{\def\alloclist@{}\def\fontlist@{}}
\let\enddocument\bye

\let\proclaim\undefined
\let\footnote\undefined
\let\=\undefined
\let\>\undefined

\catcode`\@=\active

%
%
\def\next{AMSPPT}\ifx\styname\next \endinput\fi
\catcode`\@=11
\def\styname{AMSPPT}
\def\styversion{2.0}
{\W@{\styname.STY - Version \styversion}\W@{}}
\hyphenation{acad-e-my acad-e-mies af-ter-thought anom-aly anom-alies
an-ti-deriv-a-tive an-tin-o-my an-tin-o-mies apoth-e-o-ses apoth-e-o-sis
ap-pen-dix ar-che-typ-al as-sign-a-ble as-sist-ant-ship as-ymp-tot-ic
asyn-chro-nous at-trib-uted at-trib-ut-able bank-rupt bank-rupt-cy
bi-dif-fer-en-tial blue-print busier busiest cat-a-stroph-ic
cat-a-stroph-i-cally con-gress cross-hatched data-base de-fin-i-tive
de-riv-a-tive dis-trib-ute dri-ver dri-vers eco-nom-ics econ-o-mist
elit-ist equi-vari-ant ex-quis-ite ex-tra-or-di-nary flow-chart
for-mi-da-ble forth-right friv-o-lous ge-o-des-ic ge-o-det-ic geo-met-ric
griev-ance griev-ous griev-ous-ly hexa-dec-i-mal ho-lo-no-my ho-mo-thetic
ideals idio-syn-crasy in-fin-ite-ly in-fin-i-tes-i-mal ir-rev-o-ca-ble
key-stroke lam-en-ta-ble light-weight mal-a-prop-ism man-u-script
mar-gin-al meta-bol-ic me-tab-o-lism meta-lan-guage me-trop-o-lis
met-ro-pol-i-tan mi-nut-est mol-e-cule mono-chrome mono-pole mo-nop-oly
mono-spline mo-not-o-nous mul-ti-fac-eted mul-ti-plic-able non-euclid-ean
non-iso-mor-phic non-smooth par-a-digm par-a-bol-ic pa-rab-o-loid
pa-ram-e-trize para-mount pen-ta-gon phe-nom-e-non post-script pre-am-ble
pro-ce-dur-al pro-hib-i-tive pro-hib-i-tive-ly pseu-do-dif-fer-en-tial
pseu-do-fi-nite pseu-do-nym qua-drat-ics quad-ra-ture qua-si-smooth
qua-si-sta-tion-ary qua-si-tri-an-gu-lar quin-tes-sence quin-tes-sen-tial
re-arrange-ment rec-tan-gle ret-ri-bu-tion retro-fit retro-fit-ted
right-eous right-eous-ness ro-bot ro-bot-ics sched-ul-ing se-mes-ter
semi-def-i-nite semi-ho-mo-thet-ic set-up se-vere-ly side-step sov-er-eign
spe-cious spher-oid spher-oid-al star-tling star-tling-ly
sta-tis-tics sto-chas-tic straight-est strange-ness strat-a-gem strong-hold
sum-ma-ble symp-to-matic syn-chro-nous topo-graph-i-cal tra-vers-a-ble
tra-ver-sal tra-ver-sals treach-ery turn-around un-at-tached un-err-ing-ly
white-space wide-spread wing-spread wretch-ed wretch-ed-ly Brown-ian
Eng-lish Euler-ian Feb-ru-ary Gauss-ian Grothen-dieck Hamil-ton-ian
Her-mit-ian Jan-u-ary Japan-ese Kor-te-weg Le-gendre Lip-schitz
Lip-schitz-ian Mar-kov-ian Noe-ther-ian No-vem-ber Rie-mann-ian
Schwarz-schild Sep-tem-ber}
\Invalid@\nofrills
\Invalid@\usualspace
\newif\ifnofrills@
\def\nofrills@#1#2{\relaxnext@
  \DN@{\ifx\next\nofrills
    \nofrills@true\let#2\relax\DN@\nofrills{\nextii@}%
  \else
    \nofrills@false\def#2{#1}\let\next@\nextii@\fi
\next@}}
\def\usualspace@#1{\ifnofrills@\def\usualspace{#1}\fi}
\def\addto#1#2{\csname \expandafter\eat@\string#1@\endcsname
  \expandafter{\the\csname \expandafter\eat@\string#1@\endcsname#2}}
\newdimen\bigsize@
\def\big@#1#2{{\hbox{$\left#2\vcenter to#1\bigsize@{}%
  \right.\nulldelimiterspace\z@\m@th$}}}
\def\big{\big@\@ne}
\def\Big{\big@{1.5}}
\def\bigg{\big@\tw@}
\def\Bigg{\big@{2.5}}
\def\raggedcenter@{\leftskip\z@ plus.4\hsize \rightskip\leftskip
 \parfillskip\z@ \parindent\z@ \spaceskip.3333em \xspaceskip.5em
 \pretolerance9999\tolerance9999 \exhyphenpenalty\@M
 \hyphenpenalty\@M \let\\\linebreak}
\def\upperspecialchars{\def\ss{SS}\let\i=I\let\j=J\let\ae\AE\let\oe\OE
  \let\o\O\let\aa\AA\let\l\L}
\def\uppercasetext@#1{%
  {\spaceskip1.2\fontdimen2\the\font plus1.2\fontdimen3\the\font
   \upperspecialchars\uctext@#1$\m@th\aftergroup\eat@$}}
\def\uctext@#1$#2${\endash@#1-\endash@$#2$\uctext@}
\def\endash@#1-#2\endash@{\uppercase{#1}\if\notempty{#2}--\endash@#2\endash@\fi}
\def\runaway@#1{\DN@{#1}\ifx\envir@\next@
  \Err@{You seem to have a missing or misspelled \string\end#1 ...}%
  \let\envir@\empty\fi}
\newif\iftemp@
\def\notempty#1{TT\fi\def\test@{#1}\ifx\test@\empty\temp@false
  \else\temp@true\fi \iftemp@}
\font@\tensmc=cmcsc10
\font@\sevenex=cmex7
\font@\sevenit=cmti7
\font@\eightrm=cmr8 
\font@\sixrm=cmr6 
\font@\eighti=cmmi8     \skewchar\eighti='177 
\font@\sixi=cmmi6       \skewchar\sixi='177   
\font@\eightsy=cmsy8    \skewchar\eightsy='60 
\font@\sixsy=cmsy6      \skewchar\sixsy='60   
\font@\eightex=cmex8
\font@\eightbf=cmbx8 
\font@\sixbf=cmbx6   
\font@\eightit=cmti8 
\font@\eightsl=cmsl8 
\font@\eightsmc=cmcsc8
\font@\eighttt=cmtt8 
\loadmsam
\loadmsbm
\loadeufm
\input amssym.tex\relax
\newtoks\tenpoint@
\def\tenpoint{\normalbaselineskip12\p@
 \abovedisplayskip12\p@ plus3\p@ minus9\p@
 \belowdisplayskip\abovedisplayskip
 \abovedisplayshortskip\z@ plus3\p@
 \belowdisplayshortskip7\p@ plus3\p@ minus4\p@
 \textonlyfont@\rm\tenrm \textonlyfont@\it\tenit
 \textonlyfont@\sl\tensl \textonlyfont@\bf\tenbf
 \textonlyfont@\smc\tensmc \textonlyfont@\tt\tentt
 \ifsyntax@ \def\big##1{{\hbox{$\left##1\right.$}}}%
  \let\Big\big \let\bigg\big \let\Bigg\big
 \else
  \textfont\z@=\tenrm  \scriptfont\z@=\sevenrm  \scriptscriptfont\z@=\fiverm
  \textfont\@ne=\teni  \scriptfont\@ne=\seveni  \scriptscriptfont\@ne=\fivei
  \textfont\tw@=\tensy \scriptfont\tw@=\sevensy \scriptscriptfont\tw@=\fivesy
  \textfont\thr@@=\tenex \scriptfont\thr@@=\sevenex
        \scriptscriptfont\thr@@=\sevenex
  \textfont\itfam=\tenit \scriptfont\itfam=\sevenit
        \scriptscriptfont\itfam=\sevenit
  \textfont\bffam=\tenbf \scriptfont\bffam=\sevenbf
        \scriptscriptfont\bffam=\fivebf
  \setbox\strutbox\hbox{\vrule height8.5\p@ depth3.5\p@ width\z@}%
  \setbox\strutbox@\hbox{\lower.5\normallineskiplimit\vbox{%
        \kern-\normallineskiplimit\copy\strutbox}}%
 \setbox\z@\vbox{\hbox{$($}\kern\z@}\bigsize@=1.2\ht\z@
 \fi
 \normalbaselines\rm\ex@.2326ex\jot3\ex@\the\tenpoint@}
\newtoks\eightpoint@
\def\eightpoint{\normalbaselineskip10\p@
 \abovedisplayskip10\p@ plus2.4\p@ minus7.2\p@
 \belowdisplayskip\abovedisplayskip
 \abovedisplayshortskip\z@ plus2.4\p@
 \belowdisplayshortskip5.6\p@ plus2.4\p@ minus3.2\p@
 \textonlyfont@\rm\eightrm \textonlyfont@\it\eightit
 \textonlyfont@\sl\eightsl \textonlyfont@\bf\eightbf
 \textonlyfont@\smc\eightsmc \textonlyfont@\tt\eighttt
 \ifsyntax@\def\big##1{{\hbox{$\left##1\right.$}}}%
  \let\Big\big \let\bigg\big \let\Bigg\big
 \else
  \textfont\z@=\eightrm \scriptfont\z@=\sixrm \scriptscriptfont\z@=\fiverm
  \textfont\@ne=\eighti \scriptfont\@ne=\sixi \scriptscriptfont\@ne=\fivei
  \textfont\tw@=\eightsy \scriptfont\tw@=\sixsy \scriptscriptfont\tw@=\fivesy
  \textfont\thr@@=\eightex \scriptfont\thr@@=\sevenex
   \scriptscriptfont\thr@@=\sevenex
  \textfont\itfam=\eightit \scriptfont\itfam=\sevenit
   \scriptscriptfont\itfam=\sevenit
  \textfont\bffam=\eightbf \scriptfont\bffam=\sixbf
   \scriptscriptfont\bffam=\fivebf
 \setbox\strutbox\hbox{\vrule height7\p@ depth3\p@ width\z@}%
 \setbox\strutbox@\hbox{\raise.5\normallineskiplimit\vbox{%
   \kern-\normallineskiplimit\copy\strutbox}}%
 \setbox\z@\vbox{\hbox{$($}\kern\z@}\bigsize@=1.2\ht\z@
 \fi
 \normalbaselines\eightrm\ex@.2326ex\jot3\ex@\the\eightpoint@}
\parindent1pc
\normallineskiplimit\p@
\newdimen\indenti \indenti=2pc
\def\pageheight#1{\vsize#1}
\def\pagewidth#1{\hsize#1%
   \captionwidth@\hsize \advance\captionwidth@-2\indenti}
\pagewidth{30pc} \pageheight{47pc}
\def\topmatter{%
 \ifx\undefined\msafam
 \else\font@\eightmsa=msam8 \font@\sixmsa=msam6
   \ifsyntax@\else \addto\tenpoint{\textfont\msafam=\tenmsa
              \scriptfont\msafam=\sevenmsa \scriptscriptfont\msafam=\fivemsa}%
     \addto\eightpoint{\textfont\msafam=\eightmsa \scriptfont\msafam=\sixmsa
              \scriptscriptfont\msafam=\fivemsa}%
   \fi
 \fi
 \ifx\undefined\msbfam
 \else\font@\eightmsb=msbm8 \font@\sixmsb=msbm6
   \ifsyntax@\else \addto\tenpoint{\textfont\msbfam=\tenmsb
         \scriptfont\msbfam=\sevenmsb \scriptscriptfont\msbfam=\fivemsb}%
     \addto\eightpoint{\textfont\msbfam=\eightmsb \scriptfont\msbfam=\sixmsb
         \scriptscriptfont\msbfam=\fivemsb}%
   \fi
 \fi
 \ifx\undefined\eufmfam
 \else \font@\eighteufm=eufm8 \font@\sixeufm=eufm6
   \ifsyntax@\else \addto\tenpoint{\textfont\eufmfam=\teneufm
       \scriptfont\eufmfam=\seveneufm \scriptscriptfont\eufmfam=\fiveeufm}%
     \addto\eightpoint{\textfont\eufmfam=\eighteufm
       \scriptfont\eufmfam=\sixeufm \scriptscriptfont\eufmfam=\fiveeufm}%
   \fi
 \fi
 \ifx\undefined\eufbfam
 \else \font@\eighteufb=eufb8 \font@\sixeufb=eufb6
   \ifsyntax@\else \addto\tenpoint{\textfont\eufbfam=\teneufb
      \scriptfont\eufbfam=\seveneufb \scriptscriptfont\eufbfam=\fiveeufb}%
    \addto\eightpoint{\textfont\eufbfam=\eighteufb
      \scriptfont\eufbfam=\sixeufb \scriptscriptfont\eufbfam=\fiveeufb}%
   \fi
 \fi
 \ifx\undefined\eusmfam
 \else \font@\eighteusm=eusm8 \font@\sixeusm=eusm6
   \ifsyntax@\else \addto\tenpoint{\textfont\eusmfam=\teneusm
       \scriptfont\eusmfam=\seveneusm \scriptscriptfont\eusmfam=\fiveeusm}%
     \addto\eightpoint{\textfont\eusmfam=\eighteusm
       \scriptfont\eusmfam=\sixeusm \scriptscriptfont\eusmfam=\fiveeusm}%
   \fi
 \fi
 \ifx\undefined\eusbfam
 \else \font@\eighteusb=eusb8 \font@\sixeusb=eusb6
   \ifsyntax@\else \addto\tenpoint{\textfont\eusbfam=\teneusb
       \scriptfont\eusbfam=\seveneusb \scriptscriptfont\eusbfam=\fiveeusb}%
     \addto\eightpoint{\textfont\eusbfam=\eighteusb
       \scriptfont\eusbfam=\sixeusb \scriptscriptfont\eusbfam=\fiveeusb}%
   \fi
 \fi
 \ifx\undefined\eurmfam
 \else \font@\eighteurm=eurm8 \font@\sixeurm=eurm6
   \ifsyntax@\else \addto\tenpoint{\textfont\eurmfam=\teneurm
       \scriptfont\eurmfam=\seveneurm \scriptscriptfont\eurmfam=\fiveeurm}%
     \addto\eightpoint{\textfont\eurmfam=\eighteurm
       \scriptfont\eurmfam=\sixeurm \scriptscriptfont\eurmfam=\fiveeurm}%
   \fi
 \fi
 \ifx\undefined\eurbfam
 \else \font@\eighteurb=eurb8 \font@\sixeurb=eurb6
   \ifsyntax@\else \addto\tenpoint{\textfont\eurbfam=\teneurb
       \scriptfont\eurbfam=\seveneurb \scriptscriptfont\eurbfam=\fiveeurb}%
    \addto\eightpoint{\textfont\eurbfam=\eighteurb
       \scriptfont\eurbfam=\sixeurb \scriptscriptfont\eurbfam=\fiveeurb}%
   \fi
 \fi
 \ifx\undefined\cmmibfam
 \else \font@\eightcmmib=cmmib8 \font@\sixcmmib=cmmib6
   \ifsyntax@\else \addto\tenpoint{\textfont\cmmibfam=\tencmmib
       \scriptfont\cmmibfam=\sevencmmib \scriptscriptfont\cmmibfam=\fivecmmib}%
    \addto\eightpoint{\textfont\cmmibfam=\eightcmmib
       \scriptfont\cmmibfam=\sixcmmib \scriptscriptfont\cmmibfam=\fivecmmib}%
   \fi
 \fi
 \ifx\undefined\cmbsyfam
 \else \font@\eightcmbsy=cmbsy8 \font@\sixcmbsy=cmbsy6
   \ifsyntax@\else \addto\tenpoint{\textfont\cmbsyfam=\tencmbsy
      \scriptfont\cmbsyfam=\sevencmbsy \scriptscriptfont\cmbsyfam=\fivecmbsy}%
    \addto\eightpoint{\textfont\cmbsyfam=\eightcmbsy
      \scriptfont\cmbsyfam=\sixcmbsy \scriptscriptfont\cmbsyfam=\fivecmbsy}%
   \fi
 \fi
 \let\topmatter\relax}
\def\chapterno@{\uppercase\expandafter{\romannumeral\chaptercount@}}
\newcount\chaptercount@
\def\chapter{\nofrills@{\afterassignment\chapterno@
                        CHAPTER \global\chaptercount@=}\chapter@
 \DNii@##1{\leavevmode\hskip-\leftskip
   \rlap{\vbox to\z@{\vss\centerline{\eightpoint
   \chapter@##1\unskip}\baselineskip2pc\null}}\hskip\leftskip
   \nofrills@false}%
 \FN@\next@}
\newbox\titlebox@
\def\title{\nofrills@{\uppercasetext@}\title@%
 \DNii@##1\endtitle{\global\setbox\titlebox@\vtop{\tenpoint\bf
 \raggedcenter@\ignorespaces
 \baselineskip1.3\baselineskip\title@{##1}\endgraf}%
 \ifmonograph@ \edef\next{\the\leftheadtoks}\ifx\next\empty
    \leftheadtext{##1}\fi
 \fi
 \edef\next{\the\rightheadtoks}\ifx\next\empty \rightheadtext{##1}\fi
 }\FN@\next@}
\newbox\authorbox@
\def\author#1\endauthor{\global\setbox\authorbox@
 \vbox{\tenpoint\smc\raggedcenter@\ignorespaces
 #1\endgraf}\relaxnext@ \edef\next{\the\leftheadtoks}%
 \ifx\next\empty\leftheadtext{#1}\fi}
\newbox\affilbox@
\def\affil#1\endaffil{\global\setbox\affilbox@
 \vbox{\tenpoint\raggedcenter@\ignorespaces#1\endgraf}}
\newcount\addresscount@
\addresscount@\z@
\def\address#1\endaddress{\global\advance\addresscount@\@ne
  \expandafter\gdef\csname address\number\addresscount@\endcsname
  {\vskip12\p@ minus6\p@\noindent\eightpoint\smc\ignorespaces#1\par}}
\def\email{\nofrills@{\eightpoint{\it E-mail\/}:\enspace}\email@
  \DNii@##1\endemail{%
  \expandafter\gdef\csname email\number\addresscount@\endcsname
  {\def\usualspace{{\it\enspace}}\smallskip\noindent\eightpoint\email@
  \ignorespaces##1\par}}%
 \FN@\next@}
\def\thedate@{}
\def\date#1\enddate{\gdef\thedate@{\tenpoint\ignorespaces#1\unskip}}
\def\thethanks@{}
\def\thanks#1\endthanks{\gdef\thethanks@{\eightpoint\ignorespaces#1.\unskip}}
\def\thekeywords@{}
\def\keywords{\nofrills@{{\it Key words and phrases.\enspace}}\keywords@
 \DNii@##1\endkeywords{\def\thekeywords@{\def\usualspace{{\it\enspace}}%
 \eightpoint\keywords@\ignorespaces##1\unskip.}}%
 \FN@\next@}
\def\thesubjclass@{}
\def\subjclass{\nofrills@{{\rm1980 {\it Mathematics Subject
   Classification\/} (1985 {\it Revision\/}).\enspace}}\subjclass@
 \DNii@##1\endsubjclass{\def\thesubjclass@{\def\usualspace
  {{\rm\enspace}}\eightpoint\subjclass@\ignorespaces##1\unskip.}}%
 \FN@\next@}
\newbox\abstractbox@
\def\abstract{\nofrills@{{\smc Abstract.\enspace}}\abstract@
 \DNii@{\setbox\abstractbox@\vbox\bgroup\noindent$$\vbox\bgroup
  \def\envir@{abstract}\advance\hsize-2\indenti
  \usualspace@{{\enspace}}\eightpoint \noindent\abstract@\ignorespaces}%
 \FN@\next@}
\def\endabstract{\par\unskip\egroup$$\egroup}
\def\widestnumber#1#2{\begingroup\let\head\null\let\subhead\empty
   \let\subsubhead\subhead
   \ifx#1\head\global\setbox\tocheadbox@\hbox{#2.\enspace}%
   \else\ifx#1\subhead\global\setbox\tocsubheadbox@\hbox{#2.\enspace}%
   \else\ifx#1\key\bgroup\let\endrefitem@\egroup
        \key#2\endrefitem@\global\refindentwd\wd\keybox@
   \else\ifx#1\no\bgroup\let\endrefitem@\egroup
        \no#2\endrefitem@\global\refindentwd\wd\nobox@
   \else\ifx#1\page\global\setbox\pagesbox@\hbox{\quad\bf#2}%
   \else\ifx#1\item\setboxz@h{#2}\global\rosteritemwd\wdz@
        \global\advance\rosteritemwd by.5\parindent
   \else\message{\string\widestnumber is not defined for this option
   (\string#1)}%
\fi\fi\fi\fi\fi\fi\endgroup}
\newif\ifmonograph@
\def\Monograph{\monograph@true \let\headmark\rightheadtext
  \let\varindent@\indent \def\headfont@{\bf}\def\proclaimfont@{\smc}%
  \def\demofont@{\smc}}
\let\varindent@\noindent
\newbox\tocheadbox@    \newbox\tocsubheadbox@
\newbox\tocbox@
\def\toc{\toc@{Contents}}
\def\newtocdefs{%
   \def \title##1\endtitle
       {\penaltyandskip@\z@\smallskipamount
        \hangindent\wd\tocheadbox@\noindent{\bf##1}}%
   \def \chapter##1{%
        Chapter \uppercase\expandafter{\romannumeral##1.\unskip}\enspace}%
   \def \specialhead##1\endspecialhead
       {\par\hangindent\wd\tocheadbox@ \noindent##1\par}%
   \def \head##1 ##2\endhead
       {\par\hangindent\wd\tocheadbox@ \noindent
        \if\notempty{##1}\hbox to\wd\tocheadbox@{\hfil##1\enspace}\fi
        ##2\par}%
   \def \subhead##1 ##2\endsubhead
       {\par\vskip-\parskip {\normalbaselines
        \advance\leftskip\wd\tocheadbox@
        \hangindent\wd\tocsubheadbox@ \noindent
        \if\notempty{##1}\hbox to\wd\tocsubheadbox@{##1\unskip\hfil}\fi
         ##2\par}}%
   \def \subsubhead##1 ##2\endsubsubhead
       {\par\vskip-\parskip {\normalbaselines
        \advance\leftskip\wd\tocheadbox@
        \hangindent\wd\tocsubheadbox@ \noindent
        \if\notempty{##1}\hbox to\wd\tocsubheadbox@{##1\unskip\hfil}\fi
        ##2\par}}}
\def\toc@#1{\relaxnext@
   \def\page##1%
       {\unskip\penalty0\null\hfil
        \rlap{\hbox to\wd\pagesbox@{\quad\hfil##1}}\hfilneg\penalty\@M}%
 \DN@{\ifx\next\nofrills\DN@\nofrills{\nextii@}%
      \else\DN@{\nextii@{{#1}}}\fi
      \next@}%
 \DNii@##1{%
\ifmonograph@\bgroup\else\setbox\tocbox@\vbox\bgroup
   \centerline{\headfont@\ignorespaces##1\unskip}\nobreak
   \vskip\belowheadskip \fi
   \setbox\tocheadbox@\hbox{0.\enspace}%
   \setbox\tocsubheadbox@\hbox{0.0.\enspace}%
   \leftskip\indenti \rightskip\leftskip
   \setbox\pagesbox@\hbox{\bf\quad000}\advance\rightskip\wd\pagesbox@
   \newtocdefs
 }%
 \FN@\next@}
\def\endtoc{\par\egroup}
\let\pretitle\relax
\let\preauthor\relax
\let\preaffil\relax
\let\predate\relax
\let\preabstract\relax
\let\prepaper\relax
\def\dedicatory #1\enddedicatory{\def\preabstract{{\medskip
  \eightpoint\it \raggedcenter@#1\endgraf}}}
\def\thetranslator@{}
\def\translator#1\endtranslator{\def\thetranslator@{\nobreak\medskip
 \line{\eightpoint\hfil Translated by \uppercase{#1}\qquad\qquad}\nobreak}}
\outer\def\endtopmatter{\runaway@{abstract}%
 \edef\next{\the\leftheadtoks}\ifx\next\empty
  \expandafter\leftheadtext\expandafter{\the\rightheadtoks}\fi
 \ifmonograph@\else
   \ifx\thesubjclass@\empty\else \makefootnote@{}{\thesubjclass@}\fi
   \ifx\thekeywords@\empty\else \makefootnote@{}{\thekeywords@}\fi
   \ifx\thethanks@\empty\else \makefootnote@{}{\thethanks@}\fi
 \fi
  \pretitle
  \ifmonograph@ \topskip7pc \else \topskip4pc \fi
  \box\titlebox@
  \topskip10pt
  \preauthor
  \ifvoid\authorbox@\else \vskip2.5pc plus1pc \unvbox\authorbox@\fi
  \preaffil
  \ifvoid\affilbox@\else \vskip1pc plus.5pc \unvbox\affilbox@\fi
  \predate
  \ifx\thedate@\empty\else \vskip1pc plus.5pc \line{\hfil\thedate@\hfil}\fi
  \preabstract
  \ifvoid\abstractbox@\else \vskip1.5pc plus.5pc \unvbox\abstractbox@ \fi
  \ifvoid\tocbox@\else\vskip1.5pc plus.5pc \unvbox\tocbox@\fi
  \prepaper
  \vskip2pc plus1pc
}
\def\document{\let\fontlist@\relax\let\alloclist@\relax
  \tenpoint}
\newskip\aboveheadskip       \aboveheadskip\bigskipamount
\newdimen\belowheadskip      \belowheadskip6\p@
\def\headfont@{\smc}
\def\penaltyandskip@#1#2{\relax\ifdim\lastskip<#2\relax\removelastskip
      \ifnum#1=\z@\else\penalty@#1\relax\fi\vskip#2%
  \else\ifnum#1=\z@\else\penalty@#1\relax\fi\fi}
\def\nobreak{\penalty\@M
  \ifvmode\def\penalty@{\let\penalty@\penalty\count@@@}%
  \everypar{\let\penalty@\penalty\everypar{}}\fi}
\let\penalty@\penalty
\def\heading#1\endheading{\head#1\endhead}
\def\subheading#1{\subhead#1\endsubhead}
\def\specialheadfont@{\bf}
\outer\def\specialhead{\par\penaltyandskip@{-200}\aboveheadskip
  \begingroup\interlinepenalty\@M\rightskip\z@ plus\hsize \let\\\linebreak
  \specialheadfont@\noindent\ignorespaces}
\def\endspecialhead{\par\endgroup\nobreak\vskip\belowheadskip}
\outer\def\head#1\endhead{\par\penaltyandskip@{-200}\aboveheadskip
 {\headfont@\raggedcenter@\interlinepenalty\@M
 \ignorespaces#1\endgraf}\nobreak
 \vskip\belowheadskip
 \headmark{#1}}
\let\headmark\eat@
\newskip\subheadskip       \subheadskip\medskipamount
\def\subheadfont@{\bf}
\outer\def\subhead{\nofrills@{.\enspace}\subhead@
 \DNii@##1\endsubhead{\par\penaltyandskip@{-100}\subheadskip
  \varindent@{\usualspace@{{\subheadfont@\enspace}}%
 \subheadfont@\ignorespaces##1\unskip\subhead@}\ignorespaces}%
 \FN@\next@}
\outer\def\subsubhead{\nofrills@{.\enspace}\subsubhead@
 \DNii@##1\endsubsubhead{\par\penaltyandskip@{-50}\medskipamount
      {\usualspace@{{\it\enspace}}%
  \it\ignorespaces##1\unskip\subsubhead@}\ignorespaces}%
 \FN@\next@}
\def\proclaimheadfont@{\bf}
\outer\def\proclaim{\runaway@{proclaim}\def\envir@{proclaim}%
  \nofrills@{.\enspace}\proclaim@
 \DNii@##1{\penaltyandskip@{-100}\medskipamount\varindent@
   \usualspace@{{\proclaimheadfont@\enspace}}\proclaimheadfont@
   \ignorespaces##1\unskip\proclaim@
  \sl\ignorespaces}%
 \FN@\next@}
\outer\def\endproclaim{\let\envir@\relax\par\rm
  \penaltyandskip@{55}\medskipamount}
\def\demoheadfont@{\it}
\def\demo{\runaway@{proclaim}\nofrills@{.\enspace}\demo@
     \DNii@##1{\par\penaltyandskip@\z@\medskipamount
  {\usualspace@{{\demoheadfont@\enspace}}%
  \varindent@\demoheadfont@\ignorespaces##1\unskip\demo@}\rm
  \ignorespaces}\FN@\next@}

\def\qed{\ifhmode\unskip\nobreak\fi\quad\ifmmode\square\else$\m@th\square$\fi}

\def\definition{\runaway@{proclaim}%
  \nofrills@{.\proclaimheadfont@\enspace}\definition@
        \DNii@##1{\penaltyandskip@{-100}\medskipamount
        {\usualspace@{{\proclaimheadfont@\enspace}}%
        \varindent@\proclaimheadfont@\ignorespaces##1\unskip\definition@}%
        \rm \ignorespaces}\FN@\next@}

\newdimen\rosteritemwd
\newcount\rostercount@
\newif\iffirstitem@
\let\plainitem@\item
\newtoks\everypartoks@
\def\par@{\everypartoks@\expandafter{\the\everypar}\everypar{}}
\def\roster{\edef\leftskip@{\leftskip\the\leftskip}%
 \relaxnext@
 \rostercount@\z@  
 \def\item{\FN@\rosteritem@}%
 \DN@{\ifx\next\runinitem\let\next@\nextii@\else
  \let\next@\nextiii@\fi\next@}%
 \DNii@\runinitem  
  {\unskip  
   \DN@{\ifx\next[\let\next@\nextii@\else
    \ifx\next"\let\next@\nextiii@\else\let\next@\nextiv@\fi\fi\next@}%
   \DNii@[####1]{\rostercount@####1\relax
    \enspace{\rm(\number\rostercount@)}~\ignorespaces}%
   \def\nextiii@"####1"{\enspace{\rm####1}~\ignorespaces}%
   \def\nextiv@{\enspace{\rm(1)}\rostercount@\@ne~}%
   \par@\firstitem@false  
   \FN@\next@}%
 \def\nextiii@{\par\par@  
  \penalty\@m\smallskip\vskip-\parskip
  \firstitem@true}%
 \FN@\next@}
\def\rosteritem@{\iffirstitem@\firstitem@false\else\par\vskip-\parskip\fi
 \leftskip3\parindent\noindent  
 \DNii@[##1]{\rostercount@##1\relax
  \llap{\hbox to2.5\parindent{\hss\rm(\number\rostercount@)}%
   \hskip.5\parindent}\ignorespaces}%
 \def\nextiii@"##1"{%
  \llap{\hbox to2.5\parindent{\hss\rm##1}\hskip.5\parindent}\ignorespaces}%
 \def\nextiv@{\advance\rostercount@\@ne
  \llap{\hbox to2.5\parindent{\hss\rm(\number\rostercount@)}%
   \hskip.5\parindent}}%
 \ifx\next[\let\next@\nextii@\else\ifx\next"\let\next@\nextiii@\else
  \let\next@\nextiv@\fi\fi\next@}

\newif\ifnextRunin@
\def\endroster{\relaxnext@
 \par\leftskip@  
 \penalty-50 \vskip-\parskip\smallskip  
 \DN@{\ifx\next\Runinitem\let\next@\relax
  \else\nextRunin@false\let\item\plainitem@  
   \ifx\next\par 
    \DN@\par{\everypar\expandafter{\the\everypartoks@}}%
   \else  
    \DN@{\noindent\everypar\expandafter{\the\everypartoks@}}%
  \fi\fi\next@}%
 \FN@\next@}
\newcount\rosterhangafter@
\def\Runinitem#1\roster\runinitem{\relaxnext@
 \rostercount@\z@ 
 \def\item{\FN@\rosteritem@}%
 \def\runinitem@{#1}%
 \DN@{\ifx\next[\let\next\nextii@\else\ifx\next"\let\next\nextiii@
  \else\let\next\nextiv@\fi\fi\next}%
 \DNii@[##1]{\rostercount@##1\relax
  \def\item@{{\rm(\number\rostercount@)}}\nextv@}%
 \def\nextiii@"##1"{\def\item@{{\rm##1}}\nextv@}%
 \def\nextiv@{\advance\rostercount@\@ne
  \def\item@{{\rm(\number\rostercount@)}}\nextv@}%
 \def\nextv@{\setbox\z@\vbox  
  {\ifnextRunin@\noindent\fi  
  \runinitem@\unskip\enspace\item@~\par  
  \global\rosterhangafter@\prevgraf}%
  \firstitem@false  
  \ifnextRunin@\else\par\fi  
  \hangafter\rosterhangafter@\hangindent3\parindent
  \ifnextRunin@\noindent\fi  
  \runinitem@\unskip\enspace 
  \item@~\ifnextRunin@\else\par@\fi  
  \nextRunin@true\ignorespaces}%
 \FN@\next@}
\def\footmarkform@#1{$\m@th^{#1}$}
\let\thefootnotemark\footmarkform@
\def\makefootnote@#1#2{\insert\footins
 {\interlinepenalty\interfootnotelinepenalty
 \eightpoint\splittopskip\ht\strutbox\splitmaxdepth\dp\strutbox
 \floatingpenalty\@MM\leftskip\z@\rightskip\z@\spaceskip\z@\xspaceskip\z@
 \leavevmode{#1}\footstrut\ignorespaces#2\unskip\lower\dp\strutbox
 \vbox to\dp\strutbox{}}}
\newcount\footmarkcount@
\footmarkcount@\z@
\def\footnotemark{\let\@sf\empty\relaxnext@
 \ifhmode\edef\@sf{\spacefactor\the\spacefactor}\/\fi
 \DN@{\ifx[\next\let\next@\nextii@\else
  \ifx"\next\let\next@\nextiii@\else
  \let\next@\nextiv@\fi\fi\next@}%
 \DNii@[##1]{\footmarkform@{##1}\@sf}%
 \def\nextiii@"##1"{{##1}\@sf}%
 \def\nextiv@{\iffirstchoice@\global\advance\footmarkcount@\@ne\fi
  \footmarkform@{\number\footmarkcount@}\@sf}%
 \FN@\next@}
\def\footnotetext{\relaxnext@
 \DN@{\ifx[\next\let\next@\nextii@\else
  \ifx"\next\let\next@\nextiii@\else
  \let\next@\nextiv@\fi\fi\next@}%
 \DNii@[##1]##2{\makefootnote@{\footmarkform@{##1}}{##2}}%
 \def\nextiii@"##1"##2{\makefootnote@{##1}{##2}}%
 \def\nextiv@##1{\makefootnote@{\footmarkform@{\number\footmarkcount@}}{##1}}%
 \FN@\next@}
\def\footnote{\let\@sf\empty\relaxnext@
 \ifhmode\edef\@sf{\spacefactor\the\spacefactor}\/\fi
 \DN@{\ifx[\next\let\next@\nextii@\else
  \ifx"\next\let\next@\nextiii@\else
  \let\next@\nextiv@\fi\fi\next@}%
 \DNii@[##1]##2{\footnotemark[##1]\footnotetext[##1]{##2}}%
 \def\nextiii@"##1"##2{\footnotemark"##1"\footnotetext"##1"{##2}}%
 \def\nextiv@##1{\footnotemark\footnotetext{##1}}%
 \FN@\next@}
\def\adjustfootnotemark#1{\advance\footmarkcount@#1\relax}
\def\footnoterule{\kern-3\p@
  \hrule width 5pc\kern 2.6\p@} 
\def\captionfont@{\smc}
\def\topcaption#1#2\endcaption{%
  {\dimen@\hsize \advance\dimen@-\captionwidth@
   \rm\raggedcenter@ \advance\leftskip.5\dimen@ \rightskip\leftskip
  {\captionfont@#1}%
  \if\notempty{#2}.\enspace\ignorespaces#2\fi
  \endgraf}\nobreak\bigskip}
\def\botcaption#1#2\endcaption{%
  \nobreak\bigskip
  \setboxz@h{\captionfont@#1\if\notempty{#2}.\enspace\rm#2\fi}%
  {\dimen@\hsize \advance\dimen@-\captionwidth@
   \leftskip.5\dimen@ \rightskip\leftskip
   \noindent \ifdim\wdz@>\captionwidth@ 
   \else\hfil\fi 
  {\captionfont@#1}\if\notempty{#2}.\enspace\rm#2\fi\endgraf}}
\def\@ins{\par\begingroup\def\vspace##1{\vskip##1\relax}%
  \def\captionwidth##1{\captionwidth@##1\relax}%
  \setbox\z@\vbox\bgroup} 
\def\block{\RIfMIfI@\nondmatherr@\block\fi
       \else\ifvmode\vskip\abovedisplayskip\noindent\fi
        $$\def\endblock{\par\egroup$$}\fi
  \vbox\bgroup\advance\hsize-2\indenti\noindent}
\def\endblock{\par\egroup}
\def\cite#1{{\rm[{\citefont@\m@th#1}]}}
\def\citefont@{\rm}
\def\refsfont@{\eightpoint}
\outer\def\Refs{\runaway@{proclaim}%
 \relaxnext@ \DN@{\ifx\next\nofrills\DN@\nofrills{\nextii@}\else
  \DN@{\nextii@{References}}\fi\next@}%
 \DNii@##1{\penaltyandskip@{-200}\aboveheadskip
  \line{\hfil\headfont@\ignorespaces##1\unskip\hfil}\nobreak
  \vskip\belowheadskip
  \begingroup\refsfont@\sfcode`.=\@m}%
 \FN@\next@}

\newbox\nobox@            \newbox\keybox@           \newbox\bybox@
\newbox\paperbox@         \newbox\paperinfobox@     \newbox\jourbox@
\newbox\volbox@           \newbox\issuebox@         \newbox\yrbox@
\newbox\pagesbox@         \newbox\bookbox@          \newbox\bookinfobox@
\newbox\publbox@          \newbox\publaddrbox@      \newbox\finalinfobox@
\newbox\edsbox@           \newbox\langbox@
\newif\iffirstref@        \newif\iflastref@
\newif\ifprevjour@        \newif\ifbook@            \newif\ifprevinbook@
\newif\ifquotes@          \newif\ifbookquotes@      \newif\ifpaperquotes@
\newdimen\bysamerulewd@
\setboxz@h{\refsfont@\kern3em}
\bysamerulewd@\wdz@
\newdimen\refindentwd
\setboxz@h{\refsfont@ 00. }
\refindentwd\wdz@
\outer\def\ref{\begingroup \noindent\hangindent\refindentwd
 \firstref@true \def\nofrills{\def\refkern@{\kern3sp}}%
 \ref@}
\def\ref@{\book@false \bgroup\let\endrefitem@\egroup \ignorespaces}
\def\moreref{\endrefitem@\endref@\firstref@false\ref@}%
\def\transl{\endrefitem@\endref@\firstref@false
  \book@false
  \prepunct@
  \setboxz@h\bgroup \aftergroup\unhbox\aftergroup\z@
    \def\endrefitem@{\unskip\refkern@\egroup}\ignorespaces}%
\def\emptyifempty@{\dimen@\wd\currbox@
  \advance\dimen@-\wd\z@ \advance\dimen@-.1\p@
  \ifdim\dimen@<\z@ \setbox\currbox@\copy\voidb@x \fi}
\let\refkern@\relax
\def\endrefitem@{\unskip\refkern@\egroup
  \setboxz@h{\refkern@}\emptyifempty@}\ignorespaces
\def\refdef@#1#2#3{\edef\next@{\noexpand\endrefitem@
  \let\noexpand\currbox@\csname\expandafter\eat@\string#1box@\endcsname
    \noexpand\setbox\noexpand\currbox@\hbox\bgroup}%
  \toks@\expandafter{\next@}%
  \if\notempty{#2#3}\toks@\expandafter{\the\toks@
  \def\endrefitem@{\unskip#3\refkern@\egroup
  \setboxz@h{#2#3\refkern@}\emptyifempty@}#2}\fi
  \toks@\expandafter{\the\toks@\ignorespaces}%
  \edef#1{\the\toks@}}
\refdef@\no{}{. }
\refdef@\key{[\m@th}{] }
\refdef@\by{}{}
\def\bysame{\by\hbox to\bysamerulewd@{\hrulefill}\thinspace
   \kern0sp}
\def\manyby{\message{\string\manyby is no longer necessary; \string\by
  can be used instead, starting with version 2.0 of \styname.STY}\by}
\refdef@\paper{\ifpaperquotes@``\fi\it}{}
\refdef@\paperinfo{}{}
\def\jour{\endrefitem@\let\currbox@\jourbox@
  \setbox\currbox@\hbox\bgroup
  \def\endrefitem@{\unskip\refkern@\egroup
    \setboxz@h{\refkern@}\emptyifempty@
    \ifvoid\jourbox@\else\prevjour@true\fi}%
\ignorespaces}
\refdef@\vol{\ifbook@\else\bf\fi}{}
\refdef@\issue{no. }{}
\refdef@\yr{}{}
\refdef@\pages{}{}
\def\page{\endrefitem@\def\pp@{\def\pp@{pp.~}p.~}\let\currbox@\pagesbox@
  \setbox\currbox@\hbox\bgroup\ignorespaces}
\def\pp@{pp.~}
\def\book{\endrefitem@ \let\currbox@\bookbox@
 \setbox\currbox@\hbox\bgroup\def\endrefitem@{\unskip\refkern@\egroup
  \setboxz@h{\ifbookquotes@``\fi}\emptyifempty@
  \ifvoid\bookbox@\else\book@true\fi}%
  \ifbookquotes@``\fi\it\ignorespaces}
\def\inbook{\endrefitem@
  \let\currbox@\bookbox@\setbox\currbox@\hbox\bgroup
  \def\endrefitem@{\unskip\refkern@\egroup
  \setboxz@h{\ifbookquotes@``\fi}\emptyifempty@
  \ifvoid\bookbox@\else\book@true\previnbook@true\fi}%
  \ifbookquotes@``\fi\ignorespaces}
\refdef@\eds{(}{, eds.)}
\def\ed{\endrefitem@\let\currbox@\edsbox@
 \setbox\currbox@\hbox\bgroup
 \def\endrefitem@{\unskip, ed.)\refkern@\egroup
  \setboxz@h{(, ed.)}\emptyifempty@}(\ignorespaces}
\refdef@\bookinfo{}{}
\refdef@\publ{}{}
\refdef@\publaddr{}{}
\refdef@\finalinfo{}{}
\refdef@\lang{(}{)}

\let\refdef@\relax 
\def\ppunbox@#1{\ifvoid#1\else\prepunct@\unhbox#1\fi}
\def\nocomma@#1{\ifvoid#1\else\changepunct@3\prepunct@\unhbox#1\fi}
\def\changepunct@#1{\ifnum\lastkern<3 \unkern\kern#1sp\fi}
\def\prepunct@{\count@\lastkern\unkern
  \ifnum\lastpenalty=0
    \let\penalty@\relax
  \else
    \edef\penalty@{\penalty\the\lastpenalty\relax}%
  \fi
  \unpenalty
  \let\refspace@\ \ifcase\count@,
\or;\or.\or 
  \or\let\refspace@\relax
  \else,\fi
  \ifquotes@''\quotes@false\fi \penalty@ \refspace@
}
\def\transferpenalty@#1{\dimen@\lastkern\unkern
  \ifnum\lastpenalty=0\unpenalty\let\penalty@\relax
  \else\edef\penalty@{\penalty\the\lastpenalty\relax}\unpenalty\fi
  #1\penalty@\kern\dimen@}
\def\endref{\endrefitem@\lastref@true\endref@
  \par\endgroup \prevjour@false \previnbook@false }
\def\endref@{%
\iffirstref@
  \ifvoid\nobox@\ifvoid\keybox@\indent\fi
  \else\hbox to\refindentwd{\hss\unhbox\nobox@}\fi
  \ifvoid\keybox@
  \else\ifdim\wd\keybox@>\refindentwd
         \box\keybox@
       \else\hbox to\refindentwd{\unhbox\keybox@\hfil}\fi\fi
  \kern4sp\ppunbox@\bybox@
\fi 
  \ifvoid\paperbox@
  \else\prepunct@\unhbox\paperbox@
    \ifpaperquotes@\quotes@true\fi\fi
  \ppunbox@\paperinfobox@
  \ifvoid\jourbox@
    \ifprevjour@ \nocomma@\volbox@
      \nocomma@\issuebox@
      \ifvoid\yrbox@\else\changepunct@3\prepunct@(\unhbox\yrbox@
        \transferpenalty@)\fi
      \ppunbox@\pagesbox@
    \fi 
  \else \prepunct@\unhbox\jourbox@
    \nocomma@\volbox@
    \nocomma@\issuebox@
    \ifvoid\yrbox@\else\changepunct@3\prepunct@(\unhbox\yrbox@
      \transferpenalty@)\fi
    \ppunbox@\pagesbox@
  \fi 
  \ifbook@\prepunct@\unhbox\bookbox@ \ifbookquotes@\quotes@true\fi \fi
  \nocomma@\edsbox@
  \ppunbox@\bookinfobox@
  \ifbook@\ifvoid\volbox@\else\prepunct@ vol.~\unhbox\volbox@
  \fi\fi
  \ppunbox@\publbox@ \ppunbox@\publaddrbox@
  \ifbook@ \ppunbox@\yrbox@
    \ifvoid\pagesbox@
    \else\prepunct@\pp@\unhbox\pagesbox@\fi
  \else
    \ifprevinbook@ \ppunbox@\yrbox@
      \ifvoid\pagesbox@\else\prepunct@\pp@\unhbox\pagesbox@\fi
    \fi \fi
  \ppunbox@\finalinfobox@
  \iflastref@
    \ifvoid\langbox@.\ifquotes@''\fi
    \else\changepunct@2\prepunct@\unhbox\langbox@\fi
  \else
    \ifvoid\langbox@\changepunct@1%
    \else\changepunct@3\prepunct@\unhbox\langbox@
      \changepunct@1\fi
  \fi
}
\outer\def\enddocument{%
 \runaway@{proclaim}%
\ifmonograph@ 
\else
 \nobreak
 \thetranslator@
 \count@\z@ \loop\ifnum\count@<\addresscount@\advance\count@\@ne
 \csname address\number\count@\endcsname
 \csname email\number\count@\endcsname
 \repeat
\fi
 \vfill\supereject\end}
\def\folio{{\foliofont@\ifnum\pageno<\z@ \romannumeral-\pageno
 \else\number\pageno \fi}}
\def\foliofont@{\eightrm}
\def\headlinefont@{\eightpoint}
\def\leftheadline{\rlap{\folio}\hfill \iftrue\topmark\fi \hfill}
\def\rightheadline{\hfill \expandafter
  \hfill \llap{\folio}}
\newtoks\leftheadtoks
\newtoks\rightheadtoks
\def\leftheadtext{\nofrills@{\uppercasetext@}\lht@
  \DNii@##1{\leftheadtoks\expandafter{\lht@{##1}}%
    \mark{\the\leftheadtoks\noexpand\else\the\rightheadtoks}
    \ifsyntax@\setboxz@h{\def\\{\unskip\space\ignorespaces}%
        \headlinefont@##1}\fi}%
  \FN@\next@}
\def\rightheadtext{\nofrills@{\uppercasetext@}\rht@
  \DNii@##1{\rightheadtoks\expandafter{\rht@{##1}}%
    \mark{\the\leftheadtoks\noexpand\else\the\rightheadtoks}%
    \ifsyntax@\setboxz@h{\def\\{\unskip\space\ignorespaces}%
        \headlinefont@##1}\fi}%
  \FN@\next@}
\headline={\def\chapter#1{\chapterno@. }%
  \def\\{\unskip\space\ignorespaces}\headlinefont@
  \ifodd\pageno \rightheadline \else \leftheadline\fi}
\def\NoRunningHeads{\global\runheads@false\global\let\headmark\eat@}

\def\logo@{\baselineskip2pc \hbox to\hsize{\hfil\eightpoint Typeset by
 \AmSTeX}}
\newif\iffirstpage@     \firstpage@true
\newif\ifrunheads@      \runheads@true
\output={\output@}
\def\output@{\shipout\vbox{%
 \iffirstpage@ \global\firstpage@false
  \pagebody \logo@ \makefootline%
 \else \ifrunheads@ \makeheadline \pagebody
       \else \pagebody \makefootline \fi
 \fi}%
 \advancepageno \ifnum\outputpenalty>-\@MM\else\dosupereject\fi}
\tenpoint
\catcode`\@=\active

 \NoBlackBoxes





 \define\zz{ {{\bold{Z}_2}}}


 \define\calc{\Cal C}

 \define\calq{\Cal Q}
 
 \define\calz{\Cal Z}



 \define\cycd#1#2{{\calz}_{#1}(#2)}
 \define\cyc#1#2{{\calz}^{#1}(#2)}
 \define\cych#1{{{\calz}^{#1}}}
 \define\cycp#1#2{{\calz}^{#1}(\bbp(#2))}
 \define\cyf#1#2{{\cyc{#1}{#2}}^{fix}}
 \define\cyfd#1#2{{\cycd{#1}{#2}}^{fix}}

 \define\crl#1#2{{\calz}_{\bbr}^{#1}{(#2)}}

 \define\crd#1#2{\widetilde{\calz}_{\bbr}^{#1}{(#2)}}

 \define\crld#1#2{{\calz}_{#1,\bbr}{(#2)}}
 \define\crdd#1#2{\widetilde{\calz}_{#1,{\bbr}}{(#2)}}

 \define\crlh#1{{\calz}_{\bbr}^{#1}}
 \define\crdh#1{{\widetilde{\calz}_{\bbr}^{#1}}}
 \define\cyav#1#2{{{\cyc{#1}{#2}}^{av}}}
 \define\cyavd#1#2{{\cycd{#1}{#2}}^{av}}

 \define\cyaa#1#2{{\cyc{#1}{#2}}^{-}}
 \define\cyaad#1#2{{\cycd{#1}{#2}}^{-}}

 \define\cyq#1#2{{\calq}^{#1}(#2)}
 \define\cyqd#1#2{{\calq}_{#1}(#2)}

 \define\cqt#1#2{{\calz}_{\bbh}^{#1}{(#2)}}
 \define\cqtav#1#2{{\calz}^{#1}{(#2)}^{av}}
 \define\cqtrd#1#2{\widetilde{\calz}_{\bbh}^{#1}{(#2)}}

 \define\cyct#1#2{{\calz}^{#1}(#2)_\zz}
 \define\cyft#1#2{{\cyc{#1}{#2}}^{fix}_\zz}
\define\cxg#1#2{G^{#1}_{\bbc}(#2)}
 \define\reg#1#2{G^{#1}_{\bbr}(#2)}

 \define\cyaat#1#2{{\cyc{#1}{#2}}^{-}_\zz}



 \define\fflag#1#2{{#1}={#1}_{#2} \supset {#1}_{{#2}-1} \supset
 \ldots \supset {#1}_{0} }
 
 \define\vect#1{ {\Cal{V}ect}_{#1}}


 \define\chv#1#2#3{{\calc}^{#1}_{#2}(#3)}
 \define\chvd#1#2#3{{\calc}_{#1,#2}(#3)}
 \define\chm#1#2{{\calc}_{#1}(#2)}



 \define\Claim#1{\subheading{Claim #1}}

 \define\xrightarrow#1{\overset{#1}\to{\rightarrow}}


\hfuzz1pc 


\define\bbn{\Bbb N}
\define\bbz{\Bbb Z}

\define\bbr{\Bbb R}
\define\bbc{\Bbb C}
\define\bbh{\Bbb H}
\define\bbp{\Bbb P}

\define\bw{\bold w}

\define\bc{\bold C}

\define\cc{\Cal C}

\define\ce{\Cal E}
\define\ch{\Cal H}

\define\cf{\Cal F}

\define\cs{\Cal S}
\define\cz{\Cal Z}
\define\co#1{\Cal O_{#1}}
\define\ct{\Cal T}


\define\a{\alpha}
\redefine\b{\beta}
\define\g{\gamma}
\redefine\d{\delta}
\define\r{\rho}
\define\s{\sigma}
\define\z{\zeta}
\define\x{\xi}
\define\la{\lambda}
\define\e{\epsilon}
\redefine\D{\Delta}
\define\G{\Gamma}

\define\p#1{{\bbp}^{#1}}

\define\equdef{\overset\text{def}\to=}

\define\blbx{\hfill  $\square$}
\redefine\qed{\blbx}

\define\pf{\subheading{Proof}}
\define\Lemma#1{\subheading{Lemma #1}}
\define\Theorem#1{\subheading{Theorem #1}}
\define\Prop#1{\subheading{Proposition #1}}
\define\Cor#1{\subheading{Corollary #1}}
\define\Note#1{\subheading{Note #1}}
\define\Def#1{\subheading{Definition #1}}
\define\Remark#1{\subheading{Remark #1}}
\define\Ex#1{\subheading{Example #1}}
\define\arr{\longrightarrow}

\redefine\Xi{X_{\infty}}

\define\jac#1#2{\left(\!\!\!\left(
\frac{\partial #1}{\partial #2}
\right)\!\!\!\right)}
\define\restrict#1{\left. #1 \right|_{t_{p+1} = \dots = t_n = 0}}

\define\SP#1#2{{\roman SP}^{#1}(#2)}

\define\coc#1#2#3{\cc^{#1}(#2;\, #3)}
\define\zoc#1#2#3{\cz^{#1}(#2;\, #3)}
\define\zyc#1#2#3{\cz^{#1}(#2 \times #3)}

\define\Div{\roman{ Div}}
\define\ar#1{\overset{#1}\to{\longrightarrow}}

\define\th#1#2{{\Bbb H}^{#1}(#2)}
\define\hth#1#2{\widehat{\Bbb H}^{#1}(#2)}


\define\bad#1#2{\cf_{#2}(#1)}

\define\pch#1{\bbp_{\bbc}(\bbh^{#1})}

\def\l{\ell}

\def\I#1#2{I_{#1, #2}}

\def\D{D}

\def\L{L}

\def\z2t{\text{$\bbz_2\ct$}}

\def\<{\left<}
\def\>{\right>}
\def\[{\left[}
\def\]{\right]}

\def\wt{\widetilde}

\def\SS{2}
\def\TT{3}
\def\AA{4}
\def\FF{5}
\def\CC{6}
\def\DD{7}
\def\EE{8}
\def\HH{9}
\def\II{10}
\def\JJ{11}
\def\KK{12}

\def\supp{\operatorname{supp \ }}

\redefine\and{\qquad\text{and}\qquad}

\def\dbar{\overline{\partial}}

\def\cn{\Cal N}

\def\th{^{\text{th}}} 
\def\for{\qquad\text{ for }}
\def\foralll{\qquad\text{ for all }}

\define\pr{\operatorname{pr}}

\def\HC{T}
\def\Y{\Omega}

\def\TP{C}
\def\TPO{{\TP}^o}
\def\TPF#1#2{{\bold \TP}_{#1}(#2)}
\def\bd{M}
\def\bdy{\d}

\document  

\ \vskip .3in

\centerline{\bf \titfont Boundaries of Varieties in Projective
Manifolds}

\bigskip  
\centerline{by}
\bigskip

\centerline{\bf  \aufont  F. Reese Harvey  \ and \ 
H. Blaine Lawson, Jr.\footnote{\footfont Research
 partially supported by the NSF}}
              
\vskip .4in
 
\centerline{\bf Abstract} \medskip
  \font\abstractfont=cmr10 at 10 pt

{{\parindent=.9in\narrower\abstractfont \noindent 
We give a characterization of the boundaries of holomorphic chains
in complex projective space.  This extends previous work of 
the authors  and complements results of Dolbeault and Henkin.

}}

\vskip .5in            
               
\centerline{Table of Contents}
\medskip   
       
\hskip 1in  1.        Introduction
 
\hskip 1in  2.      Statement of the Problem        
  
\hskip 1in  3.      Discussion of the Main Result  

\hskip 1in  4.      Newton Hierarchies

\hskip 1in  5.      Newton Hierarchies with Coefficients in a Vector
                  Bundle

\hskip 1in  6.      Vector Bundle-valued Currents

\hskip 1in  7.     The Projective $\dbar$-Problem

\hskip 1in  8.      Curves in $\bbp^2$
 
\hskip 1in  9.     Curves in $\bbp^n$
  
\hskip 1in \!\! 10.  Boundaries of Higher Dimension
  
\hskip 1in  \!\! 11.  The General Result in Terms of 
                     Moment Conditions
 
\hskip 1in  \!\! 12.  The Relation to Dolbeault-Henkin

\medskip

\vfill\eject

%
%

%

\noindent
{\bf \S 1.  Introduction.}  
Let $\bd\subset X$ be a compact oriented submanifold of dimension
$2p-1$ in a quasi-projective manifold $X$. The question  addressed here
is: under what conditions does $\bd$ bound a complex submanifold of
dimension $p$, or more generally, a holomorphic $p$-chain in $X$?


When $X$ is Stein, much is
known.  For example if $p>1$, a necessary and sufficient condition is
that $\bd$ be {\bf maximally complex}, i.e., that  $\dim_{\bbc}(T_x\bd\cap
iT_x\bd)=p-1$ for all $x\in \bd$ (see [HL$_1$]).   This represents a
geometric generalization of the Bochner-Hartogs extension theorem for CR
functions [B].  When $p=1$, a necessary and sufficient condition is that
$\int_\bd\omega =0$ for all holomorphic 1-forms $\omega \in
\Omega^1(\bbc^n)$ (see [HL$_1$]). This is called the {\bf moment
condition}.

For a connected curve $\g$ in $\bbc^n$ there is another, quite different
characterization: $\g$ bounds a 1-dimensional analytic
subvariety in $\bbc^n$ if and only if it is not  polynomially convex 
(i.e., the polynomials are not uniformly  dense in  $C(\g)$). When this
occurs, the subvariety with boundary $\g$ coincides with the polynomial
hull of $\g$.  This work goes back to classical results of J. Wermer
[W$_1$] for smooth curves with generalizations to  progressively less
regular curves by E. Bishop, H. Royden, G. Stolzenberg, H. Alexander and
M. Lawrence (see [L]).

 In $\bbc^n$ there is yet another characterization
of boundaries of holomorphic chains  in terms of the
positivity of linking numbers due to Alexander and Wermer [A], [AW$_2$],
[W$_2$].

For $X =\bbp^n-\bbp^{n-k}, \ k>1$
there are results in [HL$_2$] which extend those of [HL$_1$].  When $p>
k$, maximal complexity suffices, and when $p=k$ there is a moment condition
which is necessary and sufficient. (See the discussion in \S  2.)

However, the case where $p<k$, which corresponds essentially to taking
$X=\bbp^n$, is quite difficult.  There is beautiful work of Dolbeault and
Henkin [DH$_{1,2}$] which recasts the problem in terms of  certain
structure in a Cauchy-Radon tranform of $\bd$.  Nevertheless, results for
boundaries in $\bbp^n$ similar to those outlined above
have not been established. This is the first in a series of articles
devoted to that problem. In general terms what we prove here is
the following.\medskip
 
\noindent
{\bf The Main Result.}{\sl
Consider a maximally complex submanifold $\bd$
of dimension $2p-1$ in $\bbp^n$.  Fix a projection 
$\pi:\bbp^n - - -> \bbp^p$ which is  proper on $\bd$, and 
a ``base point'' $a\in \bbp^p-\pi(\bd)$. Then $\bd$ bounds a holomorphic
$p$-chain which is positive and $\ell$-sheeted over $a$ if and only if
$\bd$ satisfies a sequence  of (explicit and algorithmically
computable) non-linear moment conditions which depend on $\ell$.}

\medskip

 This generalizes  the  results in [HL$_{1,2}$] where $\ell=0$, and
complements  those in [DH$_{1,2}$]. An expanded discussion of the main
result, including a detailed exposition of the moment conditions, is
given in \S 3.   These moment conditions involve certain universal
Newton polynomials.  They  constitute  a complete set of  computable 
obstructions  to solving the central  problem. The main theorem and its
proof appear in \S\S 8--11.

In their paper [DH$_{1,2}$] Dolbeault and Henkin also fix a  projection
$\pi$ as above and consider certain naturally associated Cauchy
transforms $T_*(\bd, \pi)$.
They then vary the projection in a neighborhood $U$ of
$\pi$ on the Grassmannian.  Their result asserts that $\bd$ bounds a
holomorphic $p$-chain iff second derivatives of $T_*(\bd, \pi)$ can be
written as second derivatives of a linear combination of solutions of
the ``shock equation'' on $U$. This result is quite general, but it
leaves open the difficult problem of determining, for a given boundary 
$\bd$, whether $T_*(\bd, \pi)$ can be written in such a fashion.

One of the  major difficulties in establishing $\bbp^n$-versions of this
problem is the lack of uniqueness of the solution.  If $dT=\bd$, then
$d(T+C)=\bd$ for any algebraic $p$-cycle $C$ in $\bbp^n$. Our condition of
``positivity with $\ell$ sheets over a point'' vastly reduces the 
ambuguity in the solution and enables the derivation of the moment
conditions.  The Dolbeault-Henkin theorem  is much less explicit but
automatically incorporates the full order of non-uniqueness.

In pondering this problem one naturally wonders if boundaries 
of varieties in $\bbp^n$ can be characterized by the vanishing 
of the integrals of certain forms with singularities.
 For curves  this leads to the following question. Given a compact
irreducible complex analytic curve $C\subset \bbp^n$ with boundary $dC\neq
\emptyset$, does there exist a divisor $D$ in $\bbp^n$ such that $D\cap
C=\emptyset$? If this were always true then  a given real curve $\g
\subset \bbp^n$ would bound an analytic variety iff it satisfied the
standard moment condition in some affine chart. This would amount to the
condition that $\int_{\g} \omega =0$ for all rational 1-forms with poles
on a fixed divisor $D$. However, in  beautiful work which employs the
Rosenlicht generalized jacobians for singular curves, Bruno Fabr\'e
[Fa$_*$] has  shown that this is not the case. There exist many such
varieties $C\subset\bbp^2$ (proper subsets of algebraic curves!) which
have the property that they meet every divisor.

In considering this problem it is also natural to ask
if there exist  extensions of the work of Wermer
and others from the affine to the projective case.
The answer is yes.  In [HL$_3$] we introduce 
the concept of the {\bf projective hull} $\widehat \g$ of a subset
$\g\subset\bbp^n$. When $\g$ is a  smooth  curve which 
bounds an analytic variety $C\subset\bbp^n$, one has $C\subset\widehat \g$
and there is much evidence for the conjecture that at least for real
analytic $\gamma$ one has
$$
\widehat {\g}\ =\ \cases
\Sigma \qquad \text{ when $C$ is a subset of an 
irreducible algebraic curve } 
\Sigma \\
 C \qquad \text{ otherwise. } 
\endcases
$$

In fact it has been shown in complete generality that for any compact
subset $K\subset \bbp^n$, the set ${\widehat K}^- - K$ is {\sl 1-concave}
in the sense of  [DL].  In particular, on any domain in $\bbp^n-K$ where
the Hausdorff 2-measure of ${\widehat K}^- - K$ is finite, ${\widehat K}^-
- K$ is a 1-dimensional analytic subvariety.  Furthermore, for any finite
union of  smooth pluripolar curves  (e.g., real analytic curves)
$\gamma\subset \bbp^2$, the Hausdorff dimension of $\widehat{\gamma}$ is
$\leq 2$. These results ([HL$_3$]) make plausible certain specific
 analogues of the work of Alexander and Wermer [AW$_2$],
[W$_2$] for submanifolds of projective space. (See [HL$_4$].) 

Recall that the polynomial hull of a set $Y\subset \bbc^n$
gives a concrete realization of the Gelfand spectrum of the Banach algebra
obtained by taking the uniform closure of the the polynomials in $C(Y)$
(cf. [Ho] or [AW$_1$]). Interestingly there exists a {\bf 
projective Gelfand transformation} defined for graded Banach algebras,
which for subsets of $\bbp^n$ gives a concrete realization of the 
projective hull mentioned above. It is related to the usual Gelfand
Transform much as Grothendieck's $Proj$ is related to the spectrum of a
ring. All this is  introduced and discussed in [HL$_3$].

A special but quite interesting instance of the problem considered
here is  that of characterizing the boundary values of meromophic
mappings $F:\Omega\to X$ where $\Omega$ is a domain in $\bbc^p$
and $X$ is a projective variety. If $f:d\Omega\to X$  is a given 
function, consider its graph $\bd = \text{graph}(f) \subset
\Omega \times X\subset  \bbp^p\times\bbp^N\subset \bbp^{pN+p+N}$ and
look for a variety $V$ with  $dV = M$. Any such $V$, if it
exists, will lie in $\bbp^p\times X$ but in general it will not be the
graph of a function over $\Omega$. Nevertheless, there is a
relatively simple solution to this problem which is presented in 
[HL$_4$].

\vskip .3in

%
%

\def\Co{C_0}

\centerline{\bf \S \SS. Statement of the Problem.} \medskip
 
Suppose that $\bd$ is a compact oriented submanifold of
class $C^1$ and dimension $2p-1\geq 1$ in complex projective $n$-space
$\bbp^n$. This paper addresses the following:

\medskip\noindent
{\bf Main Problem.}  Under what conditions on $\bd$  does there exist
a  holomorphic $p$-chain in $\bbp^n$ with  boundary $\bd$? 
That is, under what conditions does there exist a holomorphic 
$p$-chain $\HC$ in $\bbp^n-\bd$ such that, as a current in $\bbp^n$, 
$\HC$ has finite mass and satisfies $d\HC=\bd$?

\medskip

Recall that a {\bf holomorphic $p$-chain} in a complex manifold $Y$ is a
current $\HC$ of dimension $2p$ which can be written as a locally finite
sum  $$
\HC=\sum_j n_j [V_j]
$$ 
where for each $j$, $n_j\in\bbz$ and $[V_j]$ is the current given
by integration over (the regular points of) a canonically oriented
irreducible analytic subvariety of dimension $p$ in $Y$. 

\medskip

There are two fundamental issues which we shall not discuss
here since they are well covered in the literature.

\medskip\noindent
{\bf Boundary Regularity.} It is shown in [HL$_1$], [H]  that if $d\HC=\bd$
as above, then one has boundary regularity almost everywhere in the
following sense.  There exists a compact subset $\Sigma\subset \bd$ of
$(2p-1)$-dimensional measure 0, such that if $U\subset \bd-\Sigma$ is an
open set where $\bd$ is of class $C^k$, $1\leq k \leq \omega$, then there is
an  open set $\wt U \subset X$ with $\wt U \cap \bd = U$ in which  $\HC$ is a
complex manifold with $C^k$-boundary $\bd$.

\medskip\noindent
{\bf Allowable Singularities on $\bd$.} The regularity conditions on $\bd$
can be weakened to allow ``scar sets''.  (See [H] 
for a complete discussion.)  
The solutions we shall discuss here will also exist for these more general
boundaries, however we will restrict attention to the smooth case.

\medskip\noindent
{\bf Some History.}  Whenever the submanifold $\bd$ lies in the complement
of a linear subspace of codimension $\leq p$, our problem is well
understood. To state the results we need two basic concepts.

\Def {\SS.1} The submanifold $\bd$ is said to be {\bf maximally complex}
if for each $a\in \bd$, 
$$
\dim_{\bbc}(T_a\bd\cap iT_a\bd)\ =\ p-1.
$$
This is equivalent to the condition that
$$
\bd = \bd_{p-1,p}+\bd_{p,p-1}
\tag{\SS.1}
$$
where  $\bd = \sum_{r+2=2p-1}\bd_{r,s}$ is the Dolbeault decomposition of
the current $\bd$. That is, the Dolbeault components $\bd_{r,s}=0$ for 
$|r-s|>1$. This condition is automatic when $p=1$.

\Def {\SS.2} The submanifold $\bd$ is said to satisfy the {\bf moment
condition } on an open neighborhood $\Omega\subset \bbp^n$ if
$$
\int_\bd \omega\ =\ 0
\qquad\text{for all }\ \ \omega \in \ce^{p,p-1}(\Omega) \text{ \ with \ }
\dbar \omega=0.  $$

\Def {\SS.3} Let $\bd$ be a $2p-1$-dimensional current
with compact support in a complex manifold $\Omega$.  Then $\bd$ is said to
{\bf bound a holomorphic $p$-chain in $\Omega$} if there exists a
holomorphic $p$-chain $\HC$ with finite mass in $\Omega-\supp(\bd)$ and 
$\supp(\HC)\subset\subset\Omega$  with the property that $d\HC=\bd$ as
currents in $\Omega$. 
By [HS] any such $\HC$ is unique modulo holomorphic chains with
compact support in $\Omega$

\Theorem {\SS.4. [HL$_{1,2}$]}
{\sl  Let $\bd\subset \bbp^n$ be a compact oriented 
$(2p-1)$-dimensional submanifold of
class $C^1$, and suppose $\bd$ is contained in 
$$
\Omega \ =\ \bbp^n\ -\ \bbp^{n-k}
$$
where $\bbp^{n-k}$ is some linear subspace of dimension $n-k$.
Then $\bd$ is the boundary of
a holomorphic $p$-chain in $\Omega$ if either

\roster
\item  \ \  $k<p$ and $\bd$ is maximally complex, or

\item \ \  $k=p$ and $\bd$ is maximally complex and satisfies the moment
condition in $\Omega$
\endroster}
\medskip

Note that when $k=1$, $\Omega \cong \bbc^n$ and the moment condition enters
only for curves ($p=1$). Maximally complex submanifolds of higher
dimension in $\bbc^n$ automatically bound holomorphic chains. 
As noted in \S 1, the  case of curves in $\bbc^n$ is alaso related to the
study of polynomial hulls and has a long history going back to fundamental
work of Wermer. \medskip

\Note{\SS.5} If $\bd$ bounds a holomorphic chain in $\bbp^n$, then $\bd$
bounds a holomorphic chain in $\Omega = \bbp^n-\bbp^{n-(p+1)}$ for almost
all linear subspaces $\bbp^{n-(p+1)}$.  Hence the cases $k>p$ not
considered in Theorem \SS.4 are essentially equivalent to the one being
studied here.

\medskip

In projective space maximal complexity does not guarantee that
$\bd$ bounds a holomorphic chain, even for $p>>1$. Consider the following
examples.
\Ex{\SS.6} Let $\bd$ be a smooth real hypersurface in a complex
$p$-dimensional submanifold $X^p\subset \bbp^n$, and assume that 
$\bd$ is not homologous to 0 in $X^p$.  Then $\bd $ cannot bound a
holomorphic chain in $\bbp^n$ for any such chain would necessarily have
support in $X^p$.

\Ex{\SS.7} Let $Y\subset \bbp^m$ be any projective variety and let
$C\subset \bbc^2$ be the graph of the function $w=\exp(z+\frac 1z)$.
Let $\bd$ be any real hypersurface in $Y\times C$ which is homologous to
$Y\times \g$ where $\g\subset C$ is the curve given by retricting the graph
to $|z|=1$. Embed $\bd$ into projective space by the composition
$ 
\bd \subset Y\times C \subset \bbp^m\times \bbc^2
\subset \bbp^m\times \bbp^2 \subset \bbp^{3m+2}
$
where the last map is the Segr\'e embedding.
Again $\bd $ cannot bound a
holomorphic chain in $\bbp^{3m+2}$ for any such chain  would necessarily
have support in $Y\times C$.

\medskip\noindent
{\bf The Non-uniqueness in the Projective Case.} A subtle difficulty of
our main problem is the non-uniqueness of solutions when they exist.
If $\HC$ is a holomorphic chain with boundary $\bd$ in $\bbp^n$, then
so is $\HC+\HC_0$,  where $\HC_0$ is any holomorphic chain with 
$d\HC_0  =0$, i.e., any algebraic cycle.

There are of course {\sl canonical} solutions.  When $\bd$ is connected,
there is a unique solution that is irreducible in $\bbp^n-\bd$.  More
generally there is a unique solution of least mass. However, these
``simpliest''solutions are hard to characterize with a process that
also detects solutions for complicated, highly non-connected boundaries.

\medskip\noindent
{\bf Structuring the Problem.} To make the general problem more
tractible one could narrow the focus of the question. 
Two natural ways of doing this are as follows.

\medskip\noindent
{\bf  1. Prescribe the homology class of the solution.} 
 We fix a homology class $u\in H_{2p}(\bbp^n, \bd;\,\bbz)$ such that 
$\bdy u = [\bd]$ where $\bdy:H_{2p}(\bbp^n, \bd;\,\bbz)
\to H_{2p-1}(\bd;\,\bbz)$ is the boundary map, and search for a holomorphic
chain in $u$.  There is a short exact sequence
$$
0\ \arr\ H_{2p}(\bbp^n;\,\bbz)\ \arr\ H_{2p}(\bbp^n, \bd;\,\bbz)
\ @>{\bdy}>>\ H_{2p-1}(\bd;\,\bbz)\ \arr\ 0
\tag{\SS.2}
$$
from which we conclude the following.

\Lemma{\SS.8}{\sl A class $u\in H_{2p}(\bbp^n, \bd;\,\bbz)$ with 
$\bdy u = [\bd]$ is completely determined by its intersection 
number $u\cdot \bbp^{n-p}$ with any linear subspace $\bbp^{n-p}
\subset \bbp^n-\bd$.}\medskip

As noted in \SS.5 our problem is  equivalent to looking
for solutions in open subsets of the form $\Omega =\bbp^n-\bbp^{n-p-1}$
which contain $\bd$. The  homomorphism $\ H_{2p-1}(\Omega,\bd;\,\bbz)\to 
\ H_{2p-1}(\bbp^n,\bd;\,\bbz)$ induced by inclusion, is an isomorphism by
general position arguments. Under the linear projection 
$$
\pi:\bbp^n-\bbp^{n-p-1}\ \arr\ \bbp^p
\tag{\SS.3}
$$
the class $\pi_*u$ corresponds to an integer-valued function $\Co$ with 
$d\Co = d \pi_* u = \pi_*\bd$.  In particular,
$$
\Co(a) \ =\ u\cdot \bbp^{n-p}_a\qquad\text{ where }\ \ 
\bbp^{n-p}_a\equiv \pi^{-1}(a)
\tag{\SS.4}
$$
 for any $a\in \bbp^p-\pi\bd$.

Now integer-valued functions $\Co$ with $d \Co = \pi_* \bd$ are
determined up to an integer constant.   Thus by Lemma \SS.8 we see that:
$$\aligned 
&\text{\sl Prescribing the homology class  $u$ is equivalent to }\\
&\text{\sl choosing
such a function $\Co$ for some projection $\pi$.}
\endaligned
\tag{\SS.5}
$$

\medskip\noindent
{\bf   2. Require the solution $\HC$ to be positive.} We could
restrict the question by  asking only for {\bf positive} holomorphic
chains: $\HC=\sum_in_i[V_i]$ with $n_i\in \bbz^+$, having $d\HC=\bd$.
This approach has several features: 
\roster
\item "$\bullet$"   It covers the case where $\bd$ is connected,  which is 
already of considerable interest.

\item "$\bullet$"  If one also prescribes the homology class of the solution,
then there exists at most a {\sl finite-dimensional} space of solutions.

\item"$\bullet$" The integer-valued function $\Co$ discussed above becomes
the {\sl sheeting number} of $\HC$ over $\bbp^p$ for the projection $\pi$.

\item"$\bullet$" Positive solutions in a homology class $u$ are exactly the
solutions to the Plateau Problem for $\bd$ in $u$ for any Kaehler
metric on $\bbp^n$.
\endroster

We shall impose somewhat weaker restrictions.

\medskip\noindent
{\bf Our Approach:  Prescribe the homology class of the solution  and
require positivity  only over one point for some projection $\pi$.} 
We shall find that this approach
\roster
\item "$\bullet$"    covers the case where $\bd$ is connected,

\item "$\bullet$"   still guarantees finite-dimensionality of the solution
space,

\item"$\bullet$"   recaptures Theorem \SS.4 as a special case, and 

\item"$\bullet$"   leads to explicit non-linear moment conditions necessary
and sufficient for the solution.
\endroster

In fact when solutions to this problem exist, one also captures   the
canonical solution (of smallest mass)  by minimizing the sheeting number
over the given point.

  \vskip .3in

%
%

\centerline{\bf \S \TT. Discussion of the Main Result.} \medskip
 
Before launching a complete exposition we present the main ideas
in outline form. For clarity of exposition we focus here  on the
hypersurface case.  However, the result in general codimensions
 is very much the same (see \S\S \HH\ and \II).

Fix a compact oriented $C^1$-submanifold
$\bd\subset\bbp^{n+1}$ of dimension $2n-1\geq 1$ which is maximally
complex. Fix a point $\bbp^{0} \notin \bd$.   Choose homogeneous
coordinates $[z_0:\dots:z_{n+1}]$ for $\bbp^{n+1}$ 
with $\bbp^0 = [0:\dots:0:1]$ and let
$$
\Omega\equiv\bbp^{n+1}-\bbp^{0}\ @>{\pi}>>\ \bbp^n
$$
be the projection defined by
$\pi([z_0:\dots:z_{n+1}])=[z_0:\dots:z_n]$. Let $w$ be the tautological
cross-section of $\pi^*\co{\bbp^n}(1) \cong
\co{\bbp^{n+1}}(1)\bigl|_{\Omega}$.  Then each push-forward 
$\pi_*(w^d\bd)$ with $d\geq 0 $ is a well-defined current of degree one
with values in the bundle $\co{\bbp^n}(d)$. Let $\pi_*(w^d\bd)^{0,1}$
denote the  Dolbeault component of bidegree $(0,1)$ of this current.
Then the maximal complexity of $\bd$ (cf. (\SS.1)) and the fact that 
$d\bd=0$ imply that 
$$
\dbar \pi_*(w^d\bd)^{0,1}\ =\ 0\qquad\text{ for } d=0,1,2,...
$$

Now fix a point $a\in \bbp^n-\pi(\bd)$.  Then standard arguments show that
there exists a unique generalized section $\TPO_d$ of $\co{\bbp^n}(d)$
such that
$$
\dbar \TPO_d \ =\ \pi_*(w^d\bd)^{0,1}  
\tag{\TT.1}
$$ 
$$
\TPO_d\ \ \text{vanishes to order $d$ at } a,
\tag{\TT.2}
$$
and any solution of (\TT.1) differs from $\TPO_d$ by a global
holomorphic section of $\co{\bbp^n}(d)$, i.e., a homogeneous polynomial
of degree $d$ in $(z_0,..., z_n)$.

The main idea now is the following.  Suppose there exists a holomorphic
chain $\HC$ with $d\HC=\bd$ which misses $\bbp^{0}$. Then $\TP_d
\equiv \pi_*(w^d \HC)$ is a solution of (\TT.1) and therefore can be
written in the form $$
\TP_d\  \equiv\  \pi_*(w^d \HC)\ =\ \TPO_d +p_d
$$
for a unique homogenous polynomial 
$$
p_d \in H^0(\bbp^n, \co{}(d)) 
\cong \bbc[z_0,...,z_n]_d \equiv \ch(d).
$$
 Let us suppose further that
$\HC$ is {\bf positive over $a$}, i.e., positive in $\pi^{-1}(U)$ for
some neighborhood $U$ of $a$ in $\bbp^n-\pi(\bd)$. Then
$\TP_0\bigl|_U\equiv  \ell \in\bbz^+$ is the sheeting number of $\HC$ 
over $U$, and with respect to any non-vanishing section $\s$ of $\co {U}
(1)$, there are  functions $f_1,...,f_\ell$ (locally defined and
holomorphic outside a divisor, as in the Weierstrauss Preparation
Theorem) such that  $$ \frac{\TP_d}{\s^d} \ =\
f_1^d+f_2^d+\dots+f_\ell^d. \tag{\TT.3}$$ A sequence
$\{\TP_d\}_{d=0}^\infty$ satisying (\TT.3) for some, and therefore any, 
$\s$ is called  a {\bf Newton hierarchy  of level $\ell$}. Such
hierarchies are freely determined by the first $\ell+1$ functions
$\TP_0,...,\TP_\ell$. The remaining functions can be written
recursively by  explicit universal homogeneous polynomial expressions  
$$ 
\TP_d\  =\ Q_{d,\ell}(\TP_1,...,\TP_{\ell}) \qquad\text{ for }\
d>\ell.
\tag{\TT.4}
$$
(When $\ell=0$, then $Q_{d,\ell}=0$ for all $d$.) \ 
If these identities are satisfied then there exist holomorphic
sections $f_1\s,...,f_\ell\s$ such that (\TT.3) holds. 
This leads to the following.

\Theorem{\TT.1. (The Main Result for Hypersurfaces)} {\sl  There  exists
a holomorphic chain $\HC$   with $d\HC=\bd$ which is supported in 
$\bbp^{n+1}-\bbp^{0}$  and is positive over a point  $a\in
\bbp^p-\pi(\bd)$ if and only if there exist an integer $\ell\geq 0$ and
homogeneous polynomials $p_d \in  \ch(d)$ for $d=1,...,\ell$ such that
in some neighborhood of $a$
$$
\TPO_k\  \equiv\ Q_{k,\ell}(\TP_1,...,\TP_{\ell}) \mod \ch(k)
\qquad\text{ for all } k>\ell.
\tag{\TT.5}
$$
where $\TP_d=\TPO_d+p_d$ for $d=1,...,\ell$.}
 \medskip

Condition (\TT.5) constitutes a {\bf non-linear moment condition}.
The essential point is that the coefficients of degree $>k$ in
the power series expansion of $\TP_k$ at $a$ can be written, via the
Bochner-Martinelli kenel,  as certain explicit integrals over $\bd$.
Condition (\TT.5) says that these coefficients for $k>\ell$ must equal
certain polynomial expressions in the coefficients of 
$\TP_1,...,\TP_\ell$ and $p_1,..,p_\ell$. 

The polynomials $p_1,..,p_\ell$ are initially 
free.  However, in general the equations (\TT.5) will determine a
certain number them.  When we have an $\ell=\ell_0$ so that $\HC$ exists
and has no boundary-free components, then all of $p_1,..,p_{\ell_0}$
will be determined. However, if we then increase $\ell$, solutions will
continue to exist but there will remain some freedom in choosing
$p_1,..,p_\ell$.  This corresponds to the freedom in adding
positive boundary-free components to the original solution.

 All this is best illustrated by considering the special case of curves
in $\bbp^2$ (which is essentially as complicated as the general case).

\Ex{\TT.2. (Curves in $\bbp^2$)}  Let $\g\subset \bbp^2$ be a compact
oriented $C^1$-curve without boundary and not necessarily connected.
Fix a point $\bbp^0\in \bbp^2$ not on $\g$. Choose
homogeneous coordinates $[z_0\:z_1:z_2]$ so that $\bbp_0=[0:0:1]$ and
define $\pi:\bbp^2-\bbp^0\to\bbp^1$ by setting
$\pi([z_0:z_1:z_2])=[z_0:z_1]$ as above. Suppose $a\equiv[1:0]\notin
\pi(\g)$.  Consider the affine coordinate $z=z_1/z_0$ on the open subset 
$U=\{z_0\neq0\}\subset\bbp^1$. Then $w([z_0\:z_1:z_2]) =z_2$ is the
tautological section of $\pi^*\co{}(1)$ on $\bbp^2-\bbp^0$, and the
integral   
$$
\TP_d(z) \ =\ \frac{1}{2\pi i}\int_{\g}\frac{w^d}{\zeta-z}\,d\zeta
\tag{\TT.6}
$$ 
solves the equation     \footnote{
Note that this can be thought of as a scalar equation
by writing
 $$
\TP_d(z) \ =\ \left\{{{\tau_d(z)}\over{z_0^d}}\right\}  z_0^d \ =\ 
\left\{{{1}\over{2\pi i}}\int_{\g} {{(w/z_0)^d}\over{\zeta-x}}\,d\zeta
\right\} z_0^d 
$$
with respect to the cross-section $\s=z_0$ of $\co{\bbp^1}(1)$ 
and then setting $ z_0=1$.} $\dbar \TP_d = \pi_*(w^d\g)^{0,1}$.

Expanding (\TT.6) in a neighborhood of 0 we find that
$$
\TP_d(z) \ =\  \sum_{k=0}^\infty\left\{
\frac{1}{2\pi i}\int_{\g}\frac{w^d}{\zeta^{k+1}}\,d\zeta\right\}z^k
\ \equiv\ \sum_{k=0}^\infty A_k(d)z^k
$$
and so
$$
\TPO_d(z)=\sum_{k=d+1}^\infty A_k(d)z^k
$$
where 
$$
A_k(d) \ =\ 
\frac{1}{2\pi i}\int_{\g}\frac{w^d}{\zeta^{k+1}}\,d\zeta
\tag3.7
$$
are classical moments of the curve $\g$.

We now have free polynomials $p_d(z)  = c_0(d)+c_1(d)z+\dots+c_d(d)z^d$
for $d=1,...,\ell$. Plugging into the equations (\TT.5) gives a
sequence of polynomial relations in the finite set $\{c_k(d)\}$ and
the moments $\{A_k(d)\}_{k>d>\ell}$. 

Consider for example the special case  $\ell=1$. Then there is only one
free polynomial
$$
p_1(z)\ =\ c_0+c_1z,
$$
which means two free constants to determine. In this case the equations
(\TT.5) have the particularly simple form
$$
\TP_d(z)\ \equiv \ \left( \TPO_1(z)\right)^d\ \equiv \ 
\left( c_0+c_1z+ \TP_1(z)\right)^d \mod \ch(d)
$$
for $d\geq 2$. The first two non-zero coefficients of $\TP_2$ give
$$
A_3(2)=2A_3(1)c_0 + 2A_2(1)c_1\and
A_4(2)=2A_4(1)c_0 + 2A_3(1)c_1+A_2(1)^2
$$
Assume $A_3(1)^2-A_2(1)A_4(1)\neq0$ and for notational convenience
let $c_k=A_k(1)$  for $k=0,1$ denote the solutions of the equations
above. Then we have the following.

\Cor{\TT.3} {\sl Necessary and sufficient
conditions that $\g$ bound a holomorphic 1-chain which is positive and
1-sheeted over 0 is that the non-linear moment conditions
$$\boxed{
A_k(d)\ =\ \sum_{j_1+\dots+ j_r=k}c_{d,J}\, A_{j_1}(1)\cdots A_{j_r}(1)
}
$$
hold for all $k>d>2$ where  $c_{d,J}$ are the obvious combinatorial
constants.  Alternatively, in recurssive form this can be written as
$$
\boxed{
A_k(d)\ =\ \sum_{i+j=k}A_i(d-1)A_j(1)}
$$
for all $k>d\geq 2$.}

\Note{\TT.4}  The case $\ell=0$ corresponds to the existence of a
solution with compact support in
$\bbc^2=\bbp^2-\overline{\pi^{-1}(a)}$.  Here our moment conditions
reduce to the linear moment conditions $A_k(d)=0$ for all $k>d>0$. 
Making the substitution $z=1/\zeta$ in (\TT.7) we recover exactly the
moment condition of [HL$_1$].

  \vskip .3in

%
%

\centerline{\bf  \S \AA. Newton Hierarchies}
\vskip .2in
\def\c{c}

The following material is classical but central to our results.
For the convenience of the reader we review this  material
in a form that is particularly adapted to our needs.

\Def{\AA.1} A sequence of complex numbers $\{ \c_d\}_{d=0}^{\infty}$
is called a {\bf Newton hierarchy} of level $\ell$ if there exist
complex numbers $b_1,\dots,b_{\ell}$ such that
$$
\c_d \ =\ b_1^d+\dots+b_{\ell}^d  \qquad\qquad\text{for all } \ d\geq 0.
$$
Note that $\c_0=\ell$. 
\medskip

The main fact is that in a Newton hierarchy of level $\ell$ the 
numbers $\c_1,...,\c_{\ell}$ can be freely prescribed and the remaining
terms in the sequence are given by explicit polynomial expressions in 
terms of these.

\Ex{\AA.2. ($\ell=1$)} \ \ $\c_d = \c_1^d$ for all $d$.  

\Ex{\AA.3. ($\ell=2$)} \ \ 
$\c_d=b_1^d+b_2^d$ for $d=0,1,2$. Then $b_1$ and $b_2$ are the roots of
the polynomial $x^2-(b_1+b_2)x+b_1b_2 = x^2-\c_1x+\frac12(\c_1^2-\c_2)$.
Thus, setting $S_1=\c_1$ and $S_2 =\frac12(\c_1^2-\c_2)$ we see that
$b_j^{d} - S_1 b_j^{d-1} + S_2b_j^{d-2} =0$ for all $d \geq 2$.
Consequently,
$$
\c_d -S_1\c_{d-1}+S_2 \c_{d-2} \ =\ 0
 \qquad\qquad\text{for all } \ d > 2.
$$

\Prop{\AA.4} {\sl Let $\{S_k(\c_1,...,\c_k)\}_{k=1}^{\infty}$
be the sequence of polynomials defined recursively by the equations
$$
\c_{\ell}-S_1 \c_{\ell-1}+S_2 \c_{\ell-2} - \dots +(-1)^{\ell}S_{\ell} \ell
\ =\ 0.
$$
Then setting
$$
\sigma_k(x_1,...,x_N)\  =\ \sum_{i_1<\dots<i_k} x_{i_1}\dots x_{i_k}
\qquad\ \ \text{and}\ \ \qquad
\c_k(x_1,...,x_N)\ =\ \sum_{i=1}^N x_{i}^k
$$
we have
$$
\sigma_k(x)\  =\ S_k(\c_1(x),...,\c_k(x))
$$
in $\bbc[x_1,...,x_N]$ for all $N\geq k$.
}

\pf We first prove the assertion for $N=k$.  This  is easily checked for
$k=1,2$. Assume inductively that it holds for all $k'<k$. Note that
$$
(t-x_1)\dots(t-x_k)\ =\
t^k-\sigma_1(x_1,...,x_k)t^{k-1}+\dots+(-1)^k\sigma_k(x_1,...,x_k). 
$$
Substituting $t=x_j$ and summing over $j$ gives
$$
\c_k(x)-\sigma_1(x)\c_{k-1}(x)+\dots+(-1)^k\sigma_k(x)k\ =\ 0.
$$
for $x=(x_1,...,x_k)$. Therefore, by induction
$$\aligned
\s_k(x)\ &=\
\frac{(-1)^{k+1}}{k}\bigg\{\c_k(x)-\sigma_1(x) \c_{k-1}(x)+
\sigma_2(x) \c_{k-2}(x)- 
\sigma_3(x) \c_{k-3}(x)+\dots\bigg\} \\ &=\
\frac{(-1)^{k+1}}{k}\bigg\{\c_k(x)-S_1(\c_1(x))\c_{k-1}(x)+
S_2(\c_1(x),\c_2(x))\c_{k-2}(x)\\
&\qquad\qquad\qquad\qquad\qquad\qquad\qquad\ \ -
S_3(\c_1(x),\c_2(x),\c_3(x))\c_{k-3}(x)+\dots\bigg\}\\
&=\ 
S_k(\c_1(x),...,\c_k(x))
 \endaligned
$$
for $x=(x_1,...,x_k)$. It remains only to establish  this identity
for $x=(x_1,...,x_N)$ where $N>k$. For this we note that 
$
P(x_1,...,x_N)\equiv \s_k(x_1,...,x_N)
-S_k(\c_1(x_1,...,x_N),...,\c_k(x_1,...,x_N))
$
is a symmetric homogeneous polynomial of degree $k$ such that 
$
P(x_1,...,x_k,0,...,0) =0
$
by the identity proved above.  By symmetry it follows that 
$
P(0,...,0,x_{i_1},0,...,0,x_{i_k},0,...,0) =0
$
for all $1\leq i_1<i_2<\dots<i_k\leq N$. It follows that all the
coefficients of $P$ are zero. \qed

\Cor{\AA.5} {\sl The sequence $\{\c_d\}_{d=0}^{\infty}$ is a Newton
hierarchy of level $\ell$ if and only if $\c_0=\ell$ and
$$
\c_k -S_1(\c_1)\c_{k-1}+S_2(\c_1,\c_2)\c_{k-2}- \dots +(-1)^{\ell}
S_{\ell}(\c_1,...,\c_{\ell})\c_{k-\ell}\ =\ 0
\tag{\AA.1}
$$
for all $k>\ell$.

In other words, fix $\ell\geq1$ and let $Q_{k,\ell}(\c_1,...,\c_\ell)$ for
$k>\ell$ be the sequence of polynomials generated recursively by setting
$\c_j=Q_{j,\ell}$ in
equation  (\AA.1). Then the sequence $\{\c_d\}_{d=0}^{\infty}$ is a Newton
hierarchy of level $\ell$ if and only if
$$
\c_k\ =\ Q_{k,\ell}(\c_1,...,\c_\ell)\qquad\text{ for all } k>\ell.
\tag{\AA.2}
$$
}

\pf Let $b_1,...,b_{\ell}$ be the roots of the polynomial 
$
p(t)=t^{\ell}-S_1(\c_1)t^{\ell-1}+S_2(\c_1,\c_2)t^{\ell-2}-\dots+
(-1)^\ell S_\ell(\c_1,...,\c_\ell).
$
Then $\c_d=b_1^d+\dots+b_{\ell}^d$ for $d=1,...,\ell$ by Proposition \AA.4.
Substituting $t=b_j$ into $t^{k-\ell} p(t)$ and summing over $j$ completes
the proof. \qed
\medskip
This gives the main result.

\Theorem{\AA.6} {\sl Each $\ell$-tuple $(\c_1,...,\c_\ell)\in \bbc^\ell$
extends to one and only one Newton hierarchy 
$\{\c_d\}_{d=0}^{\infty}$ of level $\ell$. The numbers $\c_d$ for $d>\ell$
are given by the universal polynomials (\AA.2).
}

\Note{\AA.7} Let $\cn_{\ell}$ denote the
space of Newton Hierarchies of level $\ell$.  Then
$\cn\equiv \oplus_{d\geq0} \cn_d$ has a  {\sl  vector monoid structure}.
Given  $\bc =\{c_d\}\in \cn_{\ell}$, 
$\wt{\bc} =\{\wt{c}_d\}\in \cn_{\wt \ell}$,  and $t\in \bbc$, we have  
$$
t\cdot \bc \ =\ \{t^d c_d\}_{d=0}^\infty \in \cn_{\ell}
\and
\bc +\wt{\bc} \ =\ \{c_d+\wt c_d\}_{d=0}^\infty \in \cn_{\ell+\wt\ell}.
$$

  \vskip .3in

%
%


%
       
\def\E{E}
\def\ta{c}

\centerline{\bf  \S \FF.  Newton Hierarchies with Coefficients in a
Vector Bundle }\smallskip

\centerline{\bf  and their Associated Multisections.} \vskip .2in

The notion of a Newton hierarchy extends immediately to
powers of a line bundle. Let $\pi:\E \to X$ be a holomorphic line
bundle over a complex manifold $X$ of dimension n, 
and suppose we are given holomorphic sections 
$$
\TP_d \in H^0(X;\,\co{}({\E}^d)) 
\qquad\text{ for each integer } d\geq0
\tag{\FF.1}$$
where ${\E}^d=\E\otimes\dots\otimes\E$ ($d$-times).  With respect to
a local trivialization $\s\in H^0(U;\,\co{}(\E))$ of $\E$ over an open
subset $U\subset X$ we can write   
$$
\TP_d \ =\ \ta_d \s^d \foralll d\geq 0.
$$
where $\ta_d\in \co{}(U)$. If $\wt \s\in H^0(U;\,\co{}(\E))$ is another
trivialization with $\TP_d = \wt \ta_d\wt\s^d$, then 
$$
\wt \ta_d\ =\ \left(\frac{\s}{\wt\s}\right)^d \ta_d \qquad\text{ in } U.
\tag{\FF.2}$$

\Def{\FF.1}  A sequence of sections $\{\TP_d\}_{d=0}^\infty$ as in (\FF.1)
with $\TP_0\equiv \ell \in \bbn$ is called a {\bf Newton Hierarchy of level 
$\ell$ with coefficients in } $\E$ if each point $x\in X$ has a
neighborhood $U$ and a local trivialization $\s$ of $\E\bigl|_U$
such that the sequence of functions $\{\TP_d/\s^d\}_{d=0}^\infty$ is a
Newton hierarchy of level $\ell$ at each point of $U$.  By (\FF.2) if this
condition holds for one trivialization over $U$, it holds for every
trivialization.
\medskip

\Def{\FF.2} By a {\bf multi-section of degree $\ell$ } of a line bundle
$\pi:\E\to X$ we mean an effective divisor  $D$ on the total space of $\E$
such that $\pi\bigl|_{D}$ is proper and of degree $\ell$.
\medskip

Recall that by an {\bf effective divisor} we mean a positive holomorphic
n-chain on $\E$ where $n=\dim X$, and by {\bf degree} $\ell$ we mean that
$\pi_*D=\ell[X]$ (or equivalently that
the intersection 0-cycle $D\bullet\pi^{-1}(x)$ (cf. [Fu])  has degree
$\ell$ for all $x\in X$).
The first main result of this section is the following.

\Theorem{\FF.3} {\sl Let $\pi:\E\to X$ be a holomorphic line bundle over
a complex manifold $X$.  Then there is a one-to-one correspondence between
Newton Hierarchies of level $\ell$ with coefficients in $\E$ and
multi-sections of $\E$ of degree $\ell$.

This correspondence is given as follows.  Let $w\in H^0(\E,
\co{}(\pi^*\E))$  denote the tautological cross-section $w(e)=e$.
Then for a multi-section $D$ the associated Newton hierarchy is given by
$$
\TP_d \ =\ \pi_*(w^d D) \foralll d\geq 0.
\tag{\FF.3}
$$
}\medskip
 
\pf
Since all notions are invariantly defined, it
suffices to prove the result locally on $X$.
Let $U\subset X$ be an open subset of $X$ and $\s$ a trivialization of 
$\E$ over $U$. This gives an isomorphism 
$$
U\times \bbc \ @>{\cong}>>\ \E\bigl|_U \qquad\text{ defined by }\ \ 
(x,t)\ \mapsto\  t\s(x).
\tag{\FF.4}
$$
In $U\times \bbc$ the tautological section $w$ has the form $w(t\s(x))=
t\s(x)$ and so 
$$
w/\s =t.
$$ 
Under (\FF.4) the restriction of $D$ to $\E\bigl|_U$ becomes an effective
divisor in $U\times \bbc$ which is proper and of degree $\ell$ over
$U$, which shall be denoted by $D_U$. Now it is classical that
$$
\pr_*(t^dD_U)\ =\ f_1^d +\dots + f_{\ell}^d
$$
where $\pr:U\times \bbc\to U$ is the projection, and
$\{f_1(x),...,f_\ell(x)\}\equiv \pi^{-1}(x)\bullet D_U\subset \bbc$ are the
points of $D_U$ above $x\in U$ listed to multiplicity (See [H], [HL$_1$]
for example). It follows immediately that $\ta_d \equiv \pr_*(t^dD_U)$ is a
Newton hierarchy of level $\ell$ and so $\TP_d \equiv \ta_d \,\s^d$
is a level-$\ell$ Newton hierarchy with coefficients in $\E$.

For the converse suppose $\{\TP_d\}_{d=0}^{\infty}$ is a Newton hierarchy
with coefficients in $\E$, and write $\TP_d = \ta_d \s^d$ over $U$. Define
holomorphic functions $s_k \in \co{}(U)$ by 
$$
s_k(x) \ =\ S_k(\ta_1(x),...,\ta_\ell(x)) \qquad\text{ for }\
k=1,...,\ell 
$$
where the $S_k$ are the polynomials defined in Proposition \AA.4.
Define $P\in \co{}(U\times \bbc)$ by 
$$
P(x,t)\ =\
t^\ell-s_1(x)t^{\ell-1}+s_1(x)t^{\ell-1}-\dots+(-1)^\ell s_{\ell}(x) 
\tag{\FF.5}$$
and 
$$
D_U\ \equiv \ \Div(P) \ =\ \frac 1 {2\pi}d d^C \log|P|.
$$
 This is a divisor in $U\times \bbc$ which is proper over $U$ and has the
property that $\pr_*(t^dD_U)=\ta_d$ for all $d$. \qed

\Note{\FF.5} One can check directly that the divisor of $P$ is
invariantly defined. If $\wt\s$ is another trivialization of $\E$
over $U$, and if we write $\wt t = (\s/\wt\s) t$, $\wt\ta_d = (\s/\wt\s)^d
\ta_d$ as above, then setting $\wt s_k = S_k(\wt\ta_1,...,\wt\ta_\ell)$ we
have $\wt P(x,\wt t) = (\s/\wt\s)^\ell P(x,t)$, which shows that
$\Div(\wt P)$ is the image of $\Div(\wt P)$ under the change of
trivialization.
\medskip

We now pass to an important generalization of the discussion above.
Consider the $q$-fold direct sum 
$$
\E^{\oplus q}\  =\  \E\oplus\dots\oplus \E  \ @>{\Pi}>>\ X.
$$
For each $\la=(\la_1,...,\la_q) \in \bbc^q$ consider the projection 
$$
\pi_\la:\E^{\oplus q} \ \arr\ \E\qquad \text{ given by }\ \ \ \ 
\pi_\la(e_1,\dots,e_q) \ =\ \sum_{j=1}^q\la_je_j.
$$
Note the commutative diagram
$$\CD
\E\oplus\dots\oplus \E @>{\pi_\la}>>\ \ \ \ \  \E\\
\ \ \ \ \Pi\searrow @.\!\!\!\!\swarrow \pi\\
\ \ \ \ \ \ \ \ \ \qquad \qquad \qquad \ \ \ \ \ X
\endCD
$$\smallskip

Let $w$ denote the tautological cross-section of $\pi^*\E $ over $\E$ as
above and note that
$$
\pi_{\la}^*w\ =\ \sum_{j=1}^q\la_jw_j\ \equiv\ \la \cdot {\bold w}
\tag{\FF.6}$$
where $w_j$ denotes the $j\th$ tautological cross-section of
$\Pi^*\E\to \E^{\oplus q}$, i.e.,  $w_j(e_1,...,e_q) = e_j$.

Suppose now that $D$ is a { \bf multi-section of $\E^{\oplus q}$
of degree $\ell$}, that is, 
 a positive holomorphic $n$-chain on    
$\E^{\oplus q}$ which is proper over $X$ and of degree $\ell$
(i.e., $\Pi_*D=\ell[X]$). Then for each $\la$,   $\pi_\la$ is  proper on
the support of $D$ and $D_\la\equiv (\pi_\la)_*D$ is a mult-section of $\E$
of degree $\ell$. Hence there is an associated Newton hierarchy 
$\{\TPF d {\la}\}_{d=0}^\infty$ of level   
$\ell$ on $X$ (with coefficients in $\E$) defined by 
$$
\TPF d {\la}\ =\ \pi_*(w^d D_\la)\ =\ \Pi_*\{(\la\cdot{\bold w})^d D\}.
\tag{\FF.7}
$$
Note that $\TPF d {\la}$ is a homogeneous  polynomial of degree $d$ in
$\la$ with values in $H^0(X, \co{}(E^d))$, i.e.,  
$
{\bold \TP}_d\ \in\ \text{Sym}^d(\bbc^q)^{\vee}\otimes H^0(X, \co{}(E^d)).
$
This polynomial can be re-expressed as
$$ 
\TPF d {\la}\ =\ \sum_{|\a|=d}
{\tsize\binom{d}{\a}} \la^\a\Pi_*\{w^\a D\}\ =\ \sum_{|\a|=d}
{\tsize\binom{d}{\a}}\TP_{\a}  \la^\a
\tag{\FF.8}
$$
where $w^\a=w_1^{\a_1}\dots w_1^{\a_q}$ and $|\a|=\a_1+\dots\a_q$ as
usual, and
$$
\TP_{\a}\ \equiv\Pi_*\{w^{\a} D\}
\tag{\FF.9}
$$
is the {\bf $\a\th$ moment} of $D$.

Conversely, suppose we are given elements ${\bold \TP}_d \in
\text{Sym}^d(\bbc^q)^{\vee}\otimes H^0(X, \co{}(E^d))$ for $d\geq 0$, or
equivalently, moments $\TP_\a\in H^0(X,\co{}(E^{|\a|})$ for $\a\geq
(0,...,0)$.

\Def{\FF.6} The family $\{\TP_{\a}\}_{\a\geq 0}$ is called a 
{\bf Newton family of level $\ell$ with coefficients in $E^{\oplus q}$}
if the polynomials $\TPF d {\la}= \sum_{|\a|=d}
{\tsize\binom{d}{\a}}\TP_{\a}  \la^\a$ satisfy the condition that
$$
\{\TPF d {\la}\}_{d=0}^\infty \ \text{\ is a Newton hierarchy of level
$\ell$ for all } \la\in \bbc^q.
\tag{\FF.10}$$

It is straightforward to check from
the result above and arguments in [HL$_1$, \S 7] that for any such
family there exists a multi-section $D$ of degree $\ell$  of $\E^{\oplus
q}$ such that $\TP_\a=\Pi_*(w^\a D)$ for all $\a$. Hence we obtain the
following.

\Theorem{\FF.7} {\sl Let $\pi:\E\to X$ be a holomorphic line bundle over
a complex manifold $X$.  Then there is a one-to-one correspondence between
Newton families of level $\ell$ with coefficients in $\E^{\oplus q}$ and
multi-sections of $\E^{\oplus q}$ of degree $\ell$. 

This correspondence associates to a multi-section $D$ the family of
moments $\TP_\a=\Pi_*(w^\a D)$ where $w^\a=w_1^{\a_1}\dots w_q^{\a_q}$
and $w_j$ is the $j\th$ tautological cross-section of $\Pi^*E^{\oplus q}$ 
defined above.}
 
 \vskip.3in

We now observe that {\bf  Newton families $\{\TP_{\a}\}_{\a\geq 0}$ of
level $\ell$ with coefficients in $\E^{\oplus q}$  are universally
determined by their terms of degree $\leq \ell$ } as in (\AA.2).  Indeed,
by Corollary \AA.5, our condition (\FF.10) above is equivalent to the
condition that 
$$
\TPF d {\la}\ =\ Q_{d,\ell}\left(\TPF 1 {\la},...,\TPF {\ell} {\la}\right)
\qquad\text{for all } d\geq \ell\text{ and } \la\in \bbc^q.
\tag{\FF.11}$$
These universal polynomials $Q_{d,\ell}(\x_1,...,\x_{\ell})$ introduced
in 4.5 are {\sl weighted homogeneous}, that is, 
$$
Q_{d,\ell}(t\x_1,t^2\x_2,...,t^\ell\x_{\ell})\ =\ t^d
Q_{d,\ell}(\x_1,...,\x_{\ell})\qquad\text{for all } t\in\bbc
$$
or equivalently
$$
Q_{d,\ell}(\x_1,...,\x_{\ell})\ =\ \sum_{\b_1+2\b_2+\dots+\ell\b_{\ell}=d}
c_{\b}\, \x^{\b}
$$
Let $e_1,...,e_q$ denote the standard basis of $\bbc^q$ and write
$\TP_\b\equiv \TP_{\b\cdot e}\equiv \TP_{\Sigma_i\b_i e_i}$. One sees
directly that 
$$
Q_{d,\ell}\left(\sum_i  \TP_{e_i} \la_i,\sum_{ij}
\TP_{e_i+e_j} \la_i\la_j,...,
\sum_{|\b|=\ell} {\ssize \binom d {\b}}\TP_{\b} \la^\b\right)\ =\ 
\sum_{|\a|=d} Q_{\a,\ell}(\TP_{})\la^\a
$$
is a homogeneous polynomial of degree $d$ in $\la$ whose coefficients
$Q_{\a,\ell}(\TP_{})$ are universal polynomials in the
indeterminates $\TP_{}=\{\TP_\b\}_{|\b|\leq \ell}$. This gives the
following generalization of Corollary  \AA.5 to this case.

\Prop{\FF.8} {\sl Let $\{\TP_\a\}_{\a\geq0}$ be a family of continuous
sections of a complex line bundle $E$ on a space $X$.  Then
$\{\TP_\a\}_{\a\geq0}$  is a Newton family of level $\ell$ with
coefficients in $E^{\oplus q}$  iff 
$$
\TP_\a\ =\  Q_{\a,\ell}(\TP_{}) \qquad \text{ for all } |\a| >\ell
\tag{\FF.12}
$$
where $\TP=\{\TP_\b\}_{|\b|\leq\ell}$.
}

  \vskip .3in

%
%


%

\noindent
{\bf \S \CC. Vector Bundle-Valued Currents.}  In this section we
prepare the general analytical tools for solving the main problem. 
In the next section we shall examine the concrete cases
of interest.

For a smooth vector bundle $E$ over a manifold $X$ we denote by
$\ce_p'(X,E)$ the space of $p$-dimensional currents of compact
support on $X$ with values in $E$.  By definition this is the
topological dual of the  space $\ce^p(X,E^*)$ of smooth $p$-forms
with values in $E^*$.  To any smooth map $f:Y\to X$ between
manifolds  one associates a continuous map  
$$
f_*:\ce_p'(Y, f^*E) \ \arr\ \ce_p'(X,E)
$$
given by $f_*(T)(\phi)=T(f^*\phi)$ for smooth forms $\phi\in
\ce^p(X,E^*)$

When $X$ and $Y$ are complex manifolds and $E$ and $f$ are
holomorphic, there is a decomposition of currents
$
\ce_p'(X,E)=\oplus_{r+s=p}\ce_{r,s}'(X,E)
$
preserved by $f_*$, and an operator 
$$
\dbar:\ce_{r,s}'(X,E)\ \arr\ \ce_{r,s-1}'(X,E)
$$
which commutes with $f_*$.  The space $\ce_{r,s}'(X,E)$ can be
canonically identified with the space of compactly supported
$(n-r,n-s)$-forms with distribution coefficients.
Hence we shall write
$$
{\ce}_{r,s}'(X,E)\ \equdef\ {\ce'}^{n-r,n-s}(X,E)
$$ 
Under this identification $\dbar$ is identified with the usual
operator and we have the following.
\Theorem{\CC.1 (Dolbeault)}{\sl There is a canonical isomorphism
$$
H^{p,q}({\ce'}^{*,*}(X,E))\ \cong\
H^q_{\text{cpt}}(X;\,\Omega^p\otimes E). 
\tag{\CC.1}
$$
where the right hand side denotes sheaf cohomology with
compact supports and coefficients  in the sheaf of holomorphic 
$p$-forms with values in $E$.}

  \vskip .3in

%
%


%

\def\G{T}

\noindent
{\bf \S \DD. The Projective $\dbar$-Problem.}  \ \ Let $\bbp^n$ 
denote complex projective $n$-space with homogeneous coordinates
$[z_0:\dots:z_n]$, and let  $\bbp^{n-p-1}\subset\bbp^n$ denote the
linear subspace defined by $z_0=z_1=\dots=z_p=0$. Consider the
linear projection
$$
\Y\ \equdef\ \bbp^n-\bbp^{n-p-1}\ @>{\pi}>>\ \bbp^p
\tag{\DD.1}$$
defined by $\pi([z_0:\dots:z_n]) = [z_0:\dots:z_p]$.
Note that $\pi:\Y \to \bbp^p$ is isomorphic to the vector bundle
$\co {\bbp^p}(1)^{\oplus (n-p)}$.

\Lemma{\DD.1} {\sl There is an isomorphism}
$
\pi^*\co {\bbp^p}(1)\ \cong\ \co {\bbp^n}(1)\bigl|_{\Y}.
$

\pf  This can be seen explicitly. Note that 
$$
\co {\bbp^n}(-1)
=\{([z],\zeta) \in \bbp^n\times \bbc^{n+1} : \zeta=\lambda z \
\text{  for some  } \la\in\bbc\}, 
$$
$$
\pi^*\co {\bbp^p}(-1)  =\{([z],\zeta) \in \Y\times \bbc^{p+1} :
\zeta=\lambda (z_0,\dots,z_p) \ \text{  for some  } \la\in\bbc\}.
$$
We define $\Psi:\pi^*\co {\bbp^p}(-1)\to \co
{\bbp^n}(-1)\bigl|_{\Y}$ by 
$$
\Psi([z_0:\dots:z_n], (\zeta_0,\dots,\zeta_p))
=  ([z_0:\dots:z_n], (\zeta_0,\dots,\zeta_n))
$$
 where for given
$(z_0,...,z_n)$ we have a unique $\la\in\bbc$ with  $\zeta_k=\la z_k$
for $k=0,...,p$ and we define $\zeta_k=\la z_k$ for  $k>p$.
The adjoint $\Psi^*$ gives the desired isomorphism. \qed

\Def{\DD.2}  The homogeneous coordinate functions $z_0,...,z_n$ are
holomorphic sections of $\co{\bbp^n}(1)$. By the Lemma above
each
$$
w_k\ =\ z_{p+k}\bigl|_Y
$$
is a holomorphic section of 
$
\co {\bbp^n}(1)\bigl|_{\Y}\ \cong\ \pi^*\co {\bbp^p}(1).
$
The section $w_k$ will be called the {\bf $k\th$  tautological
cross-section} of $\pi^*\co {\bbp^p}(1)$.

\medskip
It follows from Lemma \DD.1 that there are isomorphisms 
$
\pi^*\co {\bbp^p}(d)\ \cong\ \co {\bbp^n}(d)\bigl|_{\Y}
$
for all integers $d$.  Hence, by \S  \CC \ the projection $\pi$
induces continuous homomorphisms on currents:
$$
\pi_*:\ce_{*,*}'(\Y,\co{\Y}(d))\ \arr\ 
\ce_{*,*}'(\bbp^p,\co{\bbp^p}(d)).
\tag{\DD.2}
$$

\Prop{\DD.3} {\sl  Let $\G \in \ce_{p,p-1}'(\Y)$ be a $\dbar$-closed
current on $\Y$ and $\varphi\in H^0(\bbp^n, \co{\bbp^n} (d))$
a holomorphic section of $\co{\bbp^n} (d)$ for $d\geq 0$. Then
$
\dbar \pi_*(\varphi \G) \ =\ 0
$
and the associated $\dbar$-problem has solutions}
$$
\{S\in {\ce'}^0(\bbp^p,\co{\bbp^p}(d)) : 
\dbar S=\pi_*(\varphi \G)  \} 
\ \cong\ H^0(\bbp^p,\co{\bbp^p}(d)). 
$$

\pf
We have 
$
\pi_*(\varphi \G) \in \ce_{p,p-1}'(\bbp^p,\co{\bbp^p} (d))
\ \cong\ {\ce'}^{0,1}(\bbp^p,\co{\bbp^p} (d))
$
and since $H^1(\bbp^p,\co{\bbp^p}(d))=0$ there exists at least one
current $S$ with $\dbar S=\G$  by Theorem \CC.1. For any other
solution $S'$ we have $\dbar (S-S')=0$ and the result follows from
the regularity of distributions in ker$(\dbar)$. \qed

  \vskip .3in

%
%

\noindent
{\bf \S \EE. Curves in $\bbp^2$.}  \ \ 
Let $\gamma\subset \bbp^2$ be a finite disjoint union of oriented 
closed curves of class $C^1$, and suppose $\gamma$ is the (current)
boundary of a holomorphic 1-chain $\HC$. Choose a point $\bbp^0\notin
\supp(\HC)$ and consider the projection
$$
\Y\equiv \bbp^2-\bbp^0 \ @>{\pi}>>\ \bbp^1
\tag{\EE.1}
$$
defined in \S \DD. Recall that $\pi$ is isomorphic to the line bundle
$\pi:\co{}(1)\to \bbp^1$, and so $\pi^*\co{}(1)$ has a tautological
cross-section which we denote by $w$. Thus, for each $d\geq 0$ the
current $\pi_*(w^d \HC)$ is a generalized section of $\co{\bbp^1}(d)$
which satisfies the equation 
$$
\dbar \pi_*(w^d \HC)\ =\ \pi_*(w^d \gamma)^{0,1}.
\tag{\EE.2}
$$
In particular, $\pi_*(w^d \HC)$ is a holomorphic section of
$\co{}(d)$ over  $\bbp^1-\gamma$. 
 
Recall that $\HC$ is said to be {\bf positive}
over a point $a\in \bbp^1-\pi(\gamma)$ if its
restriction  to $\pi^{-1}(U)$ is  positive 
for some neighborhood $U$ of $a$.  From Theorem \FF.3 
we know that:

 \Lemma{\EE.1}{\sl If $\HC$ is positive over a point $a\in
\bbp^1-\pi\gamma$, then the holomorphic sections $\{\pi_*(w^d
\HC)\}_{d=0}^{\infty}$ form a Newton hierarchy of level $\ell$ in a
neighborhood of $a$ where
$\ell= \pi_*( \HC)$ =
the sheeting number of $\HC$ at $a$}. 
\medskip

We now look in the converse direction. Let  
$\gamma\subset \bbp^2$ be as above and ask whether $\gamma$
bounds a holomorphic 1-chain in $\bbp^2$. If it does,
then it bounds infinitely many such chains since we are free to add
any chain $S$ with $dS=0$. To eliminate this
ambiguity we first fix the homology class $u\in
H_2(\bbp^2,\gamma;\,\bbz)$ of the   solution
$\HC$.  We assume  $\partial u=[\gamma]$, and so by the short
exact sequence:
$
0\ \arr\ H_2(\bbp^2;\,\bbz)\ \arr\ H_2(\bbp^2, \gamma;\,\bbz)\ 
@>{\partial}>>\ H_1(\gamma;\,\bbz) \ \arr\ 0,
$
the class  $u$ is determined by its intersection with any projective
line $L$ which does not meet $\g$.  In particular, for  $a\in
\bbp^1-\pi\gamma$ we set $L_a=\pi^{-1}(a)$ and note that 
$$
L_a \cdot u \ =\ \pi_* u(a)\ =\ \text{the local net sheeting number 
of $\HC$ over } a.
$$
If $d\HC=\g$, then $\pi_*\HC$ is an integer-valued
function whose value at $a\in \bbp^1-\pi\gamma$ is precisely the net
sheeting number of $\HC$ over  $a$.  It satisfies the equation
$d(\pi_*\HC)=\pi_*\g$, or equivalently
$\dbar(\pi_*\HC)=(\pi_*\g)^{0,1}$ and is uniquely determined by this
equation up to an integer constant. Choosing this constant is
equivalent to choosing the homology class $u$.

Unfortunately, if there exists one solution in a given homology
class, there exist infinitely many in that class, since one can add 
algebraic cycles homologous to zero.  This
indeterminacy is greatly reduced if we require $T$   
to be positive over a chosen point  $a\in \bbp^1-\pi\gamma$. Then
there can exist at most a finite-dimensional family of solutions.
Moreover, if such solutions exist for some $\ell>0$, then there is 
a minimal $\ell_0\geq0$ where the problem is solvable, and for
$\ell_0$ the solution $\HC$ is {\bf unique}.  For each
$\ell>\ell_0$ all solutions are  of
the form $\HC+S$ where $S$ is in the finite-dimensional family of
positive algebraic 1-cycles  homologous to
 $(\ell-\ell_0)\,\bbp^1$.

With this motivation we fix a point $a\in \bbp^1-\pi\gamma$
and an integer $\ell\geq0$,
and we look for holomorphic 1-chains $\HC$ with $d\HC=\g$ which are
positive and $\ell$-sheeted over $a$.

To begin, observe that since $d\g=0$, the currents
$\pi_*(w^d\gamma)^{0,1}$ for $d\geq0$ satisfy
$$
\dbar \pi_*(w^d\gamma)^{0,1}\ =\ 0  \qquad\text{ on }\ \ \bbp^1.
$$
Therefore, by Proposition \DD.2 the space 
$$ 
\ct_{d,\gamma}\ \equdef\ \{\TP_d\in {\ce'}^0(\bbp^1, \co{}(d))\  :\
\dbar \TP_d=\pi_*(w^d \gamma)^{0,1}\}
\qquad\text{has dimension}\ \ d+1.
\tag{\EE.3}
$$
Any two elements of $\ct_{d,\gamma}$ differ by a unique element 
of $H^0(\bbp^1,\co{}(d))$, i.e., a homogeneous polynomial of degree
d in homogeneous coordinates $[z_0\:\!z_1]$ on $\bbp^1$.
This means the following.
Let $z$ be an  affine coordinate  on $\bbp^1$ 
with $z(a)=0$, choose a trivialization of
$\co{}(1)$ near 0 and for $\TP_d\in \ct_{d,\gamma}$ consider the
power series expansion 
$$
\TP_d(z)\cong\sum_{k=0}^\infty \a_k(d)z^k 
$$

\vskip.3in
\noindent
{\bf Observation \EE.3}
 {\sl 

(1)\ \ The coefficients
$\{\a_k(d)\}_{k=d+1}^\infty$ are independent of the choice of 
$\TP_d \in \ct_{d,\gamma}$ since any two elements in
$\ct_{d,\gamma}$ differ  by a polynomial 
of degree $\leq d$ in $z$.

(2)\ \ There is a unique element $\TPO_d \in \ct_{d,\gamma}$
which vanishes to order $d$ at $a$.
 }
\medskip

To solve our problem we are seeking currents $\TP_d\in
 \ct_{d,\gamma}$ which form a Newton heirarchy of level $\ell$ in a
neighborhood of $0$ (and are therefore of the form 
 $\pi_*(w^d\HC )$ for some positive divisor of degree $\ell$ above
that neighborhood). By Observation 8.3 there is some freedom in
constructing these currents which comes 
from the polynomial ambiguity in $\TP_d$ for $d\leq\ell$.
However, the equations governing Newton hierarchies precisely limit
the possibilities and will lead to a series of non-linear moment
conditions for the fixed coefficients $\a_k(d)$, $k>d>\ell$. These
will be necessary and sufficient for solving the problem.

The projection $\pi$  in (\EE.1) is called {\sl good} if
$\pi(\g)$ is an immersed curve with  normal crossings. 
Such projections are generic.

\Theorem{\EE.4} {\sl Let $\gamma\subset \bbp^2$ be a finite union
of oriented  closed curves of class $C^1$ in the complex projective
plane.  Choose a good projection $\pi:\bbp^2-\bbp^0\to\bbp^1$ from a
point $\bbp^0\notin \g$ and a point $a\in \bbp^1-\pi\g$. Then $\g$
bounds a holomorphic 1-chain $\HC$ in $\bbp^2-\bbp^0$ which is
positive over the point $a$  with sheeting number $\ell$ if
and only if there exists a sequence of homogeneous polynomials
$\{p_d(z_0,z_1)\}_{d=1}^\infty$ with deg$(p_d)=d$, such that 
$\{\TPO_d+p_d\}_{d=1}^\infty$ is a Newton hierarchy of level $\ell$
in a neighborhood of $a$.
} 

\pf
The  necessity of the condition has been established. For the
converse suppose $\{\TPO_d+p_d\}_{d=1}^\infty$ is a Newton hierarchy
of level $\ell$ in a neighborhood $U_0$ of $a$. By analyticity and
Theorem \AA.6 this must hold throughout the connected component
$U$ of $a$ in $\bbp^1-\pi(\g)$. Hence, by
Theorem \FF.3 there is a positive holomorphic chain $\ct$ in
$\pi^{-1}(U)$  (a multisection of degree $\ell$) such that 
$$
C_d\ =\ \pi_*(w^d \ct)\qquad\text{ in }\ \ U.
\tag8.6$$

The remainder of the proof now follows precisely the arguments
given in [HL$_1$, \S 6]. In [HL$_1$] we assumed  
$\TP_d\equiv 0$ in  $U$ = the component of $\infty
\in\bbp^1-\g$. Using boundary regularity,   the ``jump condition''
[HL$_1$, Lemma 5.5] and the Hadamard criteria for rationality
[HL$_1$, Theorem 4.6] we showed that for every connected component
$V$ of $\bbp^1-\g$ there exists a holomorphic 1-chain
$\ct_{V}$ in $V\times \bbc$ such that 
$$
\pi_*\left\{w^d \ct_V \right\}  \ =\ \TP_d\qquad\text{in }
V . 
$$
Furthermore, these chains extend across $\pi^{-1}(\g)$ to define 
a holomorphic 1-chain $\ct$ with $d\ct = \g$. All this was
accomplished using just  one good projection.

The arguments of [HL$_1$, \S 6] are inductive -- from component to
component of $\bbp^1-\g$. To begin one assumed  
$\ct_{U}=0$. However, the arguments apply equally well  if one
merely assumes the existence of a holomorphic chain $\ct$ in
$\pi^{-1}(U)$ which satisfies (\EE.6) in $U$.
\qed

\medskip

An alternative argument for the end of the proof can be given as
follows.  Let $\Delta_\e\subset U$ be the  disk of radius
$\e$ centered at $a$ and let $\ct_\e$ denote the
restriction of $\ct $ to $\pi^{-1}(\Delta_\e)$.
Choose $\e$ so that  $d\ct_\e \equiv \g_{\e}$ is a regular 
(oriented, analytic) curve.  Consider the new ``boundary'' curve
$\Gamma=\g-\g_\e$.  Then for each $d$ we have by (\EE.6) that 
$$
\widetilde{\TP}_d\ \equdef\     \TP_d - \pi_*\{w^d \ct_\e\} \
\equiv\ 0\qquad\text{ in }  \Delta_\e,
\tag{\EE.7}$$
and using (\EE.2) we see that
$$
d\widetilde{\TP}_d\ =\   \pi_*\{w^d \Gamma\} 
\tag{\EE.8}
$$

We now pass to the affine coordinate chart $\zeta = z_0/z_1$ on
$\bbp^1$ and the canonical trivialization given by $z_1$ over this
chart.  In this chart the solution to $\dbar \phi
=\pi_*(w^d\Gamma)$ which vanishes to order $d$ at infinity is given
explicitly by the Cauchy kernel, i.e.,  
$$
\widetilde{\TP}_d(\zeta) \ =\ \frac{1}{2\pi i}
\int_{\Gamma} \frac{w^d}{\eta-\zeta}\, d\eta\ =\ 
\sum_{k=0}^{\infty}\left\{ \frac{1}{2\pi i}
\int_{\Gamma} {w^d\eta^k}\, d\eta \right\}\zeta^{-k-1} \ \equiv\ 0
$$
in a neighborhood of infinity in the $\zeta$-plane. It follows that 
$$
\int_{\Gamma} {w^d\eta^k}\, d\eta\ =\ 0 
\qquad\text{for all } k,d\geq0.
$$
This restricted moment condition is the only hypothesis made in
[HL$_1$, \S 6] and so those arguments apply directly to prove
the existence of  a holomorphic 1-chain $\widetilde{\ct}$ with
compact support in $\bbc^2$ and $d\widetilde{\ct}  = \Gamma$. 
Hence,  $\ct \equiv\widetilde{\ct} +\ct_\e$ is a
holomorphic 1-chain in $\bbp^2-\Gamma$ with  $d\ct =
\g-\g_\e+\g_\e=\g$. Since $\ct_\e$ extends
holomorphically across $\gamma_{\e}$, it follows easily that
$\ct$ is a is a holomorphic chain in $\bbp^2-\gamma$ as desired.

By applying \AA.6 we can restate Theorem \EE.4 as follows.
Let $\ch(d)\equiv H^0(\bbp^1, \co{\bbp^1}(d))$ denote the space of
homogeneous polynomials of degree $d$ in two variables.

\Theorem{\EE.5} {\sl Let $\gamma$, $\bbp^0$, $\pi$ and $a$ be as in
Theorem \EE.4. Then $\g$ bounds a holomorphic 1-chain $\HC$ in
$\bbp^2-\bbp^0$ which is positive over  $a$  with sheeting
number $\ell$ if and only if there exist  homogeneous
polynomials $p_d\in \ch(d)$ for $d=1,...,\ell$, such
that  
$$
\TPO_k\ \equiv \ Q_{k,\ell}(\TPO_1+p_1,...,\TPO_\ell+p_\ell)
\mod  \ch(k)
\tag{\EE.9}
$$
 in a neighborhood of $a$ for all $k>\ell$. } 

\medskip

The non-linear moment conditions (\EE.9) 
represent an explicit sequence of polynomial relations
among the coefficients $\{\a_k(d)\}_{k=d+1}^\infty$ 
discussed in \EE.3 above.

 \Note{\EE.6} Theorems \EE.4 and \EE.5 continue to  hold if one
allows integer multiplicities on the connected components of $\g$.
The  proofs are essentially the same.  \medskip

  \vskip .3in

%
%

%

\noindent
{\bf \S \HH.  Curves in $\bbp^n$.}  \ \   
Let $\gamma\subset \bbp^n$ be an embedded finite union of oriented 
closed curves of class $C^1$.  Fix a codimension-2 linear subspace  
$\bbp^{n-2}$ with $\gamma\cap \bbp^{n-2}=\emptyset$ and consider the
projection 
$$
\Pi: \bbp^n-\bbp^{n-2} \ @>{}>>\ \bbp^1
\tag{\HH.1}
$$
defined in \S \DD. This projection is called {\sl good} if
$\Pi\bigl|_{\gamma}$ is an immersion with normal crossings.
Recall from \S 7 that $\Pi$ is isomorphic the vector bundle given as
the       $(n-1)$-fold direct sum 
$$
\Pi:\co{\bbp^1}(1)\oplus\dots\oplus\co{\bbp^1}(1) \ \arr\ \bbp^1,
$$
Let $\bw =(w_1,...,w_{n-1})$ denote the tautological cross section
of $\Pi^*\{\co{\bbp^1}(1)\oplus\dots\oplus\co{\bbp^1}(1)\}$, where $w_j$
is the  tautological section of the $j\th$ factor.  

For each $\a=(\a_1,...,\a_{n-1}) \geq (0,...,0)$ in $\bbz^{n-1}$ consider
the equation of $\co{\bbp^1}(d)$-valued currents on $\bbp^1$:
$$
\dbar \TP_\a\ =\ \Pi(w^\a \gamma)^{0,1}.
\tag{\HH.2}
$$
by \S7 this equation has solutions unique up to holomorphic sections of 
$\co{\bbp^1}(d)$.

Fix a point $a\in \bbp^1-\Pi(\g)$ and let  $\TPO_\a$ be the unique
solution of (\HH.2) which vanishes to order $d$ at $a$.  Then the 
general solution of (\HH.2) is of the form
$$
\TP_\a\ =\ \TPO_\a\ + p_\a
\tag{\HH.3}
$$
where $p_\a \in H^0(\bbp^1, \co{}(d))$ is a homogeneous polynomial of
degree $d$ in homogeneous coordinates on $\bbp^1$. We recall from \S \FF
\ that this set of moments $\{\TP_\a\}_{\a\geq 0}$ corresponds to a 
homogeneous polynomial family of 
distributional sections of $\co{\bbp^1}(d)$
$$
\TPF d {\la}\ =\ \sum_{|\a|=d} {\ssize \binom d {\a}}\TP_\a \la^\a,
\tag{\HH.4}
$$
which satisfy
the equation 
$$
\dbar \TPF d {\la}\ =\ \Pi_*\{(\la\cdot \bw)^d\,\gamma\}^{0,1}
\tag{\HH.5}
$$
for
$$
\la \cdot \bw \ \equiv \ \sum_j \la_j w_j \ =\ \pi_{\la}^*(w)
$$
where
$
\pi_{\la} : \co{\bbp^1}(1)^{\oplus {(n-1)}}\ \arr\
\co{\bbp^1}(1)  
$
is the $\la$-projection and  $w$ is the tautological cross-section of
$\pi^*\co{\bbp^1}(1)$ over $\co{\bbp^1}(1)$. 
Recall from Definition \FF.6 that $\{\TP_\a\}_{\a\geq 0}$  is called a
Newton family of level $\ell$ 
(with coefficients in $\co{\bbp^1}(1)^{\oplus {(n-1)}}$)
if $\{\TPF d {\la}\}_{d\geq 0}$ is a Newton
heirarchy of level $\ell$ 
(with coefficients in $\co{\bbp^1}(1)$) for all $\la\in \bbc^q$.

\Theorem{\HH.1} {\sl Let $\gamma\subset \bbp^n$ be an embedded finite
union of oriented  closed curves of class $C^1$
with possible integer multiplicities on each component.  Choose a good
projection $\Pi:\bbp^n-\bbp^{n-2}\arr \bbp^1$  from a codimension-2
linear subspace $\bbp^{n-2}\subset\bbp^n-\g$ and fix a point $a\in
\bbp^1-\Pi(\g)$.  Then $\g$ bounds a holomorphic 1-chain in
$\bbp^n-\bbp^{n-2}$ which is positive over $a$ with sheeting number
$\ell$ if and only if there exist homogeneous polynomials
$\{p_{\a}\}_{\a>0}$ with $p_{\a}\in H^0(\bbp^1,\co{}(|\a|))$, such that 
the family $\{\TP_{\a}\}_{\a\geq0}$ given in (\HH.2-3) is a Newton
family of level $\ell$ in a neighborhood of $a$.
}

\pf If $\g=d\HC$ where $\HC$ is a holomorphic 1-chain in
$\bbp^n-\bbp^{n-2}$ which is positive with degree $\ell$  over the point
$a$, then by Theorem 5.6, the functions $\TPF d {\lambda}=\Pi_*((\la\cdot
\bw)^d \HC)$ form a Newton family of degree $\ell$ in a neighborhood
of $a$. Furthermore, we have $\dbar \TPF d{\lambda}=\{d
\TPF d{\lambda}\}^{0,1}= \{\Pi_*((\la\cdot \bw)^d d\HC)\}^{0,1} =
\{\Pi_*((\la\cdot \bw)^d \g)\}^{0,1}$. Expanding as in (\HH.4) gives the
desired $\{\TP_{\a}\}_{\a\geq0}$.

Conversely, suppose that there exist polynomials
$\{p_{\a}\}_{\a\geq 0}$ such that the solutions (\HH.3) form a Newton
family of level $\ell$ over a neighborhood of $a$. Then by Theorem \EE.3
for almost every $\la$ the projection $\pi_\la(\g)$ bounds a holomorphic
chain  $\HC(\la)$ in $\bbp^2-\bbp^0$ with $d\HC(\la)=\pi_\la(\g)$.  The
holomorphic chain $\HC$ is then easily constructed from these projections
as in  [HL$_1$, \S 7]. \qed
\vskip .3in

An alternative argument for the second half of the proof can be given as
follows.  By Theorem \FF.7 the existence of the Newton family
$\{\TP_\a\}_{\a\geq0}$ implies the existence of a multi-section
$\HC_{\Delta}$ of degree $\ell$ over a disk-neighborhood $\Delta$ of $a$
in $\bbp^1-\Pi(\g)$ such that $\Pi_*\{w^\a \HC_{\Delta}\} = C^\a$ in
$\Delta$.  Now by shrinking $\Delta$ slightly we may assume that
$d\HC_{\Delta}$ is a finite set of regular oriented closed curves with
positive integer multiplicities. One now verifies, as in \S \EE \  above, 
that the curve $\g -d\HC_{\Delta}$ satisfies the moment condition of
[HL$_1$] in $\bbc^{n+1} = \bbp^{n+1}-\overline{\Pi^{-1}(a)}$.
Hence, by [HL$_1$, \S 7] there exists a holomorphic 1-chain $\HC_0$
with compact support in $\bbc^{n+1}$ such that $d\HC_0= \g
-d\HC_{\Delta}$.  Thus, $d\HC=\g$ where $\HC\equiv \HC_0+\HC_{\Delta}$
is easily seen to be a holomorphic chain in $\bbp^{n+1}-\g$ (since
$\HC_{\Delta}$ has an analytic continuation across $d\HC_{\Delta}$).

  \vskip .3in

%
%

%

\noindent
{\bf \S \II.  Boundaries of Higher Dimension.}  \ \   
Let $\bd\subset \bbp^n$ be a compact oriented embedded  
$(2p-1)$-dimensional submanifold which is maximally complex. 
Fix a linear subspace   $\bbp^{n-p-1}$ of complex codimension $p+1$ such
that  $\bd\cap \bbp^{n-p-1}=\emptyset$, and consider the projection 
$$
\Pi: \bbp^n-\bbp^{n-p-1} \ @>{}>>\ \bbp^p
\tag{\II.1}
$$
defined in \S \DD. This projection is called {\bf good} if
$\Pi\bigl|_{\bd}$ is an immersion with normal crossings
outside a $C^1$ stratified subset ${\Cal S}=\bigcup_{\mu\geq1}       
{\Cal S}_{\mu}$ where codim$_{\bbr}({\Cal S}_{\mu})=2\mu^2$.
If $\bd$ is of class $C^{2p-1}$, then by Theorem A.4 in [HL$_1$]
a generic projection is good.

As in \S \HH \  we note that $\Pi$ is isomorphic to
$
\Pi:\co{\bbp^p}(1)\oplus\dots\oplus\co{\bbp^p}(1) \ \arr\ \bbp^p
$
and let $\bw =(w_1,...,w_{n-1})$ denote the tautological cross section of
$\Pi^*\{\co{\bbp^p}(1)\oplus\dots\oplus\co{\bbp^p}(1)\}$. For
$\a=(\a_1,...,\a_{n-p}) \geq (0,...,0)$ in $\bbz^{n-p}$ consider
the equation 
$$
\dbar \TP_\a\ =\ \Pi(w^\a \bd)^{0,1}.
\tag{\II.2}
$$
of $\co{\bbp^p}(d)$-valued currents on $\bbp^p$. For $p>1$ the maximal
complexity of $\bd$ implies that 
$$
\dbar \Pi(w^\a \bd)^{0,1}\ =\ 0
$$
(cf. [HL$_1$, \S 3]).
It follows that equation (\II.2) has solutions unique up to holomorphic
sections of  $\co{\bbp^p}(d)$.

We now fix a point  $a\in \bbp^p-\Pi(\bd)$ and let  $\TPO_\a$
be the unique solution of (\II.2) which vanishes to order $d$ at $a$. 
The  general solution of (\II.2) is then of the form
$$
\TP_\a\ =\ \TPO_\a\ + p_\a
\tag{\II.3}
$$
for $p_\a \in H^0(\bbp^p, \co{}(d))$.

\Theorem{\II.1} {\sl Let $\bd\subset \bbp^n$ be as above and suppose 
 $\Pi:\bbp^n-\bbp^{n-p-1}\arr \bbp^p$ is a good projection  from a
 linear subspace $\bbp^{n-p-1}\subset\bbp^n-\bd$. Fix a point
$a\in \bbp^1-\Pi(\bd)$.  Then $\bd$ bounds a holomorphic p-chain in
$\bbp^n-\bbp^{n-2}$ which is positive over $a$ with sheeting number
$\ell$ if and only if there exist homogeneous polynomials
$\{p_{\a}\}_{\a>0}$ with $p_{\a}\in H^0(\bbp^1,\co{}(|\a|))$, such that 
the family $\{\TP_{\a}\}_{\a\geq0}$ given in (\II.2-3) is a Newton
family of level $\ell$ in a neighborhood of $a$.
}
\def\LLL{L}
\pf
Necessity is proved exactly as in  the proof of \HH.1.
The arguments for sufficiency in the proof of  Theorems \EE.1 and  \HH.1
(using results in [HL$_1$, \S 5,6])
establish the existence of a holomorphic $p$-chain $T$ with $dT=\bd$ in
$\bbp^n-\bd-\Pi^{-1}(\Pi{\Cal S})$. 

Furthermore, we have the following.  let $U$ be the neighborhood of $a$ in
$\bbp^p$ where $\{\TP_{\a}\}_{\a\geq0}$ is a Newton family of level
$\ell$, and suppose  $\LLL\subset \bbp^p$ is a projective line which meets
$U$ and for which the linear subspace $\bbp^{n-p+1}_{\LLL} \equiv
\overline{\Pi^{-1}(\LLL)}$ meets $\bd$ transversely. The restriction of
$\{\TP_{\a}\}_{\a\geq0}$  to $U\cap \LLL$ is still a Newton
family of level $\ell$. Furthermore, slicing by $\LLL$ commutes with
$\dbar$.   Hence by Theorem \HH.1 the curve 
$\bd_{\LLL}=\bd\cap \bbp^{n-p+1}_{\LLL}$ bounds a unique holomorphic 
chain $T_{\LLL}$ with compact support in $\bbp^{n-p+1}_{\LLL}
-\bbp^{n-p-1}$ which is positive and $\ell$-sheeted over $a$. 
Moreover, 
$$
T_L\ =\ T \cap \bbp^{n-p+1}_{\LLL} \text{\ \  outside } \ \ 
\Pi^{-1}(\Pi{\Cal S}) 
$$
since they give rise to the same Newton family in $U\cap \LLL$.
From this one concludes that {\bf $T$ has finite mass and bounded support
in $\bbp^n-\bbp^{n-p-1}$}.
Consequently in $\bbp^n-\bbp^{n-p-1}$  we have
$$
dT=\bd +R 
\qquad\text{with}\qquad 
\supp (R)\subset\subset \Pi^{-1}(\Pi{\Cal S}).
$$

Now for each $\la\in \bbc^{n-p}$  let $\pi_{\la}$ and $\pi$ denote the
projections defined in \FF.5 with $\Pi=\pi\circ\pi_\la$. Then 
${\pi_{\la}}_*T$ is a holomorphic $p$-chain with $d {\pi_{\la}}_*T =
{\pi_{\la}}_*\bd + R_\la$
where $R_\la = {\pi_{\la}}_*R$ satisfies
$$
\supp (R_{\la})\subset \pi_{\la}(\Pi^{-1}{\Pi(\Cal S}))=
\coprod_{\mu\geq1}  \pi^{-1}(\Pi{\Cal S}_{\mu}).
\tag{\II.4}
$$
The current $R_{\la}$ is  $d$-closed, integrally flat and of dimension
$2p-1$. Since $\pi^{-1}(\Pi{\Cal S})$ is a $C^1$-stratified set of
dimension $2p-1$, arguments of [F, 4.1.15]  show that $R_\la = c
\pi^{-1}(\Pi{\Cal S})$ for an integer-valued
function $c$ which is constant on components of the top strata of 
$\pi^{-1}(\Pi{\Cal S})$.  Unless $c=0$ this current has unbounded support
in $\bbp^{p+1}-\bbp^{p-1} = \co{\bbp^p}(1)$ because it is invariant 
under dilations in this bundle. This contradicts the boundedness of 
$T$ in $\bbp^n-\bbp^{n-p-1}$.  Thus $R=0$ and the proof is complete.
\qed

  \vskip .3in

%
%

%

\noindent
{\bf \S \JJ.  The General Result in Terms of Moment Conditions.}  \ \   
The characterization given in Theorems  \EE.1, \HH.1 and \II.1
can be reexpressed as a countable family of non-linear moment conditions.
They involve the spaces $\ch(d)\equiv H^0(\bbp^p,\co{}(d))$, and the
universal polynomials $Q_{\a,\ell}$ given in Proposition \FF.8.

\Theorem{\JJ.1} {\sl Let $\bd\subset \bbp^n$ be a compact oriented
submanifold of dimension $2p-1$ which is maximally complex.
Suppose   $\Pi:\bbp^n-\bbp^{n-p-1}\arr \bbp^p$ is a good projection  from a
 linear subspace $\bbp^{n-p-1}\subset\bbp^n-\bd$. Fix a point
$a\in \bbp^1-\Pi(\bd)$.  Then $\bd$ bounds a holomorphic p-chain in
$\bbp^n-\bbp^{n-2}$ which is positive over $a$ with sheeting number
$\ell$ if and only if there exist homogeneous
polynomials $p_{\a}\in \ch(|\a|)$ for  ${0<|\a|\leq \ell}$, such that
the canonical solutions $\TPO_{\a}$ of the equation (\II.2) satisfy
$$
\TPO_{\a}\ \equiv\ Q_{\a,\ell}(C)   \mod \ch(|\a|)
\qquad\text{ for all } \ \ |\a|>\ell
\tag{\JJ.1}
$$
in a neighbornood of $a$, where $\TP_{\a}=\TPO_{\a}+p_{\a}$ for  
$|\a|\leq \ell$.
 }

\vskip .3in
As discussed in \S \TT, when $p=1$ the equations (\JJ.1) reduce to certain
universal polynomial relations among the classical moments of the curve 
in an affine chart on $\bbp^n$.

  \vskip .3in

%
%

%

\noindent
{\bf \S \KK.  The Relation to Dolbeault-Henkin.}  \ \   
The main results here are related and complementary to those of
P. Dolbeault and G. Henkin [D], [DH$_{1,2}$]. In our notation
these authors {\bf  consider the canonical solutions $\TPO_{\a}$
for $|\a|\leq 1$ as  functions of the projection $\Pi$}. That is, they
allow the projection to vary in a neighborhood $U$ of some fixed $\Pi_0$
in the Grassmannian and consider the $\TPO_{\a}$ as functions on $U$.
Their main result asserts that $\a$ bounds a holomorphic $p$-chain
in $\bbp^n$ iff the $\TPO_{\a}$ can be expressed as a finite linear
combination of functions which satisfy a certain non-linear shock
equation on $U$.

This result is beautiful and quite general in that no positivity 
assumptions are required. On the other hand it is somewhat mysterious
in that the functions which satisfy the shock equation are produced
{\sl deus ex machina}.  There is no direct relationship to the geometry
of the  boundary $\g$.

In our result we fix one projection $\Pi_0$ and give a set
of explicit congruences for the functions $\TPO_{\a}$ (for all $\a$)
which are equivalent to the existence of the holomorphic $p$-chain. For
this we need a positivity condition which
guarantees at most a finite-dimsnional family of
possible solutions.

  \vskip .3in

%

\centerline{\bf References}

\vskip .2in

\ref \key [A] \by H. Alexander \paper Linking  and holomorphic
hulls
  \jour J. Diff. Geom.
 \vol 38    \yr 1993    \pages  151-160  \endref

\smallskip

\ref \key [AW$_1$] \by H. Alexander and J. Wermer \book Several Complex
Variables and Banach Algebras
 \publ Springer-Verlag \publaddr New York   \yr  1998 \endref

 \smallskip

\ref \key [AW$_2$] \by  H. Alexander and J. Wermer \paper Linking numbers
and boundaries of varieties
  \jour Ann. of Math.
 \vol 151  \yr  2000  \pages 125-150\endref

 \smallskip

\ref \key [B] \by S. Bochner \paper Analytic and meromorphic continuation
by means of Green's formula
  \jour Ann. of Math.
 \vol 44  \yr  1943 \pages 652-673\endref

 \smallskip

\ref\key [deR]\by  G. de Rham \book Vari\'et\'es Diff\'erentiables, formes,
courants, formes harmoniques\publ 
 Hermann\publaddr Paris \yr 1955\endref

\smallskip

\ref \key [D] \by P. Dolbeault \paper On holomorphic chains with given 
boundary in $\bbc\bbp^n$
 \jour Springer  Lecture Notes, no. 1089,
   \yr 1983   \pages 1135-1140\endref

 \smallskip

\ref \key [DH$_1$] \by P. Dolbeault and G. Henkin \paper
Surfaces de Riemann de bord donn\'e dans $\bbc\bbp^n$,
\inbook Contributions to complex analysis and analytic geometry
 \jour Aspects of Math.
\publ  Vieweg \vol 26 \yr 1994   \pages 163-187 \endref

\smallskip

\ref \key [DH$_2$] \by P. Dolbeault and G. Henkin \paper
 Cha\^ines holomorphes de bord donn\'e dans $\bbc\bbp^n$
 \jour Bull.  Soc. Math. de France \vol  125 \yr 1997   \pages 383-445
\endref

 \smallskip

\ref \key [DL] \by T.-C. Dinh and M. Lawrence \paper
 Polynomial hulls and positive currents
 \jour Ann. Fac. Sci de Toulouse \vol  12 \yr 2003   \pages 317-334
\endref

 \smallskip

\ref \key [Fa$_1$] \by B. Fabr\'e \paper
Sur l'intersection d'une surface de Riemann avec des hypersurfaces  
alg\'ebriques
  \jour C. R. Acad. Sci. Paris
 \vol  322 S\'erie I   \yr 1996    \pages 371-376   \endref

\smallskip

\ref \key [Fa$_2$] \by B. Fabr\'e \paper
On the support on complete intersection 0-cycles
  \jour The Journal of Geometric Analysis
 \vol  12    \yr 2002    \pages 601-614   \endref

\smallskip

\ref\key [F]\by   H. Federer\book Geometric Measure 
Theory\publ 
 Springer--Verlag\publaddr New York \yr 1969\endref

 \smallskip

\ref\key [Fu]\by  W. Fulton \book Intersection Theory\publ 
 Springer--Verlag\publaddr New York \yr 1984\endref

\smallskip

\ref \key [H] \by F.R. Harvey  \paper
Holomorphic chains and their boundaries  \inbook Several Complex
Variables, Proc. of Symposia in Pure Mathematics XXX Part 1
\publ A.M.S.\publaddr Providence, RI  \yr 1977    \pages 309-382  \endref

%

\smallskip

\ref \key [HL$_1$] \by F.R. Harvey and H.B. Lawson, Jr. \paper
On boundaries of complex analytic varieties, I  \jour Ann. of Math.
 \vol 102 \yr 1975  \pages 223-290\endref

\smallskip

\ref \key [HL$_2$] \by F.R. Harvey and H.B. Lawson, Jr. \paper
On boundaries of complex analytic varieties, II  \jour Ann. of Math.
 \vol 106 \yr 1977  \pages 213-238\endref

\smallskip

\ref \key [HL$_3$] \by F.R. Harvey and H.B. Lawson, Jr. \paper
Projective hulls and the projective Gelfand transformation  \jour (to
appear).\endref

\smallskip

\ref \key [HL$_4$] \by F.R. Harvey and H.B. Lawson, Jr. \paper
On the Alexander-Wermer Theorem \jour (to
appear).\endref

 \smallskip

%

\ref \key [HS] \by F.R. Harvey and B. Shiffman \paper
A characterization of holomorphic chains  \jour Ann. of Math.
 \vol 99 \yr 1974  \pages 553-587\endref

\smallskip

\ref\key [Ho]\by   L. H\"ormander \book An Introduction to Complex
Analysis in Several  Variables
\publ Van Nostrand
 \publaddr Princeton, N. J.  \yr 1966 \endref

 \smallskip

\ref \key [L] \by M. Lawrence \paper
Polynomial hulls of rectifiable curves  \jour Amer. J. Math
\vol 117 \yr 1995 \pages 405-417\endref

\smallskip

\ref \key [S] \by G. Stolzenberg \paper uniform approximation on smooth
curves
 \jour Acta math.
 \vol 115   \yr1966    \pages 185-198 \endref

 \smallskip

\ref \key [W$_1$] \by J. Wermer \paper The hull of a curve in $\bbc^n$
 \jour Ann. of Math.
 \vol 68   \yr  1958  \pages 550-561 \endref

 \smallskip

\ref \key [W$_2$] \by J. Wermer \paper The argument principle and
boundaries of analytic varieties
 \jour Operator Theory: Advances and  Applications
   \vol 127   \yr  2001  \pages 639-659 \endref

 \smallskip

 \enddocument